\documentclass[11pt]{article}

\usepackage[margin=1in]{geometry}
\usepackage[T1]{fontenc}
\usepackage[utf8]{inputenc}
\usepackage{lmodern} 
\usepackage{cite}

\usepackage{amsmath,amssymb,amsfonts,amsthm,mathtools}
\usepackage{dsfont,bm}

\numberwithin{equation}{section}

\newtheorem{theorem}{Theorem}[section]
\newtheorem{lemma}[theorem]{Lemma}

\newtheorem{corollary}[theorem]{Corollary}
\newtheorem{assumption}[theorem]{Assumption}
\theoremstyle{remark}
\newtheorem{remark}{Remark}

\usepackage{graphicx}
\usepackage{float}
\usepackage{subcaption}
\usepackage{epstopdf}
\DeclareGraphicsExtensions{.pdf,.png,.jpg,.eps}

\usepackage{algorithm}
\usepackage{algorithmic}

\usepackage{cite} 
\usepackage{hyperref}
\usepackage[capitalize,noabbrev]{cleveref}

\crefname{hypothesis}{Hypothesis}{Hypotheses}
\crefname{fact}{Fact}{Facts}

\usepackage{soul}      
\usepackage{enumerate}
\usepackage{xcolor}

\allowdisplaybreaks

\makeatletter
\renewcommand{\@maketitle}{%
  \begin{center}
    {\Large\bfseries DIGing--SGLD: Decentralized and Scalable Langevin Sampling\\
    over Time--Varying Networks\par}%
    \vskip 1em
    {\Large
      Waheed U.\ Bajwa$^{\dagger}$, 
      Mert G\"urb\"uzbalaban$^{\dagger,\ddagger}$, 
      Mustafa Ali Kutbay$^{\ddagger}$, \\
      Lingjiong Zhu$^{\S}$, 
      Muhammad Zulqarnain$^{\dagger}$\par}
    \vskip 0.5em
    {\small
      $^{\dagger}$Department of Electrical and Computer Engineering, Rutgers University, Piscataway, NJ, USA\\
      $^{\ddagger}$Department of Management Science and Information Systems, Rutgers University, Piscataway, NJ, USA\\
      $^{\S}$Department of Mathematics, Florida State University, Tallahassee, FL, USA\\[0.25em]
    }
  \end{center}
  \vskip 1.5em
  \vspace{-0.5em}
  \renewcommand{\thefootnote}{\fnsymbol{footnote}}%
  \footnotetext[1]{The authors are listed in alphabetical order. The authors can be reached at the following email addresses:
    \texttt{waheed.bajwa@rutgers.edu}, 
    \texttt{mg1366@rutgers.edu}, 
    \texttt{mustafa.kutbay@newark.rutgers.edu}, 
    \texttt{zhu@math.fsu.edu}, 
    \texttt{m.zulqarnain@rutgers.edu}.}
  \renewcommand{\thefootnote}{\arabic{footnote}}%
  \setcounter{footnote}{0}%
  \vskip 1.5em
}
\makeatother

\title{} 
\author{} 

\date{} 

\hypersetup{
  pdftitle  = {DIGing-SGLD: Decentralized and Scalable Langevin Sampling over Time-Varying Networks},
  pdfauthor = {W. U. Bajwa, M. Gurbuzbalaban, M. A. Kutbay, L. Zhu, M. Zulqarnain},
  colorlinks = true, linkcolor=blue, citecolor=blue, urlcolor=blue
}

\providecommand{\headers}[2]{}       
\providecommand{\funding}[1]{}       

\newcommand{\mgb}[1]{{\color{black} #1}} 

\definecolor{purple}{rgb}{0.5,0.0,0.5}
\newcounter{wubctr}

\begin{document}

\maketitle

\begin{abstract}
Sampling from a target distribution induced by training data is central to Bayesian learning, with Stochastic Gradient Langevin Dynamics (SGLD) serving as a key tool for scalable posterior sampling and its decentralized variants enabling learning when data are distributed across a network of agents. This paper introduces DIGing-SGLD, a decentralized SGLD algorithm designed for scalable Bayesian learning in multi-agent systems operating over time-varying networks. Existing decentralized SGLD methods are typically restricted to static network topologies, and many exhibit steady-state sampling bias caused by network effects, even when full batches are used. DIGing-SGLD overcomes these limitations by integrating Langevin-based stochastic sampling with the gradient-tracking mechanism of the DIGing algorithm, originally developed for decentralized optimization over time-varying networks, thereby enabling efficient 
sampling 
without a central coordinator. To the best of our knowledge, we provide the first finite-time, non-asymptotic Wasserstein convergence guarantees with explicit constants for DIGing-SGLD with a constant stepsize over time-varying networks, including a uniform bound over the agents' marginal distributions. Under standard strong convexity and smoothness assumptions, DIGing-SGLD achieves geometric convergence to an $\mathcal{O}(\sqrt{\eta})$ neighborhood of the target distribution, where $\eta$ is the stepsize, with dependence on the target accuracy matching the best-known rates for centralized and static-network SGLD algorithms with constant stepsize. Numerical experiments on Bayesian linear and logistic regressions validate the theoretical results and demonstrate the strong empirical performance of DIGing-SGLD under dynamically evolving network conditions.
\end{abstract}



{
\section{Introduction}
Consider a random vector $X\in\mathbb{R}^d$ with density $\pi(x)$ and, without loss of generality, write $\pi(x)\propto e^{-f(x)}$. The objective is to generate samples from $\pi$ given access to $f$. A canonical instance arises in Bayesian machine learning: given $n$ independent and identically distributed (i.i.d.) observations $Z=\{z_i\}_{i=1}^n$, each with likelihood $p(z_i\mid x)$ and prior $p(x)$, the posterior satisfies
\begin{align}
\label{eqn:data.posterior}
    \pi(x) = p(x\mid Z) \propto p(x)\prod_{i=1}^n p(z_i\mid x) \quad \Longleftrightarrow \quad f(x) = -\log p(x)-\sum_{i=1}^n \log p(z_i\mid x).
\end{align}
Posterior sampling enables principled Bayesian estimation and uncertainty quantification in a wide range of models and tasks, including logistic regression, linear and nonlinear regression, principal component analysis, and neural-network training (see, e.g., \cite{hoff2009first, wang2016towards} for representative applications).

A widely used route to sampling is Markov chain Monte Carlo (MCMC) \cite{Brooks1998JRSSD, Gelman2013}. Among MCMC methods, Langevin-type algorithms discretize the overdamped Langevin diffusion and exploit gradient information about $f$ \cite{Dalalyan}. The Unadjusted Langevin Algorithm (ULA) requires exact gradients $\nabla f(x)$ at every step \cite{DM2017}; when $f$ aggregates many data terms, computing $\nabla f(x)$ entails a full pass over the dataset and becomes computationally expensive. Stochastic Gradient Langevin Dynamics (SGLD) addresses this by replacing the full gradient with unbiased stochastic estimates built from mini-batches and by avoiding Metropolis--Hastings corrections \cite{welling2011bayesian,NEURIPS2018_335cd1b9,raginsky2017non,teh2016consistency}; this leads to scalable sampling procedures that remain effective on large datasets and high-dimensional models. These centralized methods, however, presuppose that data (or stochastic gradients) can be accessed and aggregated at a single location each iteration.

In many modern systems, data and computation are distributed across a network of devices, such as sensors, Internet-of-Things (IoT) platforms, autonomous fleets, and multi-robot systems, where privacy, bandwidth, and energy constraints preclude raw-data aggregation \cite{yang2019survey,Nokleby2020ScalingUp}. In such settings, it is natural to consider the decomposition
\begin{align}
    f(x) = \sum_{i=1}^N f_i(x)
\end{align}
over $N$ \emph{agents}, where the term \emph{agent} generically refers to any computational entity (e.g., a device, node, or processor) capable of performing local computation and communication. Each agent $i$ can compute (stochastic) gradients of its local component $f_i$ but can only exchange limited information with its immediate neighbors on a communication graph. This setting encompasses, for instance, cases in which $n$ data samples are partitioned across $N$ agents, each holding (for simplicity) $n/N$ samples. Because neither the full gradient $\nabla f(x)$ nor a centralized stochastic-gradient estimator is accessible, classical SGLD cannot be applied directly. Instead, one must design sampling dynamics that rely solely on local gradient information and constrained inter-agent communication. While decentralized and federated optimization have been extensively studied in both deterministic and stochastic regimes (see, e.g., \cite{nedic2009distributed,shi2015extra,McMahan2017,aybat2017decentralized,di2016next,arjevani2020ideal,nedic2018network}), decentralized \emph{sampling} introduces additional challenges arising from the interplay between gradient noise, network mixing, and sampling bias.\looseness=-1

Although decentralized variants of ULA and SGLD have been developed in recent years \cite{gurbuzbalaban2021decentralized,PBGG2020,EXTRALangevin,BharL4DC2025}, these approaches are restricted to \emph{static} communication graphs, where the connectivity pattern between agents remains fixed throughout the sampling process. This assumption is limiting, as in many realistic multi-agent systems the communication topology is inherently \emph{time-varying}: links may appear or disappear due to agent mobility, wireless interference, packet drops, or asynchronous operation, and networks may also be reconfigured to enhance privacy, alleviate congestion, or improve robustness against failures and attacks. Consequently, modeling the network as time-varying is more realistic than assuming a fixed topology \cite{TahbazSalehiJadbabaie2008,nedic2014distributed,nedic2017,nedic2018network,scutari2019distributed,fang2025resist}. The variability of connectivity introduces new challenges to maintaining stability and sampling accuracy, and the development of decentralized SGLD algorithms that can operate reliably under such dynamic conditions remains an open problem. To our knowledge, the only work to date that has investigated decentralized Langevin sampling with time-varying connectivity is \cite{kolesov2021decentralized}. 
However, that work builds on the perturbed push-sum framework of \cite{nedic2014distributed} and uses diminishing stepsizes to control the effects of gradient and consensus perturbations. Consequently, the $2$-Wasserstein distance bound in \cite[Theorem~4.7]{kolesov2021decentralized} decays polynomially with the number of iterations. In contrast, the gradient-tracking mechanism introduced below controls network-induced disagreement while allowing a constant stepsize; the resulting bound consists of geometrically decaying transient terms and a stepsize-dependent residual.


In this work, we propose \emph{DIGing Stochastic Gradient Langevin Dynamics} (DIGing-SGLD), a decentralized sampling algorithm designed to operate over \emph{undirected time-varying} networks while relying only on local stochastic gradients and neighbor-to-neighbor communication. DIGing-SGLD integrates stochastic gradient Langevin dynamics with the distributed inexact gradient-tracking mechanism originally introduced in the DIGing algorithm for decentralized optimization over time-varying graphs \cite{nedic2017}. Each agent maintains an auxiliary variable that tracks the evolving network-wide average of the local stochastic gradients, thereby controlling disagreement among the agents' gradient trackers. This mechanism enables agents to collaboratively approximate the global gradient and perform sampling without a central coordinator, extending decentralized Bayesian inference to dynamically changing network topologies without sacrificing convergence speed in the sense that the total number of iterations required to achieve a target accuracy $\epsilon$ matches the $\mathcal{O}\left(\frac{\log(1/\epsilon)}{\epsilon^2}\right)$ convergence rates established for SGLD-type algorithms on static graphs \cite{EXTRALangevin} and in centralized settings \cite[Theorem 4]{dalalyan2019user}.

\subsection{Relation to Prior Work}
On static graphs, \emph{decentralized SGLD} (DE-SGLD) extends SGLD by interleaving local stochastic-gradient updates with consensus steps among neighboring agents \cite{gurbuzbalaban2021decentralized,PBGG2020,cadena2021stochastic}. Convergence guarantees for DE-SGLD have been established under strong convexity and smoothness assumptions for constant stepsizes, and in certain nonconvex regimes with decaying stepsizes \cite{PBGG2020,cadena2021stochastic}. Momentum-based and event-triggered variants have also been proposed to enhance convergence behavior and reduce communication costs \cite{gurbuzbalaban2021decentralized,bhar2023distributed}. Within these approaches, using decaying stepsizes for strongly convex problems can lead to slow convergence, motivating constant-stepsize schemes in the literature. However, with constant stepsizes a key limitation is the emergence of a steady-state bias—even in the full-batch limit—arising from network-induced discrepancies among agents’ local gradients.
EXTRA-based decentralized Langevin methods~\cite{EXTRALangevin}, including \emph{EXTRA-SGLD} and \emph{generalized EXTRA Langevin dynamics}, mitigate this issue on \emph{static} graphs by incorporating bias-correction techniques from deterministic decentralized optimization, specifically the EXTRA algorithm and its generalizations \cite{shi2015extra,jakovetic2018unification}. The recent work \cite{BharL4DC2025} also employs a gradient-tracking mechanism similar in spirit to EXTRA but focuses on decentralized ULA with deterministic gradients, without addressing the stochastic-gradient setting. In particular, Bhar et al.~\cite{BharL4DC2025} introduced GT-DULA for fixed connected undirected networks with exact local gradients. Another line of research considers \emph{Metropolis-adjusted decentralized Hamiltonian Monte Carlo}, which achieves asymptotically exact sampling when local gradients are deterministic \cite{kungurtsev2023decentralized}.


\ In contrast to these static-network methods, we consider decentralized Langevin-based sampling using stochastic-gradient oracles over \emph{time-varying} networks. In this setting, stochastic-gradient error and time-varying network mixing enter the gradient-tracking and sampling recursions simultaneously, preventing a direct application of existing fixed-network gradient-tracking analyses. \ Although decentralized optimization over time-varying networks has been studied extensively (see, e.g., \cite{TahbazSalehiJadbabaie2008,nedic2014distributed,nedic2017,nedic2018network,scutari2019distributed,yang2019survey}), this literature primarily addresses optimization or consensus problems rather than sampling from a target distribution. The closest related work on decentralized Langevin sampling over time-varying networks is \cite{kolesov2021decentralized}, which analyzes stochastic-gradient Langevin dynamics on directed time-varying graphs and establishes convergence guarantees for the network average. However, the rates in \cite{kolesov2021decentralized} are based on the decentralized optimization framework of \cite{nedic2014distributed}, which is known to exhibit network-induced bias and slow convergence. 

More precisely, \cite{kolesov2021decentralized} studies decentralized
Langevin methods over directed, time-varying communication graphs.
Assuming that each local objective $f_i$ is $\mu$-strongly convex and
$L$-smooth and that the stochastic-gradient variance is bounded, they use the diminishing stepsize
$\alpha_k=\frac{c_0}{1+k}$, where
$c_0=\min\left\{\frac{1}{2L},\,\frac{\mu}{4L^2}\right\}$,
to establish the following bound for the $2$-Wasserstein distance
$\mathcal{W}_2$ between the law $\bar{\nu}_K$ of the network-average iterate
and the target distribution $\pi$:
$
\mathcal{W}_2(\bar{\nu}_K,\pi)
=
\mathcal{O}\!\left(
\frac{1}{(\mu c_0-c)K^c}
\right)
$
{for $0<c<\mu c_0$.}
{When the condition number $\kappa:=L/\mu$ is large,}
$
\mu c_0
=
\Theta\!\left(\frac{1}{\kappa^2}\right),
$
and hence the admissible exponent satisfies
$c=\mathcal{O}(1/\kappa^2)$.
{Consequently, for the guarantee
$\mathcal{W}_2(\bar{\nu}_K,\pi)\leq\epsilon$,
this upper bound yields the iteration complexity
$
K
=
\mathcal{O}\!\left(\epsilon^{-1/c}\right)
$,
whose exponent in $1/\epsilon$ scales as
$\Theta(\kappa^2)$ when $c=\Theta(1/\kappa^2)$.}
{Thus, for ill-conditioned problems, this bound can be highly conservative,
even when the stochastic-gradient error is set to zero.}
{This motivates decentralized samplers that use stochastic mini-batch gradients
with constant stepsizes over time-varying networks and admit sharper finite-time
guarantees.}
Our finite-time Wasserstein analysis of constant-stepsize DIGing-SGLD jointly controls stochastic-gradient error, gradient-tracking error, time-varying consensus dynamics, and Langevin discretization error.

Before presenting our main contributions, we note that several other studies consider related but distinct problems. Some works focus on decentralized maximum-likelihood estimation rather than full Bayesian inference \cite{barbarossa2007decentralized,4407653,blatt2004distributed}, while others address decentralized Bayesian inference by requiring agents to share local posterior distributions \cite{campbell2014approximate}, an approach that is typically communication-intensive and computationally costly. We also note that distributed sampling methods such as Consensus Monte Carlo \cite{Rendell03042021} and other ADMM-style distributed MCMC algorithms \cite{vono2022efficient,vono-distributed}, as well as parallel data-partitioning approaches \cite{dai2012sampling,chowdhury2018parallel,barbos2017clone}, all rely on a central coordinator to aggregate information. Similarly, MCMC and Langevin algorithms developed for federated learning settings depend on a central server to orchestrate parameter updates \cite{el2021federated,liang2025bayesian,sun2024federated}. In contrast, our setting considers an \emph{ad hoc} network without any central coordinator capable of aggregating information. Consequently, these methods are not directly applicable to the decentralized, coordinator-free environment studied in this paper.

\subsection{Our Contributions}
We develop and analyze DIGing-SGLD, a decentralized sampling algorithm that extends SGLD to undirected \emph{time-varying} networks. DIGing-SGLD integrates stochastic gradient Langevin dynamics with distributed inexact gradient tracking~\cite{nedic2017} to correct for network-induced drift and maintain consensus among agents under dynamic connectivity. From a theoretical standpoint, to the best of our knowledge, we establish the first finite-time, non-asymptotic Wasserstein convergence guarantees with explicit constants for decentralized SGLD with a constant stepsize over time-varying networks. Under standard assumptions that each local function $f_i$ is $\mu$-strongly convex and $L$-smooth, and that the sequence of mixing matrices satisfies a joint spectral condition over bounded time intervals, we show that the marginal distribution of each agent’s iterate converges in the $2$-Wasserstein distance at a \emph{geometric rate} to an $\mathcal{O}(\sqrt{\eta})$ neighborhood of the target distribution, where $\eta$ is the constant stepsize. 

In particular, for any fixed network size $N$, we prove that with an appropriate choice of the stepsize $\eta=\mathcal{O}(\epsilon^2)$, after $K = \mathcal{O}\!\left(\tfrac{\log(1/\epsilon)}{\epsilon^2}\right)$ iterations, every agent can sample from a distribution that lies within $\epsilon$ of the target in $\mathcal{W}_2$ distance. The resulting bounds make explicit the dependence on the stepsize, network connectivity, gradient-noise variance, problem dimension, and number of agents. Notably, the $\epsilon$-dependence of our rates matches the best known results for SGLD-type methods with constant stepsize in both centralized and decentralized static-graph settings \cite{EXTRALangevin,dalalyan2019user}, despite the additional challenges posed by time-varying networks.
More specifically, these guarantees apply to constant-stepsize decentralized Langevin sampling with stochastic gradients and gradient tracking over time-varying networks, addressing an important theoretical gap in scalable decentralized Bayesian inference. 

Beyond theoretical analysis, DIGing-SGLD addresses sampling over time-varying networks by extending the frameworks of DE-SGLD~\cite{gurbuzbalaban2021decentralized} and DIGing-based optimization~\cite{nedic2017}, thereby unifying decentralized stochastic sampling and gradient-tracking-based optimization over time-varying networks within one framework. Experiments on Bayesian linear and logistic regression with both synthetic and real datasets corroborate the theory and demonstrate that, under time-varying network topologies, DIGing-SGLD outperforms DE-SGLD. 

\subsection{Notation} 
We let $1_{m}$ denote the $m$-dimensional all-ones column vector and $I_{m}$ the $m\times m$ identity matrix, for any positive integer $m$. For $v=[v_{1}^{\top},\ldots,v_{N}^{\top}]^{\top}\in\mathbb{R}^{Nd}$ with $v_{i}\in\mathbb{R}^{d}$, its average is defined as $\bar{v}:=\tfrac{1}{N}\sum_{i=1}^{N}v_{i}\in\mathbb{R}^{d}$, and the \emph{replicated/stacked average vector} as $\mathbf{\bar{v}}:=[\bar{v}^{\top},\ldots,\bar{v}^{\top}]^{\top}=\tfrac{1}{N}\big((1_{N}1_{N}^{\top})\otimes I_{d}\big)v\in\mathbb{R}^{Nd}$. The \emph{consensus error} of $v$ is $\tilde{v}:=v-\mathbf{\bar{v}}=\mathcal{L}_{N}v$, where $\mathcal{L}_{N}:=I_{Nd}-\tfrac{1}{N}\big((1_{N}1_{N}^{\top})\otimes I_{d}\big)$ is a symmetric $Nd\times Nd$ matrix. For $a\in\mathbb{R}^{Nd}$, we define $\|a\|_{\mathcal{L}_{N}}:=\sqrt{\langle a,\mathcal{L}_{N}a\rangle}$. For a vector $x$, $\|x\|$ denotes the Euclidean norm, while for a random vector $X$, we write $\|X\|_{L_{2}}:=(\mathbb{E}\|X\|^{2})^{1/2}$.  

We denote by $\mathcal{S}_{\mu,L}(\mathbb{R}^d)$ the class of functions $g:\mathbb{R}^{d}\to\mathbb{R}$ that are $\mu$-strongly convex and $L$-smooth, i.e.,
\begin{align}
\frac{\mu}{2}\|x-y\|^{2}\le g(x)-g(y)-\nabla g(y)^{\top}(x-y)\le \frac{L}{2}\|x-y\|^{2}, \qquad \forall x,y\in\mathbb{R}^{d}.
\end{align}
Let $\mathcal{P}_{2}(\mathbb{R}^{d})$ denote the space of Borel probability measures on $\mathbb{R}^{d}$ with finite second moment. For $\mu_{1},\mu_{2}\in\mathcal{P}_{2}(\mathbb{R}^{d})$, the $2$-Wasserstein distance is defined as $\mathcal{W}_{2}(\mu_{1},\mu_{2}) := \left( \inf \mathbb{E}\!\left[\|Z_{1}-Z_{2}\|^{2}\right] \right)^{1/2}$, where the infimum is taken over all pairs of random variables $(Z_{1},Z_{2})$ defined on a common probability space with marginal distributions $\mu_{1}$ and $\mu_{2}$, respectively; see \cite{villani2008optimal} for further details.
}
\section{Preliminaries and Problem Formulation}\label{sec:background}

\vspace{0.05in}
Langevin algorithms are popular MCMC methods for obtaining samples from a given target density $\pi(x)$ of interest. 
If we consider $f(x) := -\log(\pi(x))$, classic first-order Langevin algorithms are based on discretizing 
the \emph{overdamped Langevin diffusion}:
\begin{equation}\label{eq:overdamped-2}
dX(t)=-\nabla f(X(t))dt+\sqrt{2}dW_{t},
\end{equation}
(see e.g. \cite{Dalalyan,DM2016,DM2017,dalalyan2019user,EH2021,BCESZ2022,Chewi2024}). This diffusion admits $\pi(x)$ as the unique stationary distribution under some smoothness and growth assumptions on $f$ \cite{pavliotis2014stochastic}. 
Here, $W_{t}$ is a standard $d$-dimensional Brownian motion for $t\geq 0$ initialized as $W_0 = 0$. ULA \cite{DM2017,DM2016} is a fundamental Langevin algorithm, based on the Euler-Maruyama discretization of \eqref{eq:overdamped-2}, and results in the dynamics
\begin{equation}\label{discrete:overdamped}
x_{k+1}=x_{k}-\eta \nabla f(x_k)+\sqrt{2\eta}w_{k}\,,
\end{equation}
where $\eta>0$ is the stepsize parameter, and $w_k \in \mathbb{R}^d$ is a sequence of i.i.d. standard Gaussian random vectors $\mathcal{N}(0,I_{d})$. Then, it is known that for strongly convex and smooth functions $f$, the iterates $x_k$ converge to a neighborhood of the target distribution in the 2-Wasserstein distance where the size of the neighborhood goes to zero as $\eta\to 0$ \cite{dalalyan2019user,Dalalyan}. 

The ULA algorithm works with deterministic gradients, which is often expensive to compute for Bayesian inference problems involving large data. This motivated
the development of Langevin algorithms that can support stochastic gradients. In particular, if one replaces the full gradient $\nabla f$ in \eqref{discrete:overdamped} by a stochastic unbiased estimate {$\widetilde \nabla f$} of the gradient, the resulting algorithm
is known as the \emph{stochastic gradient Langevin dynamics} (SGLD) (see, e.g., \cite{welling2011bayesian}) that has found many applications in machine learning and large-scale Bayesian data analysis due to its scalability properties. As an example, when $f$ has the form given in \eqref{eqn:data.posterior} corresponding to $n$ data samples $Z=\{z_i\}_{i=1}^n$ and the number of samples $n$ is large, computing the gradient $\nabla f(x)$ requires processing all data samples and is often computationally expensive. However, if we estimate the gradient using $b$ randomly sampled data points, i.e., using the estimator $\widetilde{\nabla} f(x):= -\nabla \log p(x)-\frac{n}{b}\sum_{\ell=1}^b \nabla \log p(z_{i_\ell}\mid x)$, where the indices $i_{\ell}$ are sampled independently and uniformly (with replacement) from $\{1,2,\dots,n\}$ and $b\ll n$, then $\widetilde{\nabla} f(x)$ is an unbiased stochastic estimate of the full gradient $\nabla f(x)$ and is significantly cheaper to compute.
This principle underlies SGLD, and SGLD 
and its variants admit various performance guarantees in a variety of metrics and under various assumptions on $f$; see e.g. \cite{dalalyan2019user,Dalalyan,raginsky2017non}. However, SGLD is still a centralized Langevin algorithm and is not applicable to decentralized sampling problems, which we discuss next.

\subsection{Decentralized Sampling Over a Network} In the context of decentralized sampling within a network of $N$ agents, the objective is to sample at any one of the agents from a target distribution with density $\pi(x) \propto e^{-f(x)}$ on $\mathbb{R}^d$,
where the potential (or objective) function $f$ admits a decomposition over the network:
$f(x):=\sum_{i=1}^Nf_i(x)$. 
In this context, both the component function $f_i$ and the (stochastic) estimates of the gradient $\nabla f_i$ are only available to agent $i\in\{1,2,\dots,N\}$.

Decentralized stochastic gradient Langevin dynamics (DE-SGLD) is a decentralized version of the SGLD algorithm for multi-agent systems described by undirected static graphs~\cite{gurbuzbalaban2021decentralized,PBGG2020,cadena2021stochastic}. Given an undirected graph $\mathcal{G}=(\mathcal{V},\mathcal{E})$, where $\mathcal{V}=\{1,2,\dots,N\}$ is the set of agents and $\mathcal{E}$ is the set of links/edges between the agents, DE-SGLD updates are as follows. Each agent $i$ has access to a subroutine that returns (stochastic) estimates of the gradient of its local component function $f_i(x)$, which contributes to the global function $f(x)=\sum_{i=1}^N f_i(x)$. Agent $i$ maintains a local iterate $x_i^{(k)}$ at step $k$ and updates it after receiving the local variables of its immediate neighbors $j\in\Omega_i:=\{j:(i,j)\in\mathcal{E}\}$, where we adopt the convention that $i\in\Omega_i$. The update consists of taking a weighted average of $\{x_j^{(k)}\}_{j\in\Omega_i}$ followed by a stochastic gradient step with respect to $f_i(x)$ and the addition of Gaussian noise~\cite{swenson2020distributed,gurbuzbalaban2021decentralized}:
\begin{equation}
x_{i}^{(k+1)}=\sum_{j\in\Omega_{i}}W_{ij}x_{j}^{(k)}-\eta\,\widetilde{\nabla} f_{i}\!\left(x_{i}^{(k)}\right)+\sqrt{2\eta}\,w_{i}^{(k+1)}.
\end{equation}
Here, as before, $\eta>0$ is the stepsize, $w_{i}^{(k)}$ are i.i.d.\ standard Gaussian random vectors with zero mean and identity covariance matrix for every $i$ and $k$, and $\widetilde\nabla f_{i} \left(x_{i}^{(k)}\right)$ is an unbiased estimator of the gradient $\nabla f_i\left(x_i^{(k)}\right)$. 
Additionally, the $W_{ij}$'s correspond to the entries of a doubly stochastic mixing matrix $W=(W_{ij})$ that determines the weights for the averaging and that respects the network structure, i.e. $W_{ij}>0$ only when $(i,j)\in \mathcal{E}$. If the Gaussian noise term is omitted in DE-SGLD, the iterations reduce to the decentralized stochastic gradient descent algorithm~\cite{swenson-journal,robust-network-asg}, which itself builds on the decentralized gradient descent (DGD) method \cite{nedic2009distributed}.

\subsection{Decentralized Sampling Over a Time-Varying Network}
Our focus in this paper is on decentralized sampling in a multi-agent system where the underlying network topology is time-varying. We consider a time-invariant set of agents $\mathcal{V}=\{1,2,\ldots,N\}$ interacting over a sequence of undirected graphs $\{\mathcal{G}(k)=(\mathcal{V},\mathcal{E}(k))\}_{k=0}^{\infty}$. While the agent set $\mathcal{V}$ remains fixed, the edge set $\mathcal{E}(k)$ may change over time. At each iteration $k$, an unordered pair $(i,j)\in\mathcal{E}(k)$ if and only if agents $i$ and $j$ can exchange information at time $k$. Since the graphs are undirected, $(i,j)\in\mathcal{E}(k)$ implies $(j,i)\in\mathcal{E}(k)$. The neighborhood of agent $i$ at time $k$, including the agent itself, is defined as $\Omega_i(k):=\{j : (j,i)\in\mathcal{E}(k)\}$.

A limitation of existing decentralized Langevin algorithms is the lack of theoretical guarantees when the underlying network is time-varying and gradients are stochastic. In contrast, decentralized deterministic optimization methods such as DIGing \cite{nedic2017}, along with extensions that handle stochastic gradients \cite{li2020s}, provide rigorous convergence guarantees for minimizing objectives of the form $f(x)=\sum_{i=1}^N f_i(x)$ over time-varying networks. However, DIGing is fundamentally an optimization algorithm and cannot be used directly for sampling from a target distribution. In this paper, we introduce and analyze DIGing-SGLD, which, as its name suggests, integrates SGLD-type updates with DIGing-based gradient tracking to achieve Langevin sampling over time-varying networks. The next section describes the proposed algorithm in detail and states the key assumptions underlying the theoretical results presented thereafter.

\section{DIGing-SGLD for Time-Varying Networks}

In this section, we propose \textit{DIGing stochastic gradient Langevin dynamics}
(DIGing-SGLD) for time-varying networks, where the edge set $\mathcal{E}(k)$
and the associated mixing matrix $W^{(k)}$ may vary with $k$. Let $x_i^{(k)}$ denote the
local iterate of agent $i$. We make the standard assumption that at iteration $k$, agent $i$ has access to a stochastic gradient oracle that provides
$
\widetilde{\nabla} f_i\!\left(x_i^{(k)},v_i^{(k+1)}\right),
$
an estimate of the gradient $\nabla f_i(x_i^{(k)})$, where $v_i^{(k+1)}$ is the random variable used to generate this stochastic gradient. We use $w_i^{(k+1)}$ to denote the Gaussian noise added to the iterate after the stochastic gradient is computed at time $k$ by agent $i$. We assume the stochastic gradient is
conditionally unbiased and has bounded conditional variance; the precise filtration and
independence assumptions are stated below.

\begin{assumption}[Noise filtration under deterministic network mixing]
\label{assumption:noise}
Assume that $x_i^{(0)}$ has finite second moment, i.e., $\mathbb E[\|x_i^{(0)}\|^2]<\infty$, for every $i=1,2,\dots,N$. We use $\mathcal F_{\mathrm{gen}}(\cdot)$ to denote the $\sigma$-algebra generated by its arguments. For $k\geq0$, define the $\sigma$-algebra $\mathcal H_k$ containing the information available before the stochastic-gradient evaluations at iteration $k$ by
\begin{equation}\label{eq:pre-oracle-filtration}
\mathcal H_k:=\mathcal F_{\mathrm{gen}}\!\left(x_1^{(0)},x_2^{(0)}, \dots, x_N^{(0)},
\{v_i^{(t)}:1\leq i\leq N,\ 1\leq t\leq k\},
\{w_i^{(t)}:1\leq i\leq N,\ 1\leq t\leq k\}\right),
\end{equation}
with the convention that the indexed sets are empty when $k=0$.
Since the mixing matrices $W^{(k)}$ are deterministic, they need not be included in $\mathcal H_k$. At time $k$, agent $i$ evaluates
$\widetilde\nabla f_i(x_i^{(k)},v_i^{(k+1)})$ and we set
\begin{equation}
\xi_i^{(k+1)}:=\widetilde\nabla f_i(x_i^{(k)},v_i^{(k+1)})-\nabla f_i(x_i^{(k)}).
\end{equation}
We assume that there exists $\sigma\geq0$ such that, for every $i$ and $k$,
\begin{equation}\label{eq:conditional-gradient-noise}
\mathbb E[\xi_i^{(k+1)}\mid\mathcal H_k]=0,
\qquad
\mathbb E[\|\xi_i^{(k+1)}\|^2\mid\mathcal H_k]\leq\sigma^2,
\end{equation}
and the vectors $\{\xi_i^{(k+1)}\}_{i=1}^N$ are conditionally independent given $\mathcal H_k$. 
For two $\sigma$-algebras $\mathcal A$ and $\mathcal B$, let $\mathcal A\vee\mathcal B$ denote the $\sigma$-algebra they generate jointly. After the stochastic-gradient evaluations at iteration $k$, define the $\sigma$-algebra $\mathcal G_k$ containing the information available before drawing the Gaussian noise by
\begin{equation}\label{eq:post-oracle-filtration}
\mathcal G_k:=\mathcal H_k\vee\mathcal F_{\mathrm{gen}}(v_1^{(k+1)},\ldots,v_N^{(k+1)}){\color{blue}.}
\end{equation}
Conditional on $\mathcal G_k$, the vectors $\{w_i^{(k+1)}\}_{i=1}^N$ are independent $\mathcal N(0,I_d)$ random vectors, and
$\mathcal H_{k+1}=\mathcal G_k\vee{\color{blue}\mathcal F_{\mathrm{gen}}}(w_1^{(k+1)},\ldots,w_N^{(k+1)})$.
To simplify notation, we suppress the dependence on $v_i^{(k+1)}$ and write $\widetilde\nabla f_i(x_i^{(k)})$ for $\widetilde\nabla f_i(x_i^{(k)},v_i^{(k+1)})$.
\end{assumption}
We are now ready to introduce the iterates of DIGing-SGLD.
In contrast to DE-SGLD, DIGing-SGLD maintains an additional auxiliary local variable $y_i^{(k)}$ at each iteration, which enables each agent to track the global gradient and mitigate the bias introduced by decentralized information exchange and the time-varying nature of the underlying network.

To initialize DIGing-SGLD, each agent $i$ evaluates $g_i^{(0)}=\widetilde\nabla f_i(x_i^{(0)},v_i^{(1)})$ and sets $y_i^{(0)}=g_i^{(0)}$. At iteration $k$, agent $i$ receives both the local iterates $x_j^{(k)}$ and the auxiliary variables $y_j^{(k)}$ from its neighbors $j\in\Omega_i(k)$ and first performs
\begin{align}
 x_i^{(k+1)}&=\sum_{j\in\Omega_i(k)}W_{ij}^{(k)}x_j^{(k)}-\eta y_i^{(k)}+\sqrt{2\eta}\,w_i^{(k+1)},\label{eqn:DIGing_algo_updates}
\end{align}
then makes one fresh oracle call $g_i^{(k+1)}=\widetilde\nabla f_i(x_i^{(k+1)},v_i^{(k+2)})$ and updates
\begin{align}
 y_i^{(k+1)}&=\sum_{j\in\Omega_i(k)}W_{ij}^{(k)}y_j^{(k)}+g_i^{(k+1)}-g_i^{(k)}.\label{eqn:DIGing_algo_updates_2}
\end{align}
Thus each iteration uses one new stochastic-gradient evaluation per agent and stores the preceding evaluation. In the absence of stochastic-gradient noise and Gaussian noise, these updates reduce to DIGing \cite{nedic2017}.
As specified in Assumption~\ref{assumption:noise}, the oracle values used to form $y_i^{(k)}$ are available before the Gaussian vector $w_i^{(k+1)}$ is drawn. The updates use the mixing matrix $W^{(k)}$, which is deterministic and therefore independent of the stochastic noise generated by the algorithm. The structural assumptions on $W^{(k)}$ are given in Assumption~\ref{assump:W}.\looseness=-1

Using the definition of $\xi_i^{(k+1)}$, we can rewrite the DIGing-SGLD updates as follows:
\begin{align}
&x_{i}^{(k+1)}=\sum_{j\in\Omega_{i}(k)}W_{ij}^{(k)}x_{j}^{(k)}-\eta y_{i}^{(k)}+\sqrt{2\eta}w_{i}^{(k+1)},\label{DIGing:x}
\\
&y_{i}^{(k+1)}=\sum_{j\in\Omega_{i}(k)}W_{ij}^{(k)}y_{j}^{(k)}+\nabla f_{i}\left(x_{i}^{(k+1)}\right)-\nabla f_{i}\left(x_{i}^{(k)}\right)
+\xi_{i}^{(k+2)}-\xi_{i}^{(k+1)}.\label{DIGing:y}
\end{align}
By introducing the function $F(x):\mathbb{R}^{Nd} \rightarrow \mathbb{R}$ defined as $F(x):=\sum_{i=1}^N f_i(x_i)$ for any $x:=\left(x_{1}^{\top},\ldots,x_{N}^{\top}\right)^{\top}\in\mathbb{R}^{Nd}$, and by stacking the local variables $x_i^{(k)},y_{i}^{(k)}$ into the vectors $x^{(k)} = \left[\left(x_1^{(k)}\right)^{\top}, \ldots, \left(x_N^{(k)}\right)^{\top}\right]^{\top} \in \mathbb{R}^{Nd}$ and $y^{(k)} = \left[\left(y_1^{(k)}\right)^{\top}, \ldots, \left(y_N^{(k)}\right)^{\top}\right]^{\top} \in \mathbb{R}^{Nd}$, we can also rewrite \eqref{DIGing:x}--\eqref{DIGing:y} in vector form as
\begin{align}
&x^{(k+1)}=\mathcal{W}^{(k)}x^{(k)}-\eta y^{(k)}+\sqrt{2\eta}w^{(k+1)},\label{x:vector:form}
\\
&y^{(k+1)}=\mathcal{W}^{(k)}y^{(k)}+\nabla F\left(x^{(k+1)}\right)-\nabla F\left(x^{(k)}\right)
+\xi^{(k+2)}-\xi^{(k+1)},\label{y:vector:form}
\end{align}
where  
$w^{(k)} := \left[\left(w_1^{(k)}\right)^{\top}, \ldots, \left(w_N^{(k)}\right)^{\top}\right]^{\top}\in\mathbb{R}^{Nd}$, 
$\xi^{(k)} := \left[\left(\xi_1^{(k)}\right)^{\top}, \ldots, \left(\xi_N^{(k)}\right)^{\top}\right]^{\top}\in\mathbb{R}^{Nd}$,
and $\mathcal{W}^{(k)}:=W^{(k)}\otimes I_{d}$ for any $k=0,1,2,\ldots$

Next, we consider the mixing matrices $W^{(k)}$ and introduce the notation
\begin{equation}
W_{B}^{(k)}:=W^{(k)}W^{(k-1)}\cdots W^{(k-B+1)},
\end{equation}
for any $k=0,1,2,\ldots$ and any $B=1,2,\ldots$, with the convention that $W_{B}^{(k)}=I_{N}$ for any $k<0$ and $W_{0}^{(k)}=I_{N}$ for any $k$. The product $W_{B}^{(k)}$ captures the aggregated connectivity of the network over the time interval from $k-B+1$ to $k$. Moreover, we introduce the notation
\begin{equation}
\mathcal{W}_{B}^{(k)}:=W_{B}^{(k)}\otimes I_{d},
\end{equation}
for any $k=0,1,2,\ldots$ and any $B=1,2,\ldots$. 
We make the following assumption on the mixing matrices.
\begin{assumption}\label{assump:W}
For any $k=0,1,2,\ldots$, the mixing matrix $W^{(k)}=\left(W_{ij}^{(k)}\right)\in\mathbb{R}^{N\times N}$ is symmetric \mgb{with non-negative entries} and 
satisfies the following properties:

(i) (Decentralized property) If $i\neq j$ and the edge $(j,i)\notin\mathcal{E}(k)$, then $W_{ij}^{(k)}=0$.

(ii) (Double stochasticity) $W^{(k)}1_{N}=1_{N}$  and $1_{N}^{\top}W^{(k)}=1_{N}^{\top}$.

(iii) (Joint spectral property) There exists a positive integer $B$ such that, for every $k=0,1,2,\ldots$, $\delta:=\sup_{k\geq B-1}\delta(k)\in \mgb{(0,1)}$, where $\delta(k):=\sigma_{\max}\left\{W_{B}^{(k)}-\frac{1}{N}1_{N}1_{N}^{\top}\right\}$.
\end{assumption}

Parts (i) and (ii) of Assumption~\ref{assump:W} are standard and are imposed even when the underlying graph is static; see, e.g., \cite{gurbuzbalaban2021decentralized,EXTRALangevin,yuan2016convergence}. In particular, part (i) ensures that the averaging operation respects the network's connectivity pattern, while part (ii) guarantees that, if all agents converge to the same vector, the mixing process preserves this agreement (consensus). Part (iii) of Assumption~\ref{assump:W} controls the spectral gap $\Delta:= 1- \delta$ uniformly over iterations and ensures that the connectivity observed by the iterates across any $B$ consecutive iterations is sufficient (i.e. $\Delta>0$). This assumption, along with closely related variants, appears frequently in the literature; see e.g. \cite{nedic2017,di2016next} and the references therein. Also, for static graphs, it is standard to assume part (iii) with $B=1$, in which case $\delta = \delta(k)$ for every $k$; see \cite{gurbuzbalaban2021decentralized,EXTRALangevin,yuan2016convergence}. \mgb{Note that $\delta=0$ corresponds to exact averaging over every admissible $B$-step block, i.e.,
$
W_B^{(k)} = N^{-1} 1_N 1_N^\top.
$
In this degenerate case, the problem is not genuinely decentralized in the sense that the block communication acts as exact consensus, effectively collapsing the network disagreement over each such block. We therefore exclude this degenerate exact-block-averaging case from consideration.} 

For analysis purposes, we will also make the following assumption throughout on the local objectives. 
\begin{assumption}
\label{assumption:f}
Each local objective function satisfies \( f_i \in \mathcal{S}_{\mu,L}(\mathbb{R}^d) \) for every \( i = 1, 2, \ldots, N \); that is, each \( f_i(x) \) is \(\mu\)-strongly convex and \(L\)-smooth.
\end{assumption}
Under Assumption~\ref{assumption:f}, the global objective $f(x)=\sum_{i=1}^{N} f_i(x)$ admits a unique minimizer $x_\ast$, and the target distribution $\pi(x) \propto e^{-f(x)}$ is strongly log-concave. This assumption has been employed in the literature to analyze the DIGing algorithm \cite{nedic2017,li2020s} for decentralized optimization. In this work, we adopt the same assumption but pursue a fundamentally different goal: \emph{we focus on sampling rather than optimization}. Assumption \ref{assumption:f} has also been used in the literature to obtain convergence
guarantees for decentralized Langevin algorithms such as DE-SGLD and EXTRA-SGLD \cite{gurbuzbalaban2021decentralized,EXTRALangevin}; however, these works are restricted to static communication graphs. In contrast, our analysis accommodates time-varying network topologies. 

\subsection{Main Results} \label{subsection:main_results}

Before we proceed to the main results of the paper, we introduce additional notation. Let $\alpha$, $\beta$, and $\lambda$ be positive scalars, and assume that $\lambda\in(\delta^{1/B},1)$ for $B\geq 1$, so that $\delta<\lambda^{B}<1$. For each iteration $k$, recall that the quantities $\tilde{x}^{(t-1)}=x^{(t-1)}-\mathbf{\bar{x}}^{(t-1)}$
and $\tilde{y}^{(t-1)}=y^{(t-1)}-\mathbf{\bar{y}}^{(t-1)}$ denote the consensus errors of the local iterates and the auxiliary variables, respectively. We now introduce the following non-negative quantities: 
\begin{align}
&\tilde{\omega}_{1}:=\frac{\lambda^{B}}{\lambda^{B}-\delta}
\sum_{t=1}^{B}\lambda^{1-t}\left\Vert\tilde{y}^{(t-1)}\right\Vert_{L_{2}},
\qquad
\hat{\omega}_{1}:=\frac{\lambda^{B}}{\lambda^{B}-\delta}\cdot 2B\sigma\sqrt{N},\label{defn:tilde:omega:1:main}
\\
&\tilde{\omega}_{2}:= 0,
\qquad \qquad \qquad \qquad \qquad \qquad \quad
\hat{\omega}_{2}:=0,\label{defn:tilde:omega:2:main}
\\
&\tilde{\omega}_{3}:=\sqrt{N} \mgb{\left(\mathbb{E}\left\Vert\bar{x}^{(0)}-x_{\ast}\right\Vert^2\right)^{1/2}}
\qquad
\hat{\omega}_{3}:=\frac{\sqrt{N}}{\lambda}\left(\sqrt{\frac{L(1+\alpha)}{\mu\alpha}+\beta}\right)\frac{1}{\mu}\left(\sigma+\sqrt{\frac{2d}{\eta}}\right),\label{defn:tilde:omega:3:main}
\\
&\tilde{\omega}_{4}:=\frac{\lambda^{B}}{\lambda^{B}-\delta}
\sum_{t=1}^{B}\lambda^{1-t}\left\Vert\tilde{x}^{(t-1)}\right\Vert_{L_{2}},
\qquad
\hat{\omega}_{4}:=\frac{\lambda^{B}}{\lambda^{B}-\delta}\cdot B\sqrt{2\eta Nd},\label{defn:tilde:omega:4:main}
\end{align}
for $t=1,\ldots,B$.
\footnote{Because $\tilde{\omega}_2=\hat{\omega}_2=0$, these terms do not contribute to the bounds in Theorems~\ref{thm:average} and~\ref{thm:undirected}.} We also introduce the scalars
\begin{align}
&\gamma_{1}:=\frac{\lambda(1-\lambda^{B})}{(\lambda^{B}-\delta)(1-\lambda)},
\qquad
\gamma_{2}:=L\left(1+\frac{1}{\lambda}\right),\label{defn:gamma:1:2:main}
\\
&\gamma_{3}:=\left(1+\frac{{N}}{\lambda}\left(\sqrt{\frac{L(1+\alpha)}{\mu\alpha}+\beta}\right)\right),
\qquad
\gamma_{4}:=\frac{\eta(1-\lambda^{B})}{(\lambda^{B}-\delta)(1-\lambda)}.\label{defn:gamma:3:4:main}
\end{align}
Note that when the stepsize $\eta>0$ is sufficiently small, we have $\gamma_{1}\gamma_{2}\gamma_{3}\gamma_{4}\in(0,1)$. \mgb{For $r,q\in[0,1)$ and $k\geq 1$, we also define
\begin{equation}\label{eq:def:Hk}
\mathsf H_k(r,q)
:=
\sum_{j=0}^{k-1} r^{k-1-j}q^j
=
\begin{cases}
\dfrac{r^k-q^k}{r-q}, & r\neq q,\\[4pt]
kr^{k-1}, & r=q,
\end{cases}
\qquad
\mathsf H_0(r,q):=0.
\end{equation}
This quantity is the exact finite convolution of the geometric sequences
$\{r^k\}$ and $\{q^k\}$. Its asymptotic decay is governed, up to the
polynomial prefactor in the degenerate case $r=q$, by the more slowly
decaying of the two geometric rates.}

We now present our main results (Theorem~\ref{thm:average} and Theorem~\ref{thm:undirected}). The proofs of these results are given in Section~\ref{sec-proof-main-results}. They show that, when the stepsize $\eta$ is chosen sufficiently small---so as to satisfy certain inequalities ensuring the stability of the algorithm (in particular, boundedness of the $L_2$ norms of the iterates)---the distributions of the agents’ iterates converge linearly (i.e., at a geometric rate) in the $2$-Wasserstein distance to a neighborhood of the target distribution. Moreover, the radius of this neighborhood scales as $\mathcal{O}(\sqrt{\eta})$ as $\eta\to 0$. 

We first establish $2$-Wasserstein convergence guarantees showing that the law, $\mathrm{Law}(\bar{x}^{(k)})$, of the average iterate $\bar{x}^{(k)}$ converges to \mgb{an $\mathcal{O}(\sqrt{\eta})$-neighborhood of} the target Gibbs distribution $\pi\propto e^{-f}$.

\begin{theorem}\label{thm:average}
Consider the DIGing-SGLD algorithm with constant stepsize $\eta>0$. Assume that $\left\Vert x^{(0)}\right\Vert_{L_{2}}$ is finite. Let $\alpha,\beta>0$ be fixed scalars, and let $\lambda \in (\delta^{1/B},1)$, where $\delta\in(0,1)$ is as given in Assumption~\ref{assump:W}. Assume that the stepsize $\eta>0$ is chosen such that the following conditions are satisfied:
\begin{equation}\label{assump:alpha:beta:main}
\sqrt{1-\frac{\eta\mu\beta}{\beta+1}}\leq\lambda<1,
\quad
\eta\leq\frac{1}{(1+\alpha)L},
\quad
\quad
\gamma_{1}\gamma_{2}\gamma_{3}\gamma_{4} \in (0,1),
\end{equation}
where $\gamma_{1},\gamma_{2},\gamma_{3},\gamma_{4}$ are defined in \eqref{defn:gamma:1:2:main}--\eqref{defn:gamma:3:4:main}. Then, for every iteration $k$, we have
\begin{align}
\mathcal{W}_{2}\left(\mathrm{Law}\left(\bar{x}^{(k)}\right),\pi\right)\leq E_{1}(k,\eta)+E_{2}(k,\eta,\delta),
\label{ineq-wasserstein-bound}
\end{align}
where $E_{1}=E_{1}(k,\eta)$
and $E_{2}=E_{2}(k,\eta,\delta)$ are given by:
{\small
\begin{align}
E_{1}&:=
(1-\mu\eta)^{k}\left(\left(\mathbb{E}\left\Vert \bar{x}^{(0)}-x_{\ast}\right\Vert^{2}\right)^{1/2}+\sqrt{2\mu^{-1}dN^{-1}}\right)
+\frac{1.65L}{\mu}\sqrt{\eta dN^{-1}},\label{E:1:defn}
\\
E_{2}&:=\eta^{1/2}\left(\frac{\eta}{\mu(1-\frac{\eta L}{2})}+\frac{(1+\eta L)^{2}}{\mu^{2}(1-\frac{\eta L}{2})^{2}}\right)^{1/2}\nonumber
\cdot
\left(\frac{3L^{2}D^{2}\eta  \delta^{-2}}{N(1-\delta^{\frac{1}{B}})^{2}}
+\frac{6dL^{2}\cdot  \delta^{-2}}{1-\delta^{\frac{2}{B}}}
\right)^{1/2}
\\
&\quad\quad
+\frac{\sqrt{\eta}\sigma}{\sqrt{\mu(1-\frac{\eta L}{2})N}}
+\left[\eta\left(\eta+\frac{(1+\eta L)^2}{\mu(1-\eta L/2)}\right)
\mathsf H_k\!\left(1-\eta\mu(1-\eta L/2),\delta^{2/B}\right)\right]^{1/2}
\frac{\sqrt{3}L\delta^{-1}}{\sqrt{N}}\left\Vert x^{(0)}\right\Vert_{L_{2}}\label{E:2:defn}
\end{align}}
where $x_{\ast}$ is the minimizer of $f$,
$\bar{x}^{(0)}=\frac{1}{N}\sum_{i=1}^{N}x_{i}^{(0)}$, and
\begin{align}
D &:= \Bigg[ 
2\left(\frac{\gamma_{1}\gamma_{2}\gamma_{3}(\tilde{\omega}_{4}+\hat{\omega}_{4})
+\gamma_{1}\gamma_{2}(\tilde{\omega}_{3}+\hat{\omega}_{3})
+\tilde{\omega}_{1}+\hat{\omega}_{1}}
{1-\gamma_{1}\gamma_{2}\gamma_{3}\gamma_{4}}\right)^{2} \nonumber \\
&\qquad\quad
+\mgb{4L^2}
\left(\frac{\gamma_{3}\gamma_{4}(\tilde{\omega}_{1}+\hat{\omega}_{1})
+\gamma_{3}(\tilde{\omega}_{4}+\hat{\omega}_{4})
+\tilde{\omega}_{3}+\hat{\omega}_{3}}
{1-\gamma_{1}\gamma_{2}\gamma_{3}\gamma_{4}}\right)^{2}
+\mgb{4} 
\sigma^{2}
\Bigg]^{1/2}.
\label{eqn:D1:main}
\end{align}
Here, $\tilde{\omega}_i$ and $\hat{\omega}_i$ are defined in \eqref{defn:tilde:omega:1:main}--\eqref{defn:tilde:omega:4:main}
for $i=1,2,3,4$,
and $\pi$ is the Gibbs distribution with probability density function proportional to $\exp(-f(x))$.
\end{theorem}

Next, we present non-asymptotic convergence guarantees for the average $2$-Wasserstein distance between the distributions of the agent-wise iterates and the Gibbs distribution.

\begin{theorem}\label{thm:undirected} 
Under the same setting as in Theorem~\ref{thm:average}, for every iteration $k$, we have 
\begin{align}
\frac{1}{N}\sum_{i=1}^{N}\mathcal{W}_{2}\left(\mathrm{Law}\left(x_{i}^{(k)}\right),\pi\right)
\leq E_{1}(k,\eta)+E_{2}(k,\eta,\delta)+E_{3}(k,\eta,\delta),
\label{ineq-wass-distance-each-node}
\end{align}
where $E_{1}=E_{1}(k,\eta)$
and $E_{2}=E_{2}(k,\eta,\delta)$ are given in \eqref{E:1:defn}--\eqref{E:2:defn} and
$E_{3}=E_{3}(k,\eta,\delta)$ is defined as:
\begin{align}
E_{3}:=\frac{\sqrt{3} \delta^{-1}\delta^{\frac{k}{B}}}{\sqrt{N}}\left\Vert x^{(0)}\right\Vert_{L_{2}}
+\frac{\sqrt{3}D\eta  \delta^{-1}}{\sqrt{N}(1-\delta^{\frac{1}{B}})}
+\frac{\sqrt{6d\eta} \delta^{-1}}{\sqrt{1-\delta^{\frac{2}{B}}}},\label{E:3:defn}
\end{align}
where $x_{\ast}$ is the minimizer of $f$,
$\bar{x}^{(0)}=\frac{1}{N}\sum_{i=1}^{N}x_{i}^{(0)}$, $D$ is given in \eqref{eqn:D1:main}, 
and $\pi$ is the Gibbs distribution with the probability density function proportional to $\exp(-f(x))$.
\end{theorem}

\begin{remark}[Interpretations of $E_{1},E_{2}$ and $E_{3}$]
In Theorem~\ref{thm:average} and Theorem~\ref{thm:undirected}, as we will elaborate further in Section \ref{sec-proof-main-results}, $E_{1}$ serves as an upper bound on the $2$-Wasserstein distance between the distribution of $x_{k}$, the $k$-th iterate of the (centralized) unadjusted Langevin algorithm, and the Gibbs distribution. Moreover, $E_{2}$ is an upper bound on the $2$-Wasserstein distance between the distribution of the average of iterates $\bar{x}^{(k)}$ and that of $x_{k}$, and captures the additional error induced by stochastic gradients and decentralized averaging relative to the centralized unadjusted Langevin dynamics. Finally, $E_{3}$ is an upper bound on the averaged $2$-Wasserstein distance between the distribution of each agent $x_{i}^{(k)}$ and that of the average of iterates $\bar{x}^{(k)}$. Controlling each of these terms allows us to bound the Wasserstein distance between the distribution of each agent and the target, by an application of the triangle inequality for the $2$-Wasserstein distance.
\end{remark}

\begin{remark}[Feasibility of parameter choices]
Theorem~\ref{thm:average} and Theorem~\ref{thm:undirected} ensure that the constraints in \eqref{assump:alpha:beta:main} can be satisfied simultaneously. Indeed, fix any $\alpha,\beta>0$, and choose $\lambda$ such that $\tfrac{\delta+1}{2}\leq \lambda^{B}<1$. Under this choice, the quantities $\tilde{\omega}_{i},\hat{\omega}_{i},\gamma_{i}$ defined in \eqref{defn:tilde:omega:1:main}--\eqref{defn:gamma:3:4:main} 
remain uniformly bounded, i.e., all are of order $\mathcal{O}(1)$ except\footnote{\mgb{Here, $\tilde{\omega}_1$ and $\tilde{\omega}_4$ also absorb transient contributions from the iterates $x^{(t)}$ and $y^{(t)}$ over the initial $B$ iterations, prior to the onset of network contraction. Consequently, the notation $\mathcal{O}(\cdot)$ suppresses constants that may depend on the initialization $\|x^{(0)}\|{L_2}$, as well as on $B$, $N$, $d$, $L$, $\sigma$, and an upper bound $\bar{\eta}$ on the stepsize. In particular, explicit upper bounds for $\tilde{\omega}_i$, $\hat{\omega}_i$, and $\gamma_i$ are derived in the supplementary material as part of the proof of Lemma~\ref{lem-bound-D}. Since these calculations are routine but somewhat lengthy, we omit them from the main text.}} for $\hat{\omega}_{3}$, which is of order $\mathcal{O}\left(1/\sqrt{\eta}\right)$. Next, select $\eta>0$ sufficiently small (independently of $\lambda$) so that  
$\gamma_{1}\gamma_{2}\gamma_{3}\gamma_{4}\leq \mgb{q_0}$ for some $\mgb{q_0}\in (0,1)$
and
$\eta\leq \frac{1}{(1+\alpha)L}$.
With this choice, $D$ in \eqref{eqn:D1:main} is of order $\mathcal{O}\left(1/\sqrt{\eta}\right)$, so that $D\sqrt{\eta} = \mathcal{O}(1)$ as $\eta\to 0$. 
\mgb{If we let
$
r:=1-\eta\mu\left(1-\frac{\eta L}{2}\right)$ and $
q:=\delta^{2/B}.
$
The transient term in $E_2$ involving $\mathsf H_k(r,q)$ satisfies,
for $r\neq q$,
$
\sqrt{\mathsf H_k(r,q)}
=
\mathcal O\!\left(\max\{r,q\}^{k/2}\right).
$
Consequently, the upper bound $E_1+E_2+E_3$ on the Wasserstein distance given in \eqref{ineq-wass-distance-each-node} scales like
\[
E_1+E_2+E_3
=
\mathcal O\!\left((1-\mu\eta)^k\right)
+
\mathcal O\!\left(
\max\left\{
\sqrt{1-\eta\mu(1-\eta L/2)},
\delta^{1/B}
\right\}^{k}
\right)
+
\mathcal O(\sqrt{\eta}),
\]
}
which can be made arbitrarily small by choosing the stepsize $\eta>0$ small enough and $k$ large enough.\footnote{\mgb{If the rate $r=q$, the corresponding convolution term is
$
\sqrt{\mathsf H_k(r,r)}
=
\sqrt{k}\,r^{(k-1)/2},
$
so the same conclusion holds with this additional polynomial
prefactor.}} Finally, once $\eta$ is fixed, we can always choose $\lambda$ so that  
$\sqrt{1-\frac{\eta\mu\beta}{\beta+1}}\;\leq\;\lambda<1$,
which is possible since the square-root term is strictly below $1$ for any $\eta>0$.  Therefore, feasible pairs $(\eta,\lambda)$ that satisfy the conditions \eqref{assump:alpha:beta:main} always exist. Moreover, as $k\to\infty$, the upper bound $E_1 + E_2 +E_3 = \mathcal{O}(\sqrt{\eta})$. This shows that, for any target accuracy $\epsilon>0$, by choosing the stepsize $\eta$ sufficiently small (as a function of $\epsilon$) and the number of iterates $k$ large enough, we can ensure each agent's iterates are close to the target distribution, i.e., 
$\frac{1}{N}\sum_{i=1}^{N}\mathcal{W}_{2}\!\left(\mathrm{Law}\!\left(x_{i}^{(k)}\right),\pi\right)\leq \epsilon$.
\end{remark}

In the following discussion, we 
present a family of explicit choices of the stepsize $\eta$ and the parameters $\lambda,\alpha,\beta$, as well as the corresponding number of iterations $k$ 
required to ensure that the error satisfies $E_1+E_2+E_3\leq \epsilon$. Other admissible choices of $\eta,\lambda,\alpha,\beta$ are also possible and can be obtained by numerically optimizing the upper bound $E_1+E_2+E_3$ as a function of these parameters; 
however, our goal here is to provide fully explicit, constructive parameter selections.

We first present a lemma 
that specifies an explicit admissible range of stepsizes $\eta$ together with a corresponding choice of $\lambda(\eta)$ and fixed values of $\alpha$ and $\beta$, and establishes an explicit bound on the quantity $D\sqrt{\eta}$, showing that it is $\mathcal{O}(1)$ as $\eta\to 0$. The resulting bound will be used in the subsequent corollary to derive an explicit iteration-complexity guarantee for DIGing-SGLD. The main challenge in deriving such explicit bounds and parameter choices is that $D$ depends jointly on $\eta$, $\lambda$, $\alpha$, and $\beta$, which are coupled through the nonlinear conditions in \eqref{assump:alpha:beta:main}. 
{
\begin{lemma}\label{lem-bound-D}
In the setting of Theorem~\ref{thm:average}, let
\begin{equation}\label{def:J1-corrected-main}
\kappa:=\frac{L}{\mu},\qquad
J_1:=3\kappa B^2\bigl(1+4N\sqrt{\kappa}\bigr).
\end{equation}
and define
\begin{align}
\eta_{\rm sg}
&:=\frac{1.5\bigl(\sqrt{J_1^2+(1-\delta^2)J_1}-\delta J_1\bigr)^2}
{\mu J_1(J_1+1)^2},\label{def:eta-sg-main}\\
\bar\eta
&:=\min\left\{\eta_{\rm sg},\;\frac{1}{2L},\;
\frac{3}{4\mu}\bigl(1-\delta^{2/B}\bigr)\right\}.\label{bar:eta:defn}
\end{align}
For every \(\eta\in(0,\bar\eta]\), take \(\alpha=1\), \(\beta=2L/\mu\), and
\begin{equation}\label{def-opt-lambda-eta}
\lambda(\eta):=\left(1-\frac{\eta\mu}{1.5}\right)^{1/(2B)}.
\end{equation}
Set
\begin{equation}\label{def:lambda-rho-q-main}
\underline\lambda:=\lambda(\bar\eta),\qquad
\rho_0:=1-\frac{\bar\eta\mu}{1.5}-\delta^{2/B},\qquad
q_0:=\frac{1+2N\sqrt\kappa/\underline\lambda}{1+4N\sqrt\kappa}.
\end{equation}
Then
\begin{equation}
\delta^{1/B}<\underline\lambda\leq\lambda(\eta)<1,\qquad
\rho_0>0,\qquad 0<q_0<1,
\label{ineq-lower-bound-rho-zero}
\end{equation}
and all the conditions in \eqref{assump:alpha:beta:main} are satisfied. In particular,
\begin{equation}\label{eq:small-gain-corrected-main}
0<\gamma_1\gamma_2\gamma_3\gamma_4\leq q_0<1.
\end{equation}

To give an explicit bound for \(D\sqrt\eta\), define the finite-horizon deterministic bounds
\begin{align}
X_0&:=\|x^{(0)}\|_{L_2},\qquad Y_0:=\|y^{(0)}\|_{L_2},\nonumber\\
X_{t+1}&:=X_t+\bar\eta Y_t+\sqrt{2\bar\eta Nd},\nonumber\\
Y_{t+1}&:=(1+L\bar\eta)Y_t+2LX_t+L\sqrt{2\bar\eta Nd}+2\sqrt N\,\sigma,
\qquad t=0,\ldots,B-2,\label{def:finite-horizon-bounds-main}
\end{align}
in a recursive manner. Also set
\begin{align}
C_\lambda&:=\frac{\underline\lambda^B}{\underline\lambda^B-\delta},\nonumber\\
\Gamma_1&:=\frac{B}{\underline\lambda^B-\delta},\qquad
\Gamma_2:=L\left(1+\frac1{\underline\lambda}\right),\nonumber\\
\Gamma_3&:=1+\frac{2N\sqrt\kappa}{\underline\lambda},\qquad
\Gamma_4:=\frac{B\bar\eta}{\underline\lambda^B-\delta},\label{def:Gamma-main}\\
\Omega_1&:=C_\lambda\sqrt{\bar\eta}\left(
\sum_{t=1}^{B}\underline\lambda^{1-t}Y_{t-1}+2B\sigma\sqrt N\right),\nonumber\\
\Omega_3&:=\sqrt{N\bar\eta}\,(\mathbb{E}\|\bar x^{(0)}-x_*\|^2)^{1/2}
+\frac{2\sqrt{N\kappa}}{\mu\underline\lambda}
\left(\sigma\sqrt{\bar\eta}+\sqrt{2d}\right),\nonumber\\
\Omega_4&:=C_\lambda\left(
\sqrt{\bar\eta}\sum_{t=1}^{B}\underline\lambda^{1-t}X_{t-1}
+B\bar\eta\sqrt{2Nd}\right).\label{def:Omega-main}
\end{align}
Then \(D\sqrt\eta\), with \(D\) defined in \eqref{eqn:D1:main}, satisfies
\begin{align}
D\sqrt\eta\leq \overline D:=\Bigg[&
2\left(
\frac{\Gamma_1\Gamma_2\Gamma_3\Omega_4+
\Gamma_1\Gamma_2\Omega_3+\Omega_1}{1-q_0}
\right)^2\nonumber\\
&+4L^2\left(
\frac{\Gamma_3\Gamma_4\Omega_1+\Gamma_3\Omega_4+\Omega_3}{1-q_0}
\right)^2
+4\sigma^2\bar\eta\Bigg]^{1/2}.
\label{eqn:D1:main:upperbound}
\end{align}
The constant \(\overline D\) is finite and is independent of the target accuracy \(\epsilon\).
\end{lemma}
}
\begin{proof}
\mgb{The proof proceeds in two steps. First, it verifies that the proposed choices of the stepsize $\eta$ and parameter $\lambda(\eta)$ satisfy the conditions in \eqref{assump:alpha:beta:main}, with the explicit bound
$
\gamma_1\gamma_2\gamma_3\gamma_4\le q_0<1.
$
Second, it derives uniform finite-horizon bounds for the transient quantities $\hat{\omega}_i$ and $\tilde{\omega}_i$ appearing in the definition of $D$, and combines these bounds with the stability margin $1-q_0>0$ to obtain the stated estimate on $D\sqrt{\eta}$. The calculations are routine but somewhat lengthy, and are therefore deferred to the supplementary material.}
\end{proof}

Next, we establish an iteration complexity bound that quantifies how many iterations of DIGing-SGLD are required to ensure that the average $2$-Wasserstein error does not exceed $\epsilon$. To obtain this bound, we propose a stepsize choice that depends explicitly on the target accuracy $\epsilon$ and build on the previous lemma.
{
\begin{corollary}\label{coro-param}
In the setting of Theorem~\ref{thm:undirected}, let the target accuracy \(\epsilon>0\) be given and use the parameter choices of Lemma~\ref{lem-bound-D}. Define
\begin{align*}
P_0&:=\left[\bar\eta\left(\bar\eta+
\frac{(1+\bar\eta L)^2}{\mu(1-\bar\eta L/2)}\right)\right]^{1/2},\\
\overline C_1&:=
\left(\mathbb E\|\bar x^{(0)}-x_*\|^2\right)^{1/2}
+\sqrt{2\mu^{-1}dN^{-1}},\\
\overline C_2&:=
\frac{P_0}{\sqrt{\rho_0}}
\frac{\sqrt3 L\delta^{-1}}{\sqrt N}\|x^{(0)}\|_{L_2}
+\frac{\sqrt3\delta^{-1}}{\sqrt N}\|x^{(0)}\|_{L_2},\\
\overline C_3&:=
\frac{1.65L}{\mu}\sqrt{dN^{-1}}
+\frac{\sqrt{6d}\,\delta^{-1}}{\sqrt{1-\delta^{2/B}}}
+\frac{2\sigma}{\sqrt{3\mu N}}
+\frac{2}{\mu}\left(\frac{6dL^2\delta^{-2}}{1-\delta^{2/B}}\right)^{1/2}\\
&\qquad
+\frac{\sqrt3\delta^{-1}\overline D}{\sqrt N(1-\delta^{1/B})}
+\frac{2\overline D}{\mu}
\left(\frac{3L^2\delta^{-2}}{N(1-\delta^{1/B})^2}\right)^{1/2},\\
\overline C_4&:=
\frac{2}{\sqrt{3\mu}}
\left(\frac{6dL^2\delta^{-2}}{1-\delta^{2/B}}\right)^{1/2}
+\frac{2}{\sqrt{3\mu}}
\left(\frac{3L^2\delta^{-2}}{N(1-\delta^{1/B})^2}\right)^{1/2}\overline D,
\end{align*}
where \(\rho_0\) and \(\overline D\) are defined in
\eqref{def:lambda-rho-q-main} and \eqref{eqn:D1:main:upperbound}, respectively. Let
\begin{equation}\label{def:eta-star-corrected-main}
\eta_*:=\min\{\bar\eta,\eta_{\rm noise}(\epsilon)\},\qquad
\eta_{\rm noise}(\epsilon):=
\min\left\{\frac{\epsilon^2}{9\overline C_3^2},
\frac{\epsilon}{3\overline C_4}\right\}.
\end{equation}
Then DIGing-SGLD satisfies
\[
\frac1N\sum_{i=1}^N
\mathcal W_2\!\left(\operatorname{Law}(x_i^{(k)}),\pi\right)\leq\epsilon
\]
whenever
\begin{equation}\label{def:k-star-corrected-main}
k\geq k_*(\epsilon):=
\frac{3}{\mu\eta_*}
\log\left(\frac{3(\overline C_1+\overline C_2)}{\epsilon}\right).
\end{equation}
\end{corollary}
}

\begin{proof}
The proof is based on deriving upper bounds for the error terms $E_1$, $E_2$, and $E_3$ that appear in Theorem~\ref{thm:undirected}, and showing that when $k \geq k_*(\epsilon)$, their sum does not exceed $\epsilon$, building on Lemma~\ref{lem-bound-D}. The required computations are straightforward but somewhat lengthy, and are therefore provided in the supplementary material.
\end{proof}

{
\begin{remark}
By \eqref{bar:eta:defn} and \eqref{ineq-lower-bound-rho-zero}, $\bar\eta\leq \frac{3}{4\mu}(1-\delta^{2/B})$ and therefore
$
\rho_0=1-\frac{\bar\eta\mu}{1.5}-\delta^{2/B}
\geq \frac{1-\delta^{2/B}}{2}>0.
$
Thus, the denominator of \(\overline C_2\) is positive and  \(\overline C_2\) is well-defined. The constants \(\overline C_i\), \(i=1,2,3,4\), are independent of the target accuracy \(\epsilon\). Consequently, as \(\epsilon\to0\),
\(\eta_*=\Theta(\epsilon^2)\), and hence
\[
k_*(\epsilon)=\Theta\!\left(\frac{\log(1/\epsilon)}{\epsilon^2}\right).
\]
Therefore, it suffices to have $\Theta\!\left(\frac{\log(1/\epsilon)}{\epsilon^2}\right)$ iterations of DIGing-SGLD for the average distance to the target satisfy
$
\frac1N\sum_{i=1}^N
\mathcal W_2\!\left(\operatorname{Law}(x_i^{(k)}),\pi\right)\leq\epsilon
$.
\end{remark}
}

\subsection{Proofs of the Main Results}\label{sec-proof-main-results}

In this section, we present the proofs of Theorem~\ref{thm:average} and Theorem~\ref{thm:undirected} by establishing a sequence of technical lemmas
whose proofs will be provided in Appendix~\ref{proofs:technical}. 
To prove Theorem~\ref{thm:average} and Theorem~\ref{thm:undirected}, based on the triangle inequality for the 2-Wasserstein distance, we consider the following decompositions:
\begin{align}
\frac{1}{N}\sum_{i=1}^N\mathcal{W}_2\left(\mathrm{Law}\left(x_i^{(k)}\right)\,,\, \pi \right) 
&\leq \frac{1}{N}\sum_{i=1}^N\mathcal{W}_2\left(\mathrm{Law}\left(x_i^{(k)}\right)\,,\, \mathrm{Law}\left(\overline{x}^{(k)}\right) \right) + \mathcal{W}_2\left(\mathrm{Law}\left(\overline{x}^{(k)}\right)\,,\,\pi \right)\label{decompose:0},
\end{align}
and
\begin{align}
\mathcal{W}_2\left(\mathrm{Law}\left(\overline{x}^{(k)}\right)\,,\,\pi \right)\leq \mathcal{W}_2\left(\mathrm{Law}\left(\overline{x}^{(k)}\right),\mathrm{Law}(x_{k}) \right)
+\mathcal{W}_{2}\left(\mathrm{Law}(x_{k}),\pi\right).
\label{decompose:00}
\end{align}
Here, $\overline{x}^{(k)}:=\frac{1}{N}\sum_{i=1}^{N}x_{i}^{(k)}$
is the average of the agent iterates and $x_{k}$ is defined via the iteration
\begin{equation}\label{overdamped:iterates:0}
x_{k+1} = x_{k} - \frac{\eta}{N}\nabla f\left(x_k\right) + \sqrt{2\eta}\overline{w}^{(k+1)},
\end{equation}
which correspond to the Euler-Maruyama discretization of overdamped Langevin diffusion 
\begin{equation}\label{overdamped:SDE}
dX_t = -\frac{1}{N}\nabla f(X_t)dt + \sqrt{2N^{-1}}dW_t,
\end{equation}
where $W_{t}$ is a standard $d$-dimensional Brownian motion, $\overline{w}^{(k)}:=\frac{1}{N}\sum_{i=1}^{N}w_{i}^{(k)}$, 
and $w_{i}^{(k)}$ are $\mathcal{N}(0,I_{d})$ distributed that are i.i.d. 
in both $k\in\mathbb{N}$ and $i=1,2,\ldots,N$.

The main idea of our proof technique is to bound the following three terms: (1) the $L_{2}$ distance between $x_{i}^{(k)}$ and their average $\bar{x}^{(k)}$; (2) the $L_{2}$ distance between the average iterate $\bar{x}^{(k)}$ and iterates $x_{k}$ in \eqref{overdamped:iterates:0} obtained from Euler-Maruyama discretization of overdamped diffusion \eqref{overdamped:SDE}; and (3) the $\mathcal{W}_{2}$ distance between the law of $x_{k}$ in \eqref{overdamped:iterates:0} and the Gibbs distribution $\pi$. First, we upper bound the $L_{2}$ distance between $x_{i}^{(k)}$ and their average.

\subsubsection{Uniform $L_{2}$ bounds between $x_{i}^{(k)}$ and their average $\bar{x}^{(k)}$}

In this section, we derive uniform $L_{2}$ bounds between $x_{i}^{(k)}$ and their average $\bar{x}^{(k)}$, where ``uniform'' refers to validity for all iterations 
$k$. As a first step, we derive a
uniform $L_{2}$ bound for $y^{(k)}$, 
which is a key ingredient of our analysis.
First, we recall from the notations introduced in Section~\ref{sec:background} that
\begin{equation}
\tilde{x}^{(k)}=x^{(k)}-\mathbf{\bar{x}}^{(k)},
\qquad
\tilde{y}^{(k)}=y^{(k)}-\mathbf{\bar{y}}^{(k)},
\end{equation}
where we recall from \eqref{x:vector:form}--\eqref{y:vector:form} that
$x^{(k)},y^{(k)}$ satisfy the iterations:
\begin{align}
&x^{(k+1)}=\mathcal{W}^{(k)}x^{(k)}-\eta y^{(k)}+\sqrt{2\eta}w^{(k+1)},
\\
&y^{(k+1)}=\mathcal{W}^{(k)}y^{(k)}+\nabla F\left(x^{(k+1)}\right)-\nabla F\left(x^{(k)}\right)
+\xi^{(k+2)}-\xi^{(k+1)}.
\end{align}

\begin{lemma}\label{lem:0} Let $\alpha,\beta>0$ and $\lambda \in (\delta^{1/B},1)$ be given and fixed, where $\delta:=\sup_{k\geq B-1}\delta(k)$ with $\delta(k)$ defined in Assumption~\ref{assump:W}. Assume that the conditions listed in \eqref{assump:alpha:beta:main} hold, i.e., 
\begin{equation}\label{assump:alpha:beta}
\sqrt{1-\frac{\eta\mu\beta}{\beta+1}}\leq\lambda<1,
\qquad
\eta\leq\frac{1}{(1+\alpha)L},
\qquad\text{and}\qquad
\gamma_{1}\gamma_{2}\gamma_{3}\gamma_{4}\in(0,1),\end{equation} where $\gamma_{1}, \gamma_{2},\gamma_{3}$ and $\gamma_{4}$ are defined by \eqref{defn:gamma:1:2:main}--\eqref{defn:gamma:3:4:main}. Then, for every $k$, 
\begin{equation}
\mathbb{E}\left\Vert y^{(k)}\right\Vert^{2}\leq D^{2},
\end{equation}
where the quantity $D$ is defined in \eqref{eqn:D1:main}.
\end{lemma}

The proof of Lemma~\ref{lem:0}, which is deferred to Appendix~\ref{app-sec-proof-lem:0}, relies on a sequence of technical lemmas, which we will introduce next.
For any $k=0,1,2,\ldots,$ we define:
\begin{equation}\label{q:z:k:defn}
q^{(k)}:=x^{(k)}-\mathbf{x}_{\ast},
\qquad
z^{(k)}:=\nabla F\left(x^{(k)}\right)-\nabla F\left(x^{(k-1)}\right),
\end{equation}
where $\mathbf{x}_{\ast}=\left[x_{\ast}^{\top},x_{\ast}^{\top},\ldots,x_{\ast}^{\top}\right]^{\top}$ \mgb{with the convention that $x^{(-1)}:= x^{(0)}$ so that $z^{(0)}=0$}.
Inspired by \cite{nedic2017}, we introduce the weighted $L_2$ norms
\begin{align}
\Vert q\Vert_{L_{2}}^{\lambda,K}
:=\max_{0,1,\ldots,K}\frac{1}{\lambda^{k}}\left(\mathbb{E}\left\Vert q^{(k)}\right\Vert^{2}\right)^{1/2},
\qquad
\Vert z\Vert_{L_{2}}^{\lambda,K}
:=\max_{0,1,\ldots,K}\frac{1}{\lambda^{k}}\left(\mathbb{E}\left\Vert z^{(k)}\right\Vert^{2}\right)^{1/2},
\end{align}
and the following loop
\begin{equation}\label{eqn:arrows}
q\rightarrow z\rightarrow\tilde{y}\rightarrow\tilde{x}\rightarrow q,
\end{equation}
as a proof technique. Here, each arrow means that the (weighted) $L_2$ norm of the sequence at the head of the arrow, can be controlled by the (weighted) $L_2$ norm of the sequence at the tail of the arrow; we will explain below what exactly we mean by this. For example, the arrow $q\rightarrow z$ in \eqref{eqn:arrows} means that we would like
to establish an upper bound on $\Vert z\Vert_{L_{2}}^{\lambda,K}$
using $\Vert q\Vert_{L_{2}}^{\lambda,K}$. As we shall discuss next, treating each arrow separately, will allow us to complete the loop and control the boundedness of the $\tilde{y}=\{\tilde{y}_{k}\}_{k\geq 0}$ sequence, and this in return will allow us to ensure the boundedness of the $y$ sequence in $L_2$. 
Let us first study the first arrow $q\rightarrow z$ in \eqref{eqn:arrows}. 
We have the following technical lemma.

\begin{lemma}\label{lem:first:arrow}
For any $K=0,1,2,\ldots,$ and $\lambda\in(0,1)$, we have
\begin{equation}
\Vert z\Vert_{L_{2}}^{\lambda,K}\leq L\left(1+\frac{1}{\lambda}\right)\Vert q\Vert_{L_{2}}^{\lambda,K}.
\end{equation}
\end{lemma}

The proof is deferred to Appendix~\ref{appendix-sec-lem:first:arrow}. Next, let us consider the second arrow $z\rightarrow\tilde{y}$ in \eqref{eqn:arrows}. 
Recall from \eqref{y:vector:form} and \eqref{q:z:k:defn} that
\begin{equation}
y^{(k+1)}=\mathcal{W}^{(k)}y^{(k)}+z^{(k+1)}+\xi^{(k+2)}-\xi^{(k+1)}.
\end{equation}
In a similar manner as before, we define:
\begin{align}
\Vert\tilde{y}\Vert_{L_{2}}^{\lambda,K}
:=\max_{0,1,\ldots,K}\frac{1}{\lambda^{k}}\left(\mathbb{E}\left\Vert\tilde{y}^{(k)}\right\Vert^{2}\right)^{1/2}.
\end{align}
We will now establish an upper bound
on $\Vert\tilde{y}\Vert_{L_{2}}^{\lambda,K}$ using $\Vert z\Vert_{L_{2}}^{\lambda,K}$ and $\left\Vert\tilde{y}^{(t-1)}\right\Vert_{L_{2}}$
for $t=1,2,\ldots,B$.
We have the following technical lemma.

\begin{lemma}\label{lem:second:arrow}
Let $\delta:=\sup_{k\geq B-1}\delta(k)$, where $\delta(k)$ is defined in Assumption~\ref{assump:W}.
Let $\lambda$ be such that $\delta<\lambda^{B}<1$. Then for any $K=0,1,2,\ldots$, we have
\begin{equation}
\Vert\tilde{y}\Vert_{L_{2}}^{\lambda,K}
\leq
\frac{\lambda(1-\lambda^{B})}{(\lambda^{B}-\delta)(1-\lambda)}\Vert z\Vert_{L_{2}}^{\lambda,K}
+\frac{\lambda^{B}}{\lambda^{B}-\delta}\frac{2B\sigma\sqrt{N}}{\lambda^{K}}
+\frac{\lambda^{B}}{\lambda^{B}-\delta}
\sum_{t=1}^{B}\lambda^{1-t}\left\Vert\tilde{y}^{(t-1)}\right\Vert_{L_{2}}.
\end{equation}
\end{lemma}

The proof is deferred to Appendix~\ref{sec:lem:second:arrow}. Next, we consider the third arrow $\tilde{y}\rightarrow\tilde{x}$ in \eqref{eqn:arrows}, where our aim will be to obtain an upper bound on 
$\Vert\tilde{x}\Vert_{L_{2}}^{\lambda,K}$ using 
$\Vert\tilde{y}\Vert_{L_{2}}^{\lambda,K}$ and $\left\Vert\tilde{x}^{(t-1)}\right\Vert_{L_{2}}$ for $t=1,2,\ldots,B$. We have the following technical lemma.

\begin{lemma}\label{lem:third:arrow}
Let $\delta:=\sup_{k\geq B-1}\delta(k)$, where $\delta(k)$ is defined in Assumption~\ref{assump:W}.
Let $\lambda$ be such that $\delta<\lambda^{B}<1$. Then for any $K=0,1,2,\ldots$, we have
\begin{equation}
\Vert\tilde{x}\Vert_{L_{2}}^{\lambda,K}
\leq
\frac{\eta(1-\lambda^{B})}{(\lambda^{B}-\delta)(1-\lambda)}\Vert\tilde{y}\Vert_{L_{2}}^{\lambda,K}
+\frac{\lambda^{B}}{\lambda^{B}-\delta}\frac{B\sqrt{2\eta Nd}}{\lambda^{K}}
+\frac{\lambda^{B}}{\lambda^{B}-\delta}
\sum_{t=1}^{B}\lambda^{1-t}\left\Vert\tilde{x}^{(t-1)}\right\Vert_{L_{2}}.
\end{equation}
\end{lemma}

The proof is given in Appendix~\ref{appendix-proof-lem:third:arrow}.
Finally, let us consider the last arrow $\tilde{x}\rightarrow q$ in \eqref{eqn:arrows}, 
for which we would like to establish an upper bound on $\left\Vert q\right\Vert_{L_{2}}^{\lambda,K}$ by using $\Vert\tilde{x}\Vert_{L_{2}}^{\lambda,K}$. 
We have the following result.

\begin{lemma}\label{lem:last:arrow}
Assume that the parameters $\alpha,\beta>0$ satisfy
\begin{equation}
\sqrt{1-\frac{\eta\mu\beta}{\beta+1}}\leq\lambda<1,
\qquad\text{and}\qquad
\eta\leq\frac{1}{(1+\alpha)L}.
\end{equation}
Then, for every $K=0,1,2,\ldots$, we have
\begin{align*}
\left\Vert q\right\Vert_{L_{2}}^{\lambda,K}
&\leq
\sqrt{N} \mgb{\left(\mathbb{E}\left\Vert\bar{x}^{(0)}-x_{\ast}\right\Vert^2\right)^{1/2}} 
+\left(1+\frac{{N}}{\lambda}\left(\sqrt{\frac{L(1+\alpha)}{\mu\alpha}+\beta}\right)\right)
\Vert\tilde{x}\Vert_{L_{2}}^{\lambda,K}
\\
&\qquad\qquad
+\frac{\sqrt{N}}{\lambda}\left(\sqrt{\frac{L(1+\alpha)}{\mu\alpha}+\beta}\right)\frac{1}{\mu}\left(\sigma+\sqrt{\frac{2d}{\eta}}\right)\frac{1}{\lambda^{K}}.
\end{align*}
\end{lemma}

The proof is given in Appendix~\ref{app-section-lem:last:arrow}. It follows from Lemma~\ref{lem:first:arrow}, Lemma~\ref{lem:second:arrow}, Lemma~\ref{lem:third:arrow}
and Lemma~\ref{lem:last:arrow} that
\begin{align}
&\Vert\tilde{y}\Vert_{L_{2}}^{\lambda,K}\leq\gamma_{1}\Vert z\Vert_{L_{2}}^{\lambda,K}+\omega_{1}(K),\label{recursive:ineq:1}
\\
&\Vert z\Vert_{L_{2}}^{\lambda,K}\leq\gamma_{2}\Vert q\Vert_{L_{2}}^{\lambda,K}+\omega_{2}(K),\label{recursive:ineq:2}
\\
&\Vert q\Vert_{L_{2}}^{\lambda,K}\leq\gamma_{3}\Vert\tilde{x}\Vert_{L_{2}}^{\lambda, K}+\omega_{3}(K),\label{recursive:ineq:3}
\\
&\Vert\tilde{x}\Vert_{L_{2}}^{\lambda,K}\leq\gamma_{4}\Vert\tilde{y}\Vert_{L_{2}}^{\lambda, K}+\omega_{4}(K),\label{recursive:ineq:4}
\end{align}
where $\gamma_{1},\gamma_{2},\gamma_{3},\gamma_{4}$ are defined in \eqref{defn:gamma:1:2:main}--\eqref{defn:gamma:3:4:main} and
\begin{align}
&\omega_{1}(K):=\frac{\lambda^{B}}{\lambda^{B}-\delta}\frac{2B\sigma\sqrt{N}}{\lambda^{K}}
+\frac{\lambda^{B}}{\lambda^{B}-\delta}
\sum_{t=1}^{B}\lambda^{1-t}\left\Vert\tilde{y}^{(t-1)}\right\Vert_{L_{2}},
\qquad
\omega_{2}(K):=0,\label{defn:omega:1:2:K}
\\
&\omega_{3}(K):=\sqrt{N} \mgb{\left(\mathbb{E}\left\Vert\bar{x}^{(0)}-x_{\ast}\right\Vert^2\right)^{1/2}} 
+\frac{\sqrt{N}}{\lambda}\left(\sqrt{\frac{L(1+\alpha)}{\mu\alpha}+\beta}\right)\frac{1}{\mu}\left(\sigma+\sqrt{\frac{2d}{\eta}}\right)\frac{1}{\lambda^{K}},\label{defn:omega:3:K}
\\
&\omega_{4}(K):=\frac{\lambda^{B}}{\lambda^{B}-\delta}\frac{B\sqrt{2\eta Nd}}{\lambda^{K}}
+\frac{\lambda^{B}}{\lambda^{B}-\delta}
\sum_{t=1}^{B}\lambda^{1-t}\left\Vert\tilde{x}^{(t-1)}\right\Vert_{L_{2}}.\label{defn:omega:4:K}
\end{align}
Note that the subscripts of the error terms correspond to the four recursive inequalities in \eqref{recursive:ineq:1}--\eqref{recursive:ineq:4}; in particular, $\omega_2(K)=0$ because \eqref{recursive:ineq:2} contains no additive error. As an immediate consequence of \eqref{recursive:ineq:1}--\eqref{recursive:ineq:4}, we obtain the following technical lemma. 

\begin{lemma}\label{lem:y:q}
Assume that $\delta<\lambda^{B}<1$ with $\delta(k)$ defined in Assumption~\ref{assump:W}
and \eqref{assump:alpha:beta} holds.
Then, for every $K=0,1,2,\ldots$,
\begin{align}
&\Vert\tilde{y}\Vert_{L_{2}}^{\lambda,K}\leq\frac{\gamma_{1}\gamma_{2}\gamma_{3}\omega_{4}(K)+\gamma_{1}\gamma_{2}\omega_{3}(K)
+\gamma_{1}\omega_{2}(K)+\omega_{1}(K)}{1-\gamma_{1}\gamma_{2}\gamma_{3}\gamma_{4}},\label{tilde:y:upper:bound}
\\
&\Vert q\Vert_{L_{2}}^{\lambda,K}\leq\frac{\gamma_{3}\gamma_{4}\gamma_{1}\omega_{2}(K)+\gamma_{3}\gamma_{4}\omega_{1}(K)
+\gamma_{3}\omega_{4}(K)+\omega_{3}(K)}{1-\gamma_{1}\gamma_{2}\gamma_{3}\gamma_{4}},\label{q:upper:bound}
\end{align}
where $\gamma_{1},\gamma_{2},\gamma_{3},\gamma_{4}$ are defined in \eqref{defn:gamma:1:2:main}--\eqref{defn:gamma:3:4:main} and $\omega_{1}(K),\omega_{2}(K),\omega_{3}(K),\omega_{4}(K)$ are defined in \eqref{defn:omega:1:2:K}, \eqref{defn:omega:3:K} and \eqref{defn:omega:4:K}.
\end{lemma}


The proof is provided in Appendix~\ref{app-sec-proof-lem:y:q}.
Note that although $\omega_2(K)=0$, we retain it in the bounds \eqref{tilde:y:upper:bound} and \eqref{q:upper:bound} to make the four-term recursive structure explicit. Next, we present a technical lemma that upper bounds the averaged $L_{2}$ distance between the iterates $x_{i}^{(k)}$ and the average $\bar{x}^{(k)}$.

\begin{lemma}\label{lem:1}
Assume that $\delta<\lambda^{B}<1$ with $\delta(k)$ defined in Assumption~\ref{assump:W}
and \eqref{assump:alpha:beta} holds.
Then, for any $k\geq 1$, we have
\begin{equation*}
\sum_{i=1}^{N}\mathbb{E}\left\Vert x_{i}^{(k)}-\bar{x}^{(k)}\right\Vert^{2}
\leq
3\left(\bar{\gamma}_{k}^{(k-1)}\right)^{2}\mathbb{E}\left\Vert x^{(0)}\right\Vert^{2}
+3D^{2}\eta^{2}\left(\sum_{s=0}^{k-1}\bar{\gamma}_{k-1-s}^{(k-1)}\right)^{2}
+6dN\eta\sum_{s=0}^{k-1}\left(\bar{\gamma}_{k-1-s}^{(k-1)}\right)^{2},
\end{equation*}
where $D$ is defined in \eqref{eqn:D1:main} and
\begin{equation}
\bar{\gamma}_{j}^{(k-1)}:=\left\Vert W_{j}^{(k-1)}-\frac{1}{N}1_{N}1_{N}^{\top}\right\Vert
\quad
\mbox{for}\quad
j=0,1,\ldots, \mgb{k}.
\end{equation}
\end{lemma}

The proof is given in Appendix~\ref{app-sec-lem:1}. Next, we aim to provide an upper bound for $\bar{\gamma}_{k-1-s}^{(k)}$ in Lemma~\ref{lem:1}
under Assumption~\ref{assump:W} for the mixing matrices $W^{(k)}$. The following corollary of Lemma~\ref{lem:1}, whose proof is in Appendix~\ref{app-sec-proof-lem-cor:1}, establishes this and shows that the iterates $x_{i}^{(k)}$ are close to the average $\bar{x}^{(k)}$ on average.

\begin{lemma}\label{cor:1}
Assume that $\delta<\lambda^{B}<1$ with $\delta(k)$ defined in Assumption~\ref{assump:W}
and \eqref{assump:alpha:beta} holds.
Then, for any $k \geq 1$, we have
\begin{equation}
\sum_{i=1}^{N}\mathbb{E}\left\Vert x_{i}^{(k)}-\bar{x}^{(k)}\right\Vert^{2}
\leq
3\cdot \delta^{-2}\left(\delta^{\frac{2}{B}}\right)^{k}\mathbb{E}\left\Vert x^{(0)}\right\Vert^{2}
+\frac{3D^{2}\eta^{2}\delta^{-2}}{(1-\delta^{\frac{1}{B}})^{2}}
+\frac{6dN\eta\cdot \delta^{-2}}{1-\delta^{\frac{2}{B}}},
\end{equation}
where $D$ is defined in \eqref{eqn:D1:main}.
\end{lemma}


\subsubsection{$L_{2}$ distance between $\overline{x}^{(k)}$ and $x_{k}$}

In this section, we derive bounds on 
the $L_{2}$ distance between $\overline{x}^{(k)}$ and $x_{k}$, 
which is the $k$-th iterate of the Euler discretization of an overdamped Langevin diffusion given in \eqref{overdamped:iterates:0}.

First, by taking the average of the $N$ nodes in \eqref{DIGing:x}--\eqref{DIGing:y}, and using the fact that $W^{(k)}$ is doubly stochastic, we obtain:
\begin{equation}
\bar{x}^{(k+1)}=\bar{x}^{(k)}-\eta\bar{y}^{(k)}+\sqrt{2\eta}\bar{w}^{(k+1)},
\end{equation}
where for any $k=0,1,2,\ldots$,
\begin{equation}
\bar{y}^{(k+1)}=\bar{y}^{(k)}+\frac{1}{N}\sum_{i=1}^{N}\nabla f_{i}\left(x_{i}^{(k+1)}\right)-\frac{1}{N}\sum_{i=1}^{N}\nabla f_{i}\left(x_{i}^{(k)}\right)
+\bar{\xi}^{(k+2)}-\bar{\xi}^{(k+1)},
\end{equation}
which implies that for any $k=0,1,2,\ldots$,
\begin{equation}
\bar{y}^{(k)}=\frac{1}{N}\sum_{i=1}^{N}\nabla f_{i}\left(x_{i}^{(k)}\right)+\bar{\xi}^{(k+1)}.
\end{equation}
Therefore, we have
\begin{equation}
\bar{x}^{(k+1)}=\bar{x}^{(k)}-\eta\frac{1}{N}\sum_{i=1}^{N}\nabla f_{i}\left(x_{i}^{(k)}\right)-\eta\bar{\xi}^{(k+1)}+\sqrt{2\eta}\bar{w}^{(k+1)},
\end{equation}
which can be re-written as
\begin{equation}\label{bar:x:iterates}
\bar{x}^{(k+1)}=\bar{x}^{(k)}-\eta\frac{1}{N}\nabla f\left(\bar{x}^{(k)}\right)+\eta\mathcal{E}_{k}-\eta\bar{\xi}^{(k+1)}+\sqrt{2\eta}\bar{w}^{(k+1)},
\end{equation}
where
\begin{equation}\label{E:k:defn}
\mathcal{E}_{k}:=\frac{1}{N}\sum_{i=1}^{N}\left[\nabla f_{i}\left(\bar{x}^{(k)}\right)-\nabla f_{i}\left(x_{i}^{(k)}\right)\right].
\end{equation}
In the next lemma, we provide an explicit upper bound on 
the $L_{2}$ norm of the gradient error term $\mathcal{E}_{k}$.

\begin{lemma}\label{lem:2}
Assume that $\delta<\lambda^{B}<1$ with $\delta(k)$ defined in Assumption~\ref{assump:W}
and \eqref{assump:alpha:beta} holds.
Then, for any $k$, we have
\begin{equation*}
\mathbb{E}\left\Vert\mathcal{E}_{k}\right\Vert^{2}
\leq
\frac{3L^{2} \delta^{-2}}{N}\left(\delta^{\frac{2}{B}}\right)^{k}\mathbb{E}\left\Vert x^{(0)}\right\Vert^{2}
+\frac{3L^{2}D^{2}\eta^{2} \delta^{-2}}{N(1-\delta^{\frac{1}{B}})^{2}}
+\frac{6dL^{2}\eta\cdot  \delta^{-2}}{1-\delta^{\frac{2}{B}}}.
\end{equation*} 
\end{lemma}


The proof is given in Appendix~\ref{app-sec-lem:2}. Next, we recall from \eqref{overdamped:iterates:0} that the iterates $x_{k}$ are given by:
\begin{equation}\label{overdamped:iterates}
x_{k+1}=x_{k}-\eta\frac{1}{N}\nabla f(x_{k})+\sqrt{2\eta}\bar{w}^{(k+1)},
\end{equation}
where we take $x_{0}= \bar{x}_0=\frac{1}{N}\sum_{i=1}^{N}x_{i}^{(0)}$.  
This is a Euler-Maruyama discretization (with stepsize $\eta$)
of the continuous-time overdamped Langevin diffusion \eqref{overdamped:SDE}.  
Since the $L_{2}$ bound of the gradient error term $\mathcal{E}_{k}$
can be controlled as in Lemma~\ref{lem:2}, 
we will show that the average $\bar{x}^{(k)}$ and $x_{k}$ are close to each other
in $L_{2}$ distance. Indeed, we have the following estimate.

\begin{lemma}\label{lem:3}
Assume that $\delta<\lambda^{B}<1$ with $\delta(k)$ defined in Assumption~\ref{assump:W}
and \eqref{assump:alpha:beta} holds.
We also assume $\mathbb{E}\left\Vert x^{(0)}\right\Vert^{2}<\infty$, where the expectation is taken with respect to the distribution of the initial iterate $x^{(0)}$ at time \(k=0\).
For any stepsize $\eta \in (0,2/L)$, we have for every $k$,
\begin{align*}
\mathbb{E}\left\Vert\bar{x}^{(k)}-x_{k}\right\Vert^{2}
&\leq
\eta\left(\frac{\eta}{\mu(1-\frac{\eta L}{2})}+\frac{(1+\eta L)^{2}}{\mu^{2}(1-\frac{\eta L}{2})^{2}}\right)
\left(\frac{3L^{2}D^{2}\eta  \delta^{-2}}{N(1-\delta^{\frac{1}{B}})^{2}}
+\frac{6dL^{2}\cdot  \delta^{-2}}{1-\delta^{\frac{2}{B}}}\right)
+\frac{\eta\sigma^{2}}{\mu(1-\frac{\eta L}{2})N}
\nonumber
\\
&\qquad\qquad\qquad
+\mgb{\eta\left(\eta+\frac{(1+\eta L)^{2}}{\mu(1-\frac{\eta L}{2})}\right)
\mathsf H_k\left(1-\eta\mu\left(1-\frac{\eta L}{2}\right),\delta^{\frac{2}{B}}\right)}
\frac{3L^{2} \delta^{-2}}{N}\mathbb{E}\left\Vert x^{(0)}\right\Vert^{2}.
\end{align*}
\end{lemma}


As the proof of this lemma closely parallels that of the proof of Lemma~7 in \cite{gurbuzbalaban2021decentralized}, it is omitted in the main body here and instead presented in the supplementary material.

\subsubsection{$\mathcal{W}_{2}$ distance between the law of $x_{k}$ and the Gibbs distribution $\pi$}

The $\mathcal{W}_{2}$ distance between the Euler-Maruyama discretization
$x_{k}$ in \eqref{overdamped:iterates:0} of the overdamped Langevin diffusion \eqref{overdamped:SDE}
and the Gibbs distribution $\pi\propto e^{-f}$ has been established in the literature. 
\mgb{
The comparison chain in \eqref{overdamped:iterates:0} is ULA for the potential $f=\sum_i f_i$ with effective stepsize $h=\eta/N$. Under Assumption~\ref{assumption:f}, $f$ is $N\mu$-strongly convex and $NL$-smooth (equivalently, $f/N$ is $\mu$-strongly convex and $L$-smooth). Applying the exact-gradient part of Theorem~1 in \cite{dalalyan2019user} with $h=\eta/N$ gives the following statement.
\begin{lemma}[Rescaled from Theorem~1(a) in \cite{dalalyan2019user}]\label{lem:DK}
For $0<\eta\leq2/(\mu+L)$,
\[
\mathcal W_2(\operatorname{Law}(x_k),\pi)
\leq(1-\mu\eta)^k\mathcal W_2(\operatorname{Law}(x_0),\pi)
+\frac{1.65L}{\mu}\sqrt{\frac{\eta d}{N}}.
\]
\end{lemma}
The condition $\eta\leq1/((1+\alpha)L)$ with $\alpha>0$ in Theorem~\ref{thm:average} implies $\eta\leq2/(\mu+L)$ because $\mu\leq L$.}

Now, we are finally ready to prove Theorem~\ref{thm:average} and Theorem~\ref{thm:undirected}.

\subsubsection{Completing the Proofs of Theorem~\ref{thm:average} and Theorem~\ref{thm:undirected}}

\begin{proof}[Proof of Theorem~\ref{thm:average}]
The $L_{2}$ distance between the minimizer of $f$ and Gibbs distribution $\pi$ has been studied
in the literature; see e.g. \cite{gurbuzbalaban2021decentralized}.
More precisely, we have
\begin{equation}\label{eqn:bound:Gibbs}
\mathbb{E}_{X\sim\pi}\Vert X-x_{\ast}\Vert^{2}
\leq\frac{2dN^{-1}}{\mu}, 
\end{equation}
where $x_{\ast}$ is the unique minimizer of $f(x)$; see Lemma~11 in \cite{gurbuzbalaban2021decentralized}. 
Since $x_{0}=\frac{1}{N}\sum_{i=1}^{N}x_{i}^{(0)}$, 
we have $\mathbb{E}\Vert x_{0}\Vert^{2}<\infty$.
By \eqref{eqn:bound:Gibbs}, we get
\begin{align*}
\mathcal{W}_{2}\left(\mathrm{Law}(x_{0}),\pi\right)
\leq
\left(\mathbb{E}\Vert x_{0}-x_{\ast}\Vert^{2}\right)^{1/2}
+\left(\mathbb{E}_{X\sim\pi}\Vert X-x_{\ast}\Vert^{2}\right)^{1/2}
\leq\left(\mathbb{E}\Vert x_{0}-x_{\ast}\Vert^{2}\right)^{1/2}+\sqrt{2\mu^{-1}dN^{-1}}.
\end{align*}
It then follows from Lemma~\ref{lem:DK} that
for any $\eta\leq\frac{2}{\mu+L}$, we have
\begin{equation*}
\mathcal{W}_{2}\left(\mathrm{Law}(x_{k}),\pi\right)
\leq
(1-\mu\eta)^{k}\left(\left(\mathbb{E}\Vert x_{0}-x_{\ast}\Vert^{2}\right)^{1/2}+\sqrt{2\mu^{-1}dN^{-1}}\right)
+\frac{1.65L}{\mu}\sqrt{\eta dN^{-1}}.
\end{equation*}
Moreover, since \mgb{$\mathcal{W}_{2}\left(\mathrm{Law}\left(\bar{x}^{(k)}\right),\mathrm{Law}(x_{k})\right) \leq
\left(\mathbb{E}\left\Vert\bar{x}^{(k)}-x_{k}\right\Vert^{2}\right)^{1/2}$}, it follows from Lemma~\ref{lem:3} that
\mgb{\small \begin{align*}
&\mathcal{W}_{2}\left(\mathrm{Law}\left(\bar{x}^{(k)}\right),\mathrm{Law}(x_{k})\right)
\\
&\leq
\eta^{1/2}\left(\frac{\eta}{\mu(1-\frac{\eta L}{2})}+\frac{(1+\eta L)^{2}}{\mu^{2}(1-\frac{\eta L}{2})^{2}}\right)^{1/2}
\cdot
\left(\frac{3L^{2}D^{2}\eta  \delta^{-2}}{N(1-\delta^{\frac{1}{B}})^{2}}
+\frac{6dL^{2}\cdot  \delta^{-2}}{1-\delta^{\frac{2}{B}}}
\right)^{1/2}
\\
&\quad
+\frac{\sqrt{\eta}\sigma}{\sqrt{\mu(1-\frac{\eta L}{2})N}}
+\left[\eta\left(\eta+\frac{(1+\eta L)^2}{\mu(1-\eta L/2)}\right)
\mathsf H_k\!\left(1-\eta\mu(1-\eta L/2),\delta^{2/B}\right)\right]^{1/2}
\frac{\sqrt{3}L\delta^{-1}}{\sqrt{N}}
\left(\mathbb{E}\left\Vert x^{(0)}\right\Vert^{2}\right)^{1/2}.
\end{align*}}
The result then follows from the triangular inequality for the $2$-Wasserstein distance.
The proof is complete.~
\end{proof}

\begin{proof}[Proof of Theorem~\ref{thm:undirected}]
By the Cauchy-Schwarz inequality,
\begin{align}
\frac{1}{N}\sum_{i=1}^{N}\mathcal{W}_{2}\left(\mathrm{Law}\left(x_{i}^{(k)}\right),\mathrm{Law}\left(\bar{x}^{(k)}\right)\right)
&\leq\sqrt{\frac{1}{N}\sum_{i=1}^{N}\mathcal{W}_{2}^{2}\left(\mathrm{Law}\left(x_{i}^{(k)}\right),\mathrm{Law}\left(\bar{x}^{(k)}\right)\right)}\nonumber
\\
&\leq
\sqrt{\frac{1}{N}\sum_{i=1}^{N}\mathbb{E}\left\Vert x_{i}^{(k)}-\bar{x}^{(k)}\right\Vert^{2}}.\label{by:cauchy:schwarz}
\end{align}
By Lemma~\ref{cor:1}, we have
\begin{align}
\sqrt{\frac{1}{N}\sum_{i=1}^{N}\mathbb{E}\left\Vert x_{i}^{(k)}-\bar{x}^{(k)}\right\Vert^{2}}
&\leq
\left(\frac{3\cdot  \delta^{-2}\left(\delta^{\frac{2}{B}}\right)^{k}}{N}\mathbb{E}\left\Vert x^{(0)}\right\Vert^{2}
+\frac{3D^{2}\eta^{2} \delta^{-2}}{N(1-\delta^{\frac{1}{B}})^{2}}
+\frac{6d\eta\cdot  \delta^{-2}}{1-\delta^{\frac{2}{B}}}\right)^{1/2}
\\
&\leq
\frac{\sqrt{3} \delta^{-1}\delta^{\frac{k}{B}}}{\sqrt{N}}\left(\mathbb{E}\left\Vert x^{(0)}\right\Vert^{2}\right)^{1/2}
+\frac{\sqrt{3}D\eta  \delta^{-1}}{\sqrt{N}(1-\delta^{\frac{1}{B}})}
+\frac{\sqrt{6d\eta} \delta^{-1}}{\sqrt{1-\delta^{\frac{2}{B}}}}\mgb{=E_3}. \label{deriv-E3}
\end{align}
The result then follows from Theorem~\ref{thm:average} and the triangular inequality for the $2$-Wasserstein distance.
The proof is complete.
\end{proof}

\mgb{While Corollary~\ref{coro-param} provided iteration complexity results for the average Wasserstein distance over the nodes to the target, we can also get similar complexity results for the individual agents, as the following corollary shows.}
\mgb{\begin{corollary}[Uniform individual-agent consequence]\label{cor:individual-agent}
Under the assumptions of Theorem~\ref{thm:undirected},
\begin{equation}\label{eq:max-agent-bound}
\max_{1\leq i\leq N}\mathcal W_2\!\left(\operatorname{Law}(x_i^{(k)}),\pi\right)
\leq E_1(k,\eta)+E_2(k,\eta,\delta)+\sqrt N\,E_3(k,\eta,\delta).
\end{equation}
Consequently, if the parameter choice in Corollary~\ref{coro-param} is applied with target accuracy $\epsilon/\sqrt N$, then the left-hand side of \eqref{eq:max-agent-bound} is at most $\epsilon$. For fixed $N$, this preserves the $\mathcal O(\log(1/\epsilon)/\epsilon^2)$ dependence on $\epsilon$.
\end{corollary}
\begin{proof}
Couple $x_i^{(k)}$ and $\bar x^{(k)}$ on the probability space of the algorithm. Then, for every $1\leq i\leq N$,
\[
\mathcal W_2\!\left(\operatorname{Law}(x_i^{(k)}),\operatorname{Law}(\bar x^{(k)})\right)
\leq\left(\mathbb E\|x_i^{(k)}-\bar x^{(k)}\|^2\right)^{1/2}
\leq\left(\sum_{j=1}^N\mathbb E\|x_j^{(k)}-\bar x^{(k)}\|^2\right)^{1/2}
\leq\sqrt N\,E_3,
\]
where we used \eqref{deriv-E3}. The triangle inequality with Theorem~\ref{thm:average} gives \eqref{eq:max-agent-bound}. Finally,
$E_1+E_2+\sqrt N E_3\leq\sqrt N(E_1+E_2+E_3)$ for $N\geq1$, which proves the stated complexity consequence.
\end{proof}}
\section{Numerical Experiments}

In this section, we present numerical experiments evaluating the sampling performance of DIGing-SGLD, whose iterates are given in~\eqref{eqn:DIGing_algo_updates}--\eqref{eqn:DIGing_algo_updates_2}, in comparison to the DE-SGLD iterations described in~\cite{gurbuzbalaban2021decentralized}. All experiments are conducted under a constant stepsize and over time-varying network topologies. We evaluate both methods on Bayesian linear and logistic regression using synthetic datasets, and additionally on a real-world dataset~\cite{breast_cancer_wisconsin_diagnostic_17} for Bayesian logistic regression. 

For Bayesian linear regression, we examine how closely the samples generated by each agent approximate the target posterior distribution in the $\mathcal{W}_2$ distance, showing that for fixed $N$, each agent converges to a distribution lying within an $\epsilon$-neighborhood of the true posterior. For Bayesian logistic regression, we assess convergence in terms of classification accuracy on a held-out test set. Across both problems, we demonstrate that incorporating the DIGing-based gradient-tracking mechanism~\cite{nedic2017} in DIGing-SGLD effectively corrects the network-induced drift observed in DE-SGLD over time-varying networks. This advantage is particularly pronounced for Bayesian logistic regression in classification accuracy, although the $\mathcal{W}_2$ metric for Bayesian linear regression also exhibits a consistent performance gap between DIGing-SGLD and DE-SGLD.

\begin{figure}[t]
    \centering
    \begin{subfigure}[b]{0.45\textwidth}
        \centering
        \includegraphics[width=0.75\textwidth]{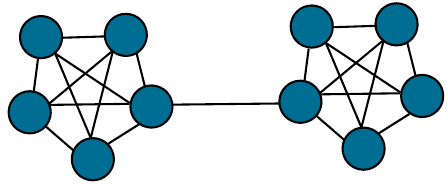}
        \caption{Example of a barbell topology for $N=10$ agents.}
        \label{fig:barbell_network}
    \end{subfigure}
    \hfill
    \begin{subfigure}[b]{0.45\textwidth}
        \centering
        \includegraphics[width=0.75\textwidth]{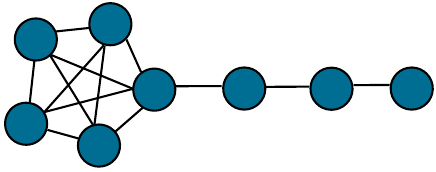}
        \caption{Example of a lollipop topology for $N=8$ agents.}
        \label{fig:lollipop_network}
    \end{subfigure}
    \caption{Illustrations of the two undirected time-varying network structures used in our experiments.}
    \label{fig:structured_network_structures}
\end{figure}

For the experimental setup, we consider two classes of undirected time-varying networks composed of $N$ agents. 
The first class is the barbell graph consisting of two cliques connected by a single edge, and the second is a lollipop graph 
$G_{{N^\prime},N^{\prime\prime}}$ comprising a clique 
$C_{N^\prime}$ with $N^\prime$ nodes and a path $P_{N^{\prime\prime}}$ of length ${N^{\prime\prime}}$ joined by a single edge to the path's terminal node.
Graphs such as barbell and lollipop graphs and their generalizations are widely used in the literature to evaluate decentralized algorithms, as they represent near–worst-case scenarios for information propagation over a network \cite{can2022randomized,aybat2017decentralized,PhysRevE.90.012816}. 
For the barbell topology, we use $N \in \{20, 30\}$ depending on the problem. The time-varying nature is introduced by first constructing two undirected complete graphs, each with $N/2$ agents. At every algorithmic iteration $k$, two random agents, one from each complete graph, are connected to form a single connected network. Figure~\ref{fig:barbell_network} illustrates an example of this topology for $N=10$. \looseness=-1
For the time-varying lollipop network, we use $N=20$ agents partitioned into two subsets, $N = N^\prime + N^{\prime\prime}$, with $N^\prime = 8$, and $N^{\prime\prime} = 12$. The set with $N^\prime$ agents forms a complete subgraph (clique). At each iteration $k$, we randomly shuffle the agents in $N^{\prime\prime}$ and construct a path graph over them. We then select a random agent from the clique and connect it to the terminal node of the path, resulting in a lollipop graph. Figure~\ref{fig:lollipop_network} shows an example of this topology for $N=8$, where $N'=5$ forms a clique and $N^{\prime\prime}= 3$ forms a path.
\looseness=-1

After generating the network topology in each iteration for both classes of time-varying networks, we construct the corresponding weight matrix $W^{(k)}$ using the Metropolis constant edge-weight rule~\cite{boyd2004fastest,xiao2007distributed}, defined as
\begin{align}
W_{ij}^{(k)} =
\begin{cases}
\frac{1}{\max\left\{\deg_i^{(k)}, \deg_j^{(k)}\right\} + \widehat{\varepsilon}}, & \text{if } (i,j) \in \mathcal{E}(k) \mgb{\text{ and } i\neq j}, \\[3pt]
0, & \text{if } (i,j) \notin \mathcal{E}(k) \text{ and } i \neq j, \\[3pt]
1 - \sum_{\ell \in \Omega_{i}(k) \mgb{\setminus\{i\}}} W_{i\ell}^{(k)}, & \text{if } i = j,
\end{cases}
\end{align}
where $\deg_i^{(k)}$ denotes the degree of agent $i$ at iteration $k$, and $\widehat{\varepsilon} > 0$ (set to $10^{-6}$ in our experiments) ensures that the Markov chain represented by $W^{(k)}$ is aperiodic. In all experiments, the time-varying networks are generated such that each topology repeats every $50$ algorithmic iterations in a cyclic manner. While this repetition is an approximation, it can be interpreted as controlling the parameter $B$ in our theoretical model, with larger repetition intervals corresponding to higher values of $B$.

\subsection{Bayesian Linear Regression}

We now present the results of DIGing-SGLD and DE-SGLD for Bayesian linear regression over the time-varying barbell and lollipop topologies. In linear regression, each data sample consists of a response variable $y\in\mathbb{R}$ and a feature vector $\widehat{z}\in\mathbb{R}^d$. To include an intercept term, we use the augmented feature vector $z=[\,\widehat{z}^\top\ \ 1\,]^\top\in\mathbb{R}^{d+1}$. In our experiments, we generate $n=100$ samples $\{(y_i,z_i)\}_{i=1}^n$ with $d=5$, and distribute them uniformly across $N=20$ agents, so that each agent holds $\bar{n}:=n/N=5$ samples.

To enable reporting of performance in terms of the $\mathcal{W}_2$ distance, we work with synthetic data generated under a Gaussian linear model. The samples are generated according to $y_i = z_i^\top x + \delta_i$, where the true parameter $x\in\mathbb{R}^{d+1}$ is drawn once from $x\sim\mathcal{N}(0,\lambda^{-1}I_{d+1})$ with $\lambda=0.1$ and then held fixed. The noise variables satisfy $\delta_i\sim\mathcal{N}(0,\sigma^2)$ with $\sigma^2=1$, and the underlying feature vectors are drawn independently as $\widehat{z}_i\sim\mathcal{N}(0,I_d)$. Under this model, the posterior distribution for $x$ given all $n$ samples is $\pi(x)\propto\exp\bigl(-\sum_{j=1}^n f_j(x)\bigr)$, where the local function at agent $j$ is the strongly convex and smooth function
\begin{align}\label{eqn:linear_reg_func_node}
    f_j(x):=\frac{\|Z_j x - y_j\|^2}{2} + \frac{\lambda\|x\|^2}{2N},
\end{align}
with $Z_j\in\mathbb{R}^{\bar{n}\times(d+1)}$ denoting the matrix of local feature vectors and $y_j\in\mathbb{R}^{\bar{n}}$ the corresponding responses. The posterior in this case also admits a closed-form expression, given by $\pi(x) = \mathcal{N}(\widehat{m},\widehat{\Sigma})$, where $\widehat{\Sigma}=(Z^\top Z+\lambda I_{d+1})^{-1}$ and $\widehat{m}=\widehat{\Sigma}Z^\top y$ for the full dataset matrix $Z\in\mathbb{R}^{n\times(d+1)}$ and full response vector $y\in\mathbb{R}^n$.

\begin{figure}[t]
    \centering
    \begin{subfigure}[b]{0.31\textwidth}
        \centering
        \includegraphics[width=\textwidth]{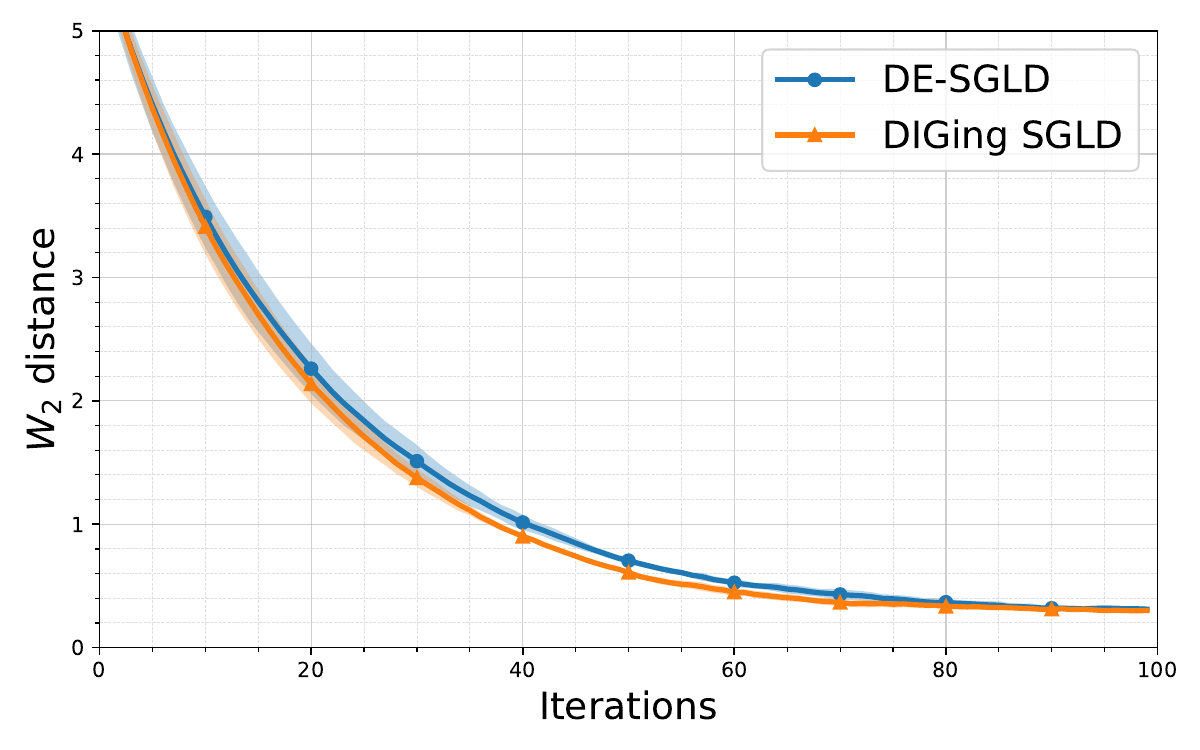}
        \caption{Time-varying barbell network with $N=20$ agents; full-batch setup ($b=5$).}
        \label{fig:LR_Fig_Experiment_LR5}
    \end{subfigure}
    \hfill
    \begin{subfigure}[b]{0.31\textwidth}
        \centering
        \includegraphics[width=\textwidth]{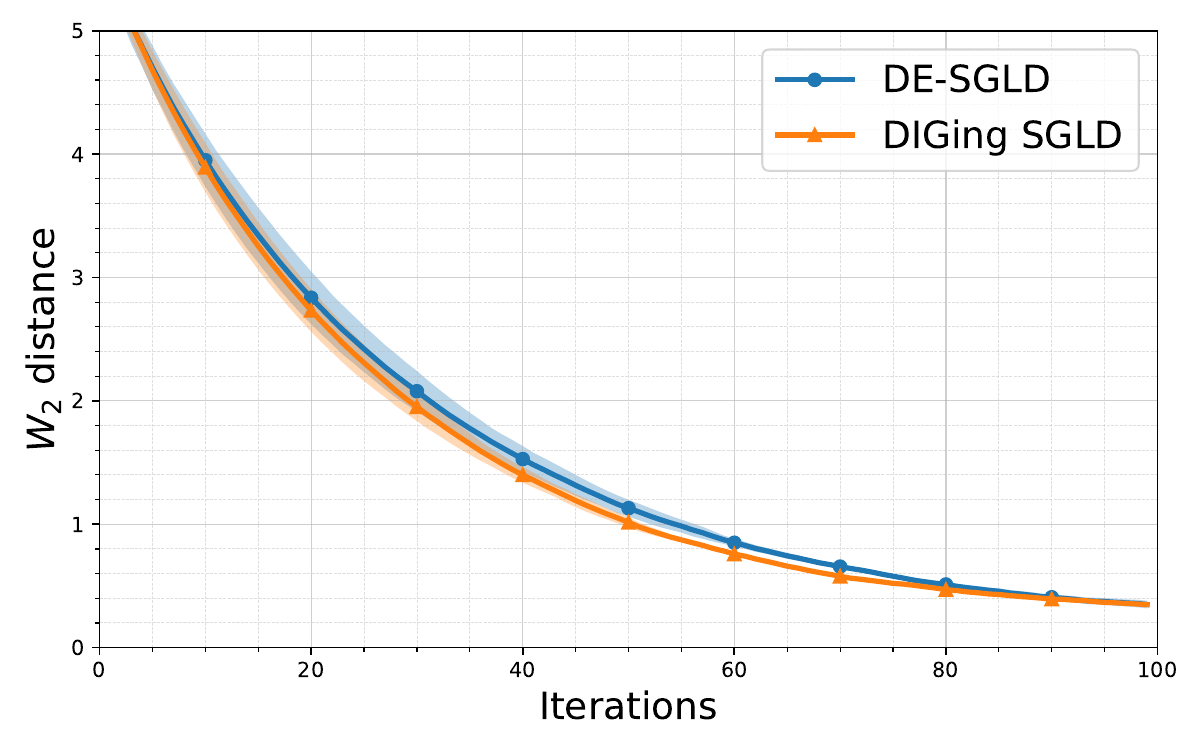}
        \caption{Time-varying barbell network with $N=20$ agents; mini-batch setup ($b=3$ out of $\bar{n}=5$).}
        \label{fig:LR_Fig_Experiment_LR5a}
    \end{subfigure}
    \hfill
    \begin{subfigure}[b]{0.31\textwidth}
        \centering
        \includegraphics[width=\textwidth]{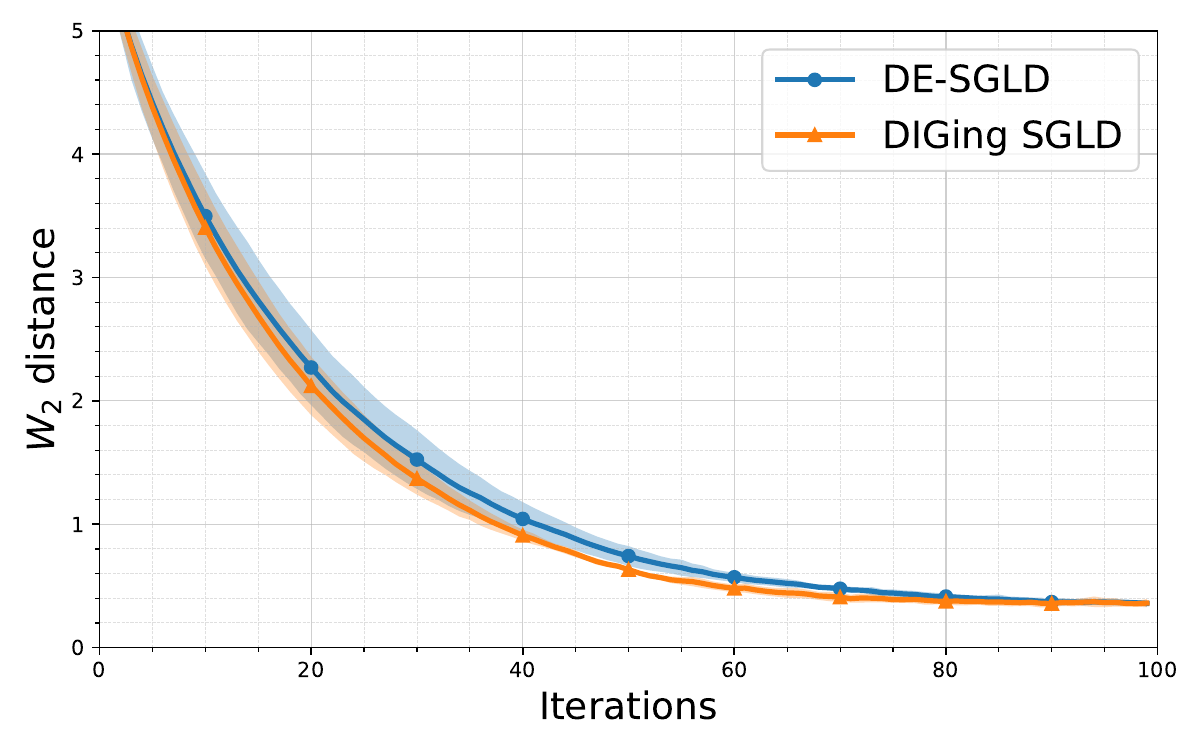}
        \caption{Time-varying lollipop network with $N=20$ agents; full-batch setup ($b=5$).}
        \label{fig:LR_Fig_Experiment_LR6}
    \end{subfigure}
    \caption{Comparison of DIGing-SGLD and DE-SGLD for Bayesian linear regression on synthetic data under time-varying barbell and lollipop network structures. Each plot displays the average Wasserstein distance across agents, with one standard deviation shown as a shaded region.}
    \label{fig:linear_regression_performance}
    \vspace{-0.15in}
\end{figure}

We conduct three experiments using this synthetic-data construction for both DIGing-SGLD and DE-SGLD. In the first experiment (Figure~\ref{fig:LR_Fig_Experiment_LR5}), each agent uses its full batch of $b=\bar{n}=5$ samples to compute the gradient of~\eqref{eqn:linear_reg_func_node} at every iteration, and the DIGing-SGLD updates are run over a time-varying barbell network. The second experiment (Figure~\ref{fig:LR_Fig_Experiment_LR5a}) uses the same data-generation procedure and barbell topology, but each agent now samples a mini-batch of $b=3$ out of its $\bar{n}=5$ local samples at each iteration. Since this results in a stochastic gradient of $f_j(x)$, and we require this gradient to be an unbiased estimate of the local gradient, we scale the mini-batch gradient by the factor $\bar{n}/b$. In the third experiment (Figure~\ref{fig:LR_Fig_Experiment_LR6}), we evaluate DIGing-SGLD and DE-SGLD over a time-varying lollipop network, with each agent again using its full batch of $b=\bar{n}=5$ samples.
\mgb{Here, we note that when full batches are used, the exact gradient update for linear regression is affine in the iterate. Consequently, the DIGing-SGLD iterates $x^{(k)}$ remain Gaussian for every $k$, and Assumption~\ref{assumption:noise} holds with $\sigma=0$. In contrast, under mini-batching, the conditional variance of the stochastic gradient estimator generally grows quadratically with $\|x\|$, because the empirical Hessian associated with a mini-batch varies across batches. Hence, the global uniform variance bound in~\eqref{eq:conditional-gradient-noise} of Assumption~\ref{assumption:noise} need not hold unless additional mechanisms or assumptions are imposed, such as projection, gradient clipping, or a moment-dependent noise analysis. Accordingly, the mini-batch experiment should be viewed as an empirical stress test that lies outside the scope of the present theoretical guarantees.}

In all three experiments, both methods are run independently for $200$ trials with different initializations. These independent runs are used to estimate the mean vector $m_j^{(k)}$ and covariance matrix $\Sigma_j^{(k)}$ of the approximate posterior at each agent $j$ and iteration $k$. Using the closed-form expression for the $\mathcal{W}_2$ distance between Gaussian distributions~\cite{givens1984class}, we \mgb{estimate} the gap between the empirical and true posteriors via
\[
\mgb{
\mathcal W_{2,j}^{(k)}=
\left(
\|\widehat m-m_j^{(k)}\|^2
+\operatorname{Tr}(\Sigma_j^{(k)})
+\operatorname{Tr}(\widehat\Sigma)
-2\operatorname{Tr}\!\left[
\left(
(\Sigma_j^{(k)})^{1/2}\widehat\Sigma(\Sigma_j^{(k)})^{1/2}
\right)^{1/2}
\right]
\right)^{1/2}.
}
\]
\mgb{This is the exact $\mathcal W_2$ distance between the target Gaussian and the moment-matched Gaussian $\mathcal{N}(m_j^{(k)}, \Sigma_j^{(k)})$}. Figure~\ref{fig:linear_regression_performance} reports the results across the three experimental setups. Each plot shows, for both DIGing-SGLD and DE-SGLD, the average Wasserstein distance over all agents, with one standard deviation of the agent-wise distances indicated as a shaded region. In each case, we fix the number of iterations to $100$ and tune the stepsize to minimize the $\mathcal{W}_2$ error at iteration $100$.

Recall that both Figure~\ref{fig:LR_Fig_Experiment_LR5} and ~\ref{fig:LR_Fig_Experiment_LR5a} concern the same barbell topology; but the former is based on exact gradients whereas the latter uses stochastic gradients with batch size $b=3$. As expected, the stochastic noise in the gradients results in slower convergence. 
Comparing Figure~\ref{fig:LR_Fig_Experiment_LR5} and Figure~\ref{fig:LR_Fig_Experiment_LR6}, we observe similar convergence behavior for both lollipop and barbell graphs.
Overall, these results illustrate that DIGing-SGLD outperforms DE-SGLD across both network structures. Although the improvement margin is modest, we should note that DE-SGLD lacks convergence guarantees for time-varying networks, further underscoring the importance of the guarantees established here for DIGing-SGLD.


\subsection{Bayesian Logistic Regression}

We next evaluate DIGing-SGLD and DE-SGLD for Bayesian logistic regression under time-varying barbell and lollipop topologies, using both synthetic datasets ($N=20$ agents) and a real-world dataset ($N=30$ agents; barbell topology only). In Bayesian logistic regression, each data sample consists of a binary class label $y\in\{0,1\}$ and a feature vector $\widehat{z}\in\mathbb{R}^d$. As in the linear regression setting, we work with the augmented feature vector $z=[\,\widehat{z}^\top\ \ 1\,]^\top\in\mathbb{R}^{d+1}$ to incorporate an intercept term. Given $z$ and a parameter vector (model) $x\in\mathbb{R}^{d+1}$, the likelihood of the class label is
\[
\mathbb{P}(y=1 \mid z,x)=\frac{1}{1+\exp(-z^\top x)}=\sigma\left(z^\top x\right),
\]
and we adopt the Gaussian prior for the model $x$: $p(x)=\mathcal{N}(0,\lambda^{-1}I_{d+1})$ with $\lambda>0$.

Given $n$ labeled training samples $\{(y_i,z_i)\}_{i=1}^n$, the posterior of $x$ is of the form $\pi(x)\propto\exp(-f(x))$. When the data are evenly distributed across the $N$ agents, so that each agent holds $\bar{n}=n/N$ samples, the global function $f(x)$ decomposes as $f(x)=\sum_{j=1}^N f_j(x)$, where each local function $f_j$ is smooth and strongly convex and is given by
\begin{align}
f_j(x)=\sum_{i=1}^{\bar{n}} \!\left[-\log\!\left(1-\sigma\left(z_{i,j}^\top x
\right)\right) + y_{i,j}\log\!\left(\frac{1-\sigma(z_{i,j}^\top x)}{\sigma(z_{i,j}^\top x)}\right)\right] + \frac{\lambda}{2N}\|x\|^2,
\end{align}
with $(z_{i,j},y_{i,j})$ denoting the $i$th sample at agent $j$.

For the Bayesian logistic regression problem, to the best of our knowledge, the $\mathcal{W}_2$ distance of the distribution of the nodes to the posterior $\pi(x)$ does not admit a closed-form expression, unlike the Bayesian linear regression setting. Even though the $\mathcal{W}_2$ distance between the histograms of two distributions can be computed using approaches such as \cite{flamary2021pot}, this would be in general expensive to compute up to high accuracy in moderate to high dimensions. Therefore, in contrast to the Bayesian linear regression setting, we compare DIGing-SGLD and DE-SGLD using classification accuracy on a separate test dataset, defined as the proportion of correctly predicted labels. As in the linear regression experiments, the accuracy at each agent and iteration is computed by averaging the per-agent accuracy over multiple independent runs, and stepsizes are hand-tuned to optimize accuracy at the final iteration (here, iteration $100$).

\subsubsection{Synthetic Data Experiments}
For the synthetic experiments, the parameter $x\in\mathbb{R}^{d+1}$ is drawn once from the prior $p(x)=\mathcal{N}(0,\lambda^{-1}I_{d+1})$ with $\lambda=0.1$, and feature vectors are generated by sampling $\widehat{z}_i\sim\mathcal{N}(0,I_d)$. Each class label is assigned by drawing 
$p_i$ from the uniform distribution on $[0,1]$, and setting $y_i=1$ if $p_i\leq\sigma(z_i^\top x)$ and $y_i=0$ otherwise. We set $d=5$ and generate $600$ samples, which are then split into training and test sets using a $70$--$30$ train--test ratio. This yields $420$ training samples, and we discard $40$ to obtain a class-balanced training set of size $n=380$. These samples are distributed evenly across $N=20$ agents, giving each agent $\bar{n}=19$ samples.

\begin{figure}[t]
    \centering
    \begin{subfigure}[b]{0.31\textwidth}
        \centering
        \includegraphics[width=\textwidth]{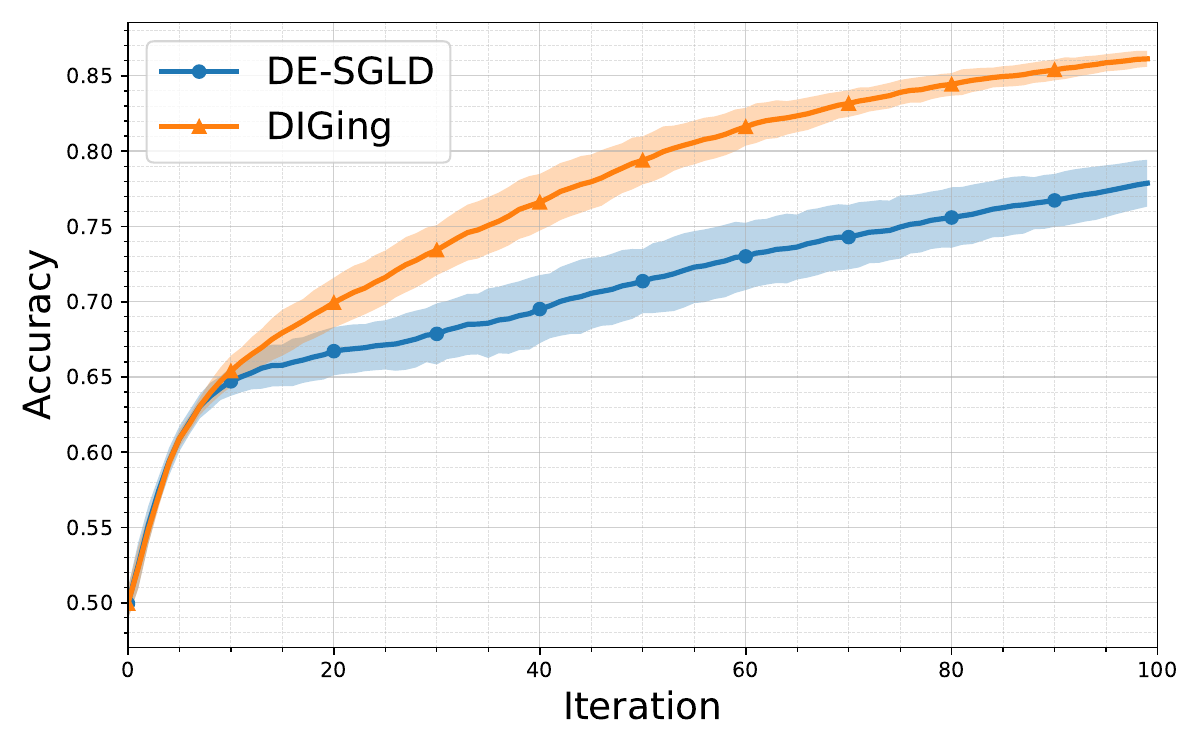}
        \caption{Time-varying barbell network ($N=20$); synthetic dataset; mini-batch setup ($b=1$ out of $\bar{n}=19$).}
        \label{fig:LO_Fig_Experiment_LO9}
    \end{subfigure}
    \hfill
    \begin{subfigure}[b]{0.31\textwidth}
        \centering
        \includegraphics[width=\textwidth]{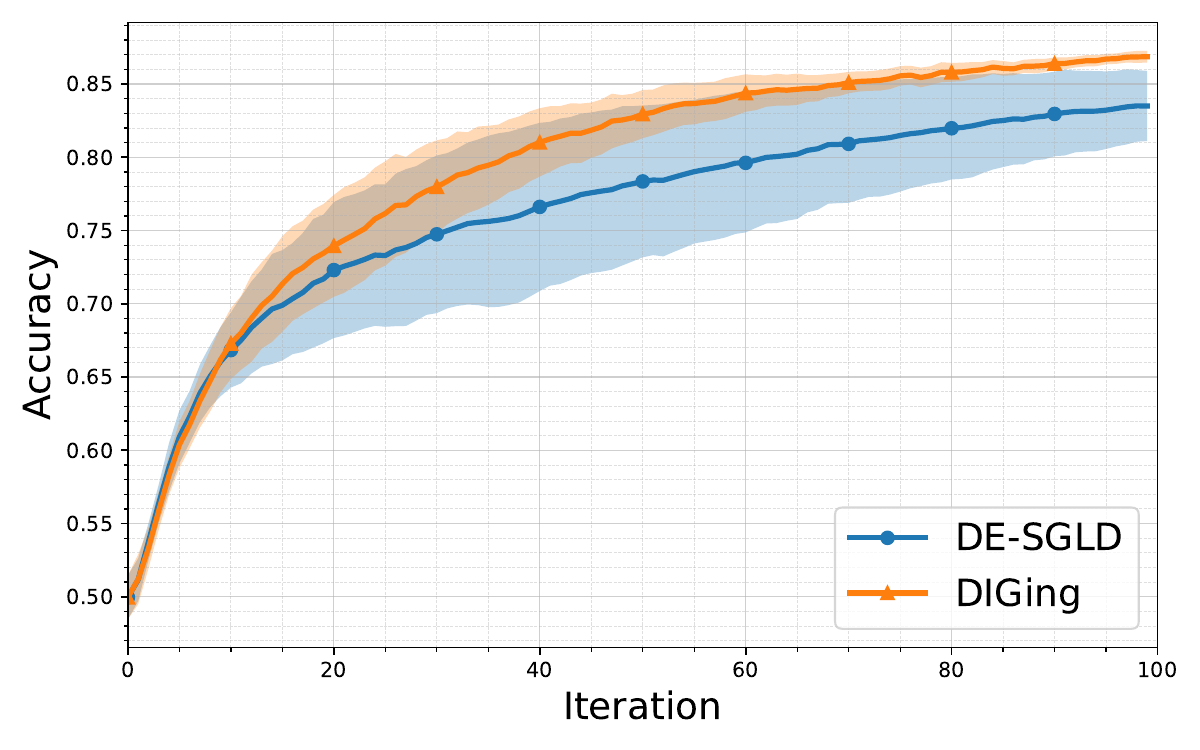}
        \caption{Time-varying lollipop network ($N=20$); synthetic dataset; mini-batch setup ($b=1$ out of $\bar{n}=19$).}
        \label{fig:LO_Fig_Experiment_LO10}
    \end{subfigure}
    \hfill
    \begin{subfigure}[b]{0.31\textwidth}
        \centering
        \includegraphics[width=\textwidth]{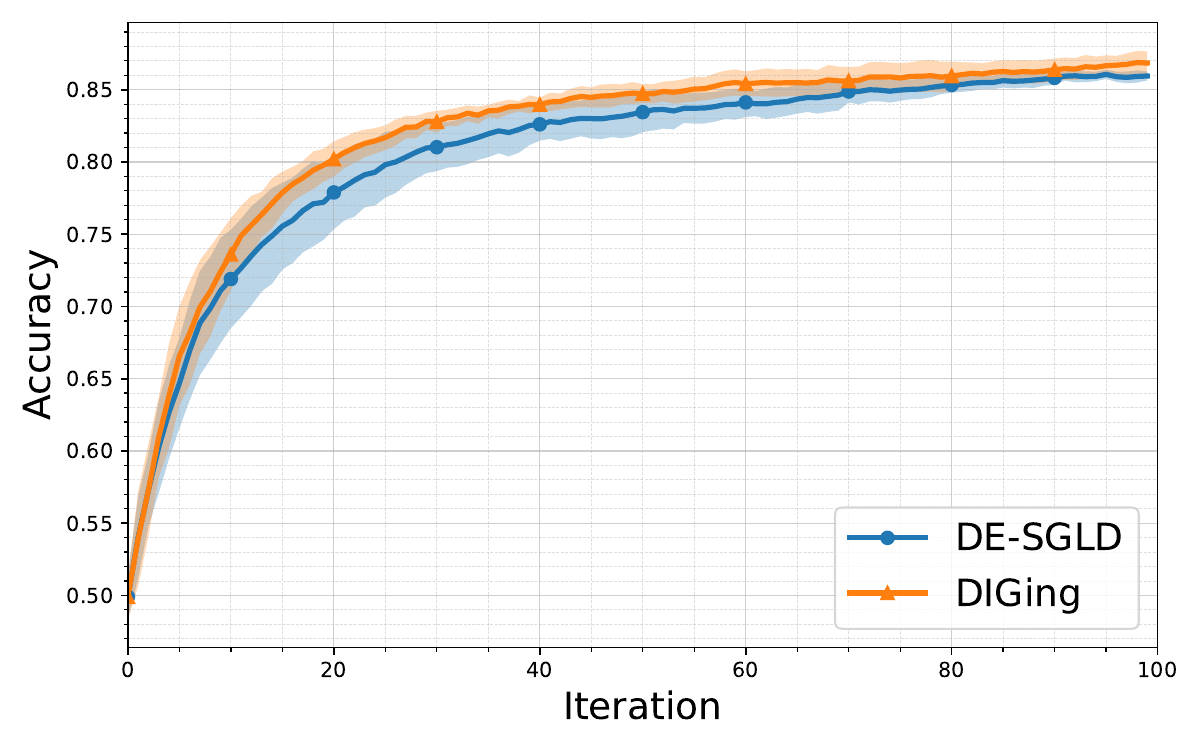}
        \caption{Time-varying barbell network ($N=30$); real dataset; full-batch setup ($b=\bar{n}=1$).}
        \label{fig:LO_Fig_Experiment_LO4a}
    \end{subfigure}
    \hfill
    \caption{Performance comparison of DIGing-SGLD and DE-SGLD for Bayesian logistic regression under time-varying barbell and lollipop network structures. Each plot displays the average classification accuracy across agents and independent trials, with one standard deviation across agents shown as a shaded region.}
    \label{fig:logistic_regression_performance}
    \vspace{-0.15in}
\end{figure}

We conduct two experiments under this setup, one using a time-varying barbell topology and one using a time-varying lollipop topology. In both experiments, at each iteration, every agent uses a mini-batch of $b=1$ sample to compute a stochastic gradient of the logistic loss. Since this yields a stochastic gradient of $f_j(x)$, we scale the gradient by $\bar{n}/b = 19$ to obtain an unbiased local gradient estimator. Each experiment is repeated independently $200$ times; accuracy is evaluated for each estimate $x_j^{(k)}$ at each agent $j$, and the plots in Figure~\ref{fig:logistic_regression_performance} report the average accuracy computed across both agents and independent trials, with one standard deviation across agents shown as a shaded region in the corresponding plots. 

The first synthetic experiment uses the time-varying barbell topology (Figure~\ref{fig:LO_Fig_Experiment_LO9}); the second uses the time-varying lollipop topology (Figure~\ref{fig:LO_Fig_Experiment_LO10}). In both cases, DIGing-SGLD consistently outperforms DE-SGLD, achieving higher accuracy throughout the training horizon. The performance gap is more pronounced for the barbell topology, which has a smaller spectral gap and therefore induces a more severe information-flow bottleneck than the lollipop topology. These results parallel our findings for Bayesian linear regression: DIGing-SGLD maintains stable performance in the time-varying network setting, whereas DE-SGLD exhibits larger variability and reduced accuracy. We again note that DE-SGLD does not admit any convergence guarantees for time-varying networks, and its weaker empirical performance in these experiments further highlights the robustness of DIGing-SGLD in dynamic network environments.

\subsubsection{Real Data Experiments}
We finally compare DIGing-SGLD and DE-SGLD on Bayesian logistic regression using real data under a time-varying barbell topology with $N=30$ agents. We use the UCI ML Breast Cancer Wisconsin (Diagnostic) dataset~\cite{breast_cancer_wisconsin_diagnostic_17}, which contains $569$ samples and $d=30$ features derived from digitized images of fine needle aspirate biopsies. Class labels are encoded as $1$ for benign and $0$ for malignant. We split the dataset using a $10$--$90$ train–test ratio, yielding $51$ training samples. After standardizing and class-balancing the training data, we obtain $n=30$ samples that are then evenly distributed across the $30$ agents, so that each agent has $\bar{n}=1$ sample. At each iteration, each agent processes its single sample ($b=\bar{n}=1$) to compute a gradient. For a fixed iteration budget, both the stepsize and the regularization parameter $\lambda$ are tuned, with the final choice of the regularization parameter being $\lambda=0.3$.

Results for this experiment are shown in Figure~\ref{fig:LO_Fig_Experiment_LO4a}. The curves display the average accuracy computed across both independent trials (200 repetitions) and agents, with one standard deviation across agents shown as a shaded region.  Consistent with our synthetic experiments, DIGing-SGLD converges reliably in the time-varying setting, while DE-SGLD converges more slowly and, importantly, has no known convergence guarantees for time-varying networks.

\section{Conclusion} 
\mgb{In this work, we introduced \emph{DIGing--SGLD}, a decentralized Langevin-based sampling algorithm that operates over time-varying networks. The method integrates distributed inexact gradient tracking (DIGing)---originating in decentralized optimization---into stochastic gradient Langevin dynamics, thereby mitigating the network-induced discrepancy between the local gradient information available to individual agents and the aggregate gradient required by the target distribution. Under strong convexity and smoothness assumptions on the component functions $f_i$ and a connectivity condition on the underlying time-varying graph, we establish finite-time, non-asymptotic guarantees in the 2-Wasserstein distance with explicit constants. In particular, the bounds show that each agent's marginal distribution approaches an $\mathcal{O}(\sqrt{\eta})$ neighborhood of the target distribution
$
\pi(x)\propto e^{-f(x)}$ where
$
f(x)=\sum_{i=1}^{n} f_i(x).
$. Thus, for a fixed positive stepsize, discretization and stochastic-gradient effects remain, while the residual neighborhood can be made arbitrarily small by decreasing $\eta$. For fixed problem and network parameters, choosing a stepsize $\eta=\mathcal{O}(\epsilon^2)$ yields an
$
\mathcal{O}\!\left(\frac{\log(1/\epsilon)}{\epsilon^2}\right)
$
iteration bound to reach $\epsilon$-accuracy, matching the dependence on $\epsilon$ of the existing constant-stepsize centralized and static-graph SGLD iteration bounds while extending the analysis to dynamic, decentralized, coordinator-free networks. Numerical experiments on Bayesian linear and logistic regression illustrate the behavior of DIGing--SGLD under changing topologies and stochastic gradient updates. Looking forward, promising research directions include extending DIGing-SGLD to non-convex objectives or directed network settings, thereby strengthening its theoretical guarantees and enhancing its applicability to large-scale decentralized Bayesian inference.}

\appendix


\section*{Acknowledgments} We gratefully acknowledge Rutgers University’s Department of Electrical and Computer Engineering, Rutgers Business School and Florida State University's Department of Mathematics for fostering the collaborative environment that enabled this work. This work was supported in part by the Office of Naval Research under grant N00014-24-1-2628. In some parts of the text, it is acknowledged that AI tools were used to edit or polish the authors’ written text for spelling and grammar.

\section{Proofs of Technical Lemmas}\label{proofs:technical}

\subsection{Proof of Lemma~\ref{lem:0}}\label{app-sec-proof-lem:0}

\begin{proof}
It follows from \eqref{tilde:y:upper:bound} and \eqref{q:upper:bound} that 
\begin{align*}
\Vert\tilde{y}\Vert_{L_{2}}^{\lambda,K}
&\leq\frac{\gamma_{1}\gamma_{2}\gamma_{3}\omega_{4}(K)+\gamma_{1}\gamma_{2}\omega_{3}(K)
+\gamma_{1}\omega_{2}(K)+\omega_{1}(K)}{1-\gamma_{1}\gamma_{2}\gamma_{3}\gamma_{4}}
\\
&=\frac{\gamma_{1}\gamma_{2}\gamma_{3}\tilde{\omega}_{4}+\gamma_{1}\gamma_{2}\tilde{\omega}_{3}
+\tilde{\omega}_{1}}{1-\gamma_{1}\gamma_{2}\gamma_{3}\gamma_{4}}
+\frac{\gamma_{1}\gamma_{2}\gamma_{3}\hat{\omega}_{4}+\gamma_{1}\gamma_{2}\hat{\omega}_{3}
+\hat{\omega}_{1}}{1-\gamma_{1}\gamma_{2}\gamma_{3}\gamma_{4}}\frac{1}{\lambda^{K}},
\end{align*}
and
\begin{align*}
\Vert q\Vert_{L_{2}}^{\lambda,K}
&\leq\frac{\gamma_{3}\gamma_{4}\gamma_{1}\omega_{2}(K)+\gamma_{3}\gamma_{4}\omega_{1}(K)
+\gamma_{3}\omega_{4}(K)+\omega_{3}(K)}{1-\gamma_{1}\gamma_{2}\gamma_{3}\gamma_{4}}
\\
&=\frac{\gamma_{3}\gamma_{4}\tilde{\omega}_{1}+\gamma_{3}\tilde{\omega}_{4}
+\tilde{\omega}_{3}}{1-\gamma_{1}\gamma_{2}\gamma_{3}\gamma_{4}}
+\frac{\gamma_{3}\gamma_{4}\hat{\omega}_{1}+\gamma_{3}\hat{\omega}_{4}
+\hat{\omega}_{3}}{1-\gamma_{1}\gamma_{2}\gamma_{3}\gamma_{4}}\frac{1}{\lambda^{K}},
\end{align*}
where $\gamma_{1},\gamma_{2},\gamma_{3},\gamma_{4}$ are defined in \eqref{defn:gamma:1:2:main}-\eqref{defn:gamma:3:4:main} and $\omega_{1}(K),\omega_{2}(K),\omega_{3}(K),\omega_{4}(K)$ are defined in \eqref{defn:omega:1:2:K}, \eqref{defn:omega:3:K} and \eqref{defn:omega:4:K} \mgb{where we used $\omega_{2}(K)=0$}. We also recall from \eqref{defn:tilde:omega:1:main}, \eqref{defn:tilde:omega:3:main} and \eqref{defn:tilde:omega:4:main} that
\begin{align}
&\tilde{\omega}_{1}:=\frac{\lambda^{B}}{\lambda^{B}-\delta}
\sum_{t=1}^{B}\lambda^{1-t}\left\Vert\tilde{y}^{(t-1)}\right\Vert_{L_{2}},
\qquad
\hat{\omega}_{1}:=\frac{\lambda^{B}}{\lambda^{B}-\delta}\cdot 2B\sigma\sqrt{N},
\\
&\tilde{\omega}_{3}:=\sqrt{N}\mgb{\left(\mathbb{E}\left\Vert\bar{x}^{(0)}-x_{\ast}\right\Vert^2\right)^{1/2}}
,
\qquad
\hat{\omega}_{3}:=\frac{\sqrt{N}}{\lambda}\left(\sqrt{\frac{L(1+\alpha)}{\mu\alpha}+\beta}\right)\frac{1}{\mu}\left(\sigma+\sqrt{\frac{2d}{\eta}}\right),
\\
&\tilde{\omega}_{4}:=\frac{\lambda^{B}}{\lambda^{B}-\delta}
\sum_{t=1}^{B}\lambda^{1-t}\left\Vert\tilde{x}^{(t-1)}\right\Vert_{L_{2}},
\qquad
\hat{\omega}_{4}:=\frac{\lambda^{B}}{\lambda^{B}-\delta}\cdot\mgb{B}\sqrt{2\eta Nd}.
\end{align}
Hence, for every $k$, we have
\begin{align*}
\left\Vert\tilde{y}^{(k)}\right\Vert_{L_{2}}
&\leq\frac{\gamma_{1}\gamma_{2}\gamma_{3}\tilde{\omega}_{4}+\gamma_{1}\gamma_{2}\tilde{\omega}_{3}
+\tilde{\omega}_{1}}{1-\gamma_{1}\gamma_{2}\gamma_{3}\gamma_{4}}\lambda^{k}
+\frac{\gamma_{1}\gamma_{2}\gamma_{3}\hat{\omega}_{4}+\gamma_{1}\gamma_{2}\hat{\omega}_{3}
+\hat{\omega}_{1}}{1-\gamma_{1}\gamma_{2}\gamma_{3}\gamma_{4}}
\\
&\leq\frac{\gamma_{1}\gamma_{2}\gamma_{3}\tilde{\omega}_{4}+\gamma_{1}\gamma_{2}\tilde{\omega}_{3}
+\tilde{\omega}_{1}}{1-\gamma_{1}\gamma_{2}\gamma_{3}\gamma_{4}}
+\frac{\gamma_{1}\gamma_{2}\gamma_{3}\hat{\omega}_{4}+\gamma_{1}\gamma_{2}\hat{\omega}_{3}
+\hat{\omega}_{1}}{1-\gamma_{1}\gamma_{2}\gamma_{3}\gamma_{4}},
\end{align*}
and
\begin{align*}
\left\Vert q^{(k)}\right\Vert_{L_{2}}
&\leq\frac{\gamma_{3}\gamma_{4}\tilde{\omega}_{1}+\gamma_{3}\tilde{\omega}_{4}
+\tilde{\omega}_{3}}{1-\gamma_{1}\gamma_{2}\gamma_{3}\gamma_{4}}\lambda^{k}
+\frac{\gamma_{3}\gamma_{4}\hat{\omega}_{1}+\gamma_{3}\hat{\omega}_{4}
+\hat{\omega}_{3}}{1-\gamma_{1}\gamma_{2}\gamma_{3}\gamma_{4}}
\\
&\leq\frac{\gamma_{3}\gamma_{4}\tilde{\omega}_{1}+\gamma_{3}\tilde{\omega}_{4}
+\tilde{\omega}_{3}}{1-\gamma_{1}\gamma_{2}\gamma_{3}\gamma_{4}}
+\frac{\gamma_{3}\gamma_{4}\hat{\omega}_{1}+\gamma_{3}\hat{\omega}_{4}
+\hat{\omega}_{3}}{1-\gamma_{1}\gamma_{2}\gamma_{3}\gamma_{4}},
\end{align*}
where we used $0<\lambda<1$.
Next, we can compute that
\begin{align*}
\mathbb{E}\left\Vert y^{(k)}\right\Vert^{2}
\leq
2\mathbb{E}\left\Vert\tilde{y}^{(k)}\right\Vert^{2}
+2\mathbb{E}\left\Vert\mathbf{\bar{y}}^{(k)}\right\Vert^{2}
=2\mathbb{E}\left\Vert\tilde{y}^{(k)}\right\Vert^{2}
+2N\mathbb{E}\left\Vert\bar{y}^{(k)}\right\Vert^{2}.
\end{align*}
Moreover, 
\begin{align*}
2N\mathbb{E}\left\Vert\bar{y}^{(k)}\right\Vert^{2}
&=2N\mathbb{E}\left\Vert\frac{1}{N}\sum_{i=1}^{N}\nabla f_{i}\left(x_{i}^{(k)}\right)+\bar{\xi}^{(k+1)}\right\Vert^{2}
\\
&=2N\mathbb{E}\left\Vert\frac{1}{N}\sum_{i=1}^{N}\left(\nabla f_{i}\left(x_{i}^{(k)}\right)-\mgb{\nabla} f_{i}(x_{\ast})\right)+\bar{\xi}^{(k+1)}\right\Vert^{2}
\\
&\leq
4N\mathbb{E}\left\Vert\frac{1}{N}\sum_{i=1}^{N}\left(\nabla f_{i}\left(x_{i}^{(k)}\right)-\mgb{\nabla} f_{i}(x_{\ast})\right)\right\Vert^{2}
+4N\mathbb{E}\left\Vert\bar{\xi}^{(k+1)}\right\Vert^{2}
\\
&\leq
\mgb{4L^2}
\mathbb{E}\sum_{i=1}^{N}\left\Vert x_{i}^{(k)}-x_{\ast}\right\Vert^{2}+\mgb{4}
\sigma^{2}
=\mgb{4L^2}
\mathbb{E}\left\Vert q^{(k)}\right\Vert^{2}+\mgb{4}
\sigma^{2}.
\end{align*}
Hence, we conclude that 
\begin{align}
\mathbb{E}\left\Vert y^{(k)}\right\Vert^{2}
&\leq
2\mathbb{E}\left\Vert\tilde{y}^{(k)}\right\Vert^{2}
+\mgb{4L^2}
\mathbb{E}\left\Vert q^{(k)}\right\Vert^{2}+\mgb{4}
\sigma^{2}
\nonumber
\\
&\leq
2\left(\frac{\gamma_{1}\gamma_{2}\gamma_{3}(\tilde{\omega}_{4}+\hat{\omega}_{4})+\gamma_{1}\gamma_{2}(\tilde{\omega}_{3}+\hat{\omega}_{3})
+\tilde{\omega}_{1}+\hat{\omega}_{1}}{1-\gamma_{1}\gamma_{2}\gamma_{3}\gamma_{4}}\right)^{2}
\nonumber
\\
&\qquad\qquad\qquad
+\mgb{4L^2}
\left(\frac{\gamma_{3}\gamma_{4}(\tilde{\omega}_{1}+\hat{\omega}_{1})+\gamma_{3}(\tilde{\omega}_{4}+\hat{\omega}_{4})
+\tilde{\omega}_{3}+\hat{\omega}_{3}}{1-\gamma_{1}\gamma_{2}\gamma_{3}\gamma_{4}}\right)^{2}+\mgb{4}
\sigma^{2}.
\end{align}
This completes the proof.
\end{proof}


\subsection{Proof of Lemma~\ref{lem:first:arrow}}\label{appendix-sec-lem:first:arrow}

\begin{proof}
The proof can be directly adapted from the proof of Lemma 3.9 in \cite{nedic2017}
by replacing the matrix notation in \cite{nedic2017} by the vector notation in our paper. 
Since $\nabla F(x)$ is $L$-Lipschitz, we have
\begin{align}
\left\Vert\nabla F\left(x^{(k+1)}\right)-\nabla F\left(x^{(k)}\right)\right\Vert
\leq L\left\Vert x^{(k+1)}-x^{(k)}\right\Vert
\leq L\left\Vert x^{(k+1)}-\mathbf{x}_{\ast}\right\Vert+L\left\Vert x^{(k)}-\mathbf{x}_{\ast}\right\Vert.
\label{F:L:ineq}
\end{align}
By the definition of $z$ and $q$, it follows from \eqref{F:L:ineq} that
\begin{equation}\label{z:q:two:sides}
\lambda^{-(k+1)}\left\Vert z^{(k+1)}\right\Vert
\leq
L\lambda^{-(k+1)}\left\Vert q^{(k+1)}\right\Vert
+\frac{L}{\lambda}\lambda^{-k}\left\Vert q^{(k)}\right\Vert.
\end{equation}
By taking the maximum on both hand sides of \eqref{z:q:two:sides} over $k=0,1,\ldots,K-1$, we conclude that
\begin{equation}
\Vert z\Vert_{L_{2}}^{\lambda,K}
\leq L\Vert q\Vert_{L_{2}}^{\lambda,K}+\frac{L}{\lambda}\Vert q\Vert_{L_{2}}^{\lambda,K-1}
\leq L\left(1+\frac{1}{\lambda}\right)\Vert q\Vert_{L_{2}}^{\lambda,K}.    
\end{equation}
This completes the proof.
\end{proof}


\subsection{Proof of Lemma~\ref{lem:second:arrow}}\label{sec:lem:second:arrow}

\begin{proof}
First, we recall from \eqref{y:vector:form} and \eqref{q:z:k:defn} that
\begin{equation}
y^{(k+1)}=\mathcal{W}^{(k)}y^{(k)}+z^{(k+1)}+\xi^{(k+2)}-\xi^{(k+1)}.\label{ineq-y-recursion}
\end{equation}
Therefore, for any $k\geq B-1$, we have
\begin{align*}
&\left\Vert\tilde{y}^{(k+1)}\right\Vert_{L_{2}}
=\left\Vert \mathcal{L}_{N}y^{(k+1)}\right\Vert_{L_{2}}
\\
&\leq\left\Vert \mathcal{L}_{N}\mathcal{W}_{B}^{(k)}y^{(k+1-B)}\right\Vert_{L_{2}}
+\left\Vert \mathcal{L}_{N}\mathcal{W}_{B-1}^{(k)}z^{(k+2-B)}\right\Vert_{L_{2}}+\cdots+\left\Vert \mathcal{L}_{N}\mathcal{W}_{1}^{(k)}z^{(k)}\right\Vert_{L_{2}}+\left\Vert \mathcal{L}_{N}\mathcal{W}_{0}^{(k)}z^{(k+1)}\right\Vert_{L_{2}}
\\
&\quad\quad
+\left\Vert \mathcal{L}_{N}\mathcal{W}_{B-1}^{(k)}\xi^{(k+3-B)}\right\Vert_{L_{2}}+\cdots+\left\Vert \mathcal{L}_{N}\mathcal{W}_{1}^{(k)}\xi^{(k+1)}\right\Vert_{L_{2}}+\left\Vert \mathcal{L}_{N}\mathcal{W}_{0}^{(k)}\xi^{(k+2)}\right\Vert_{L_{2}}
\\
&\quad\quad\quad
+\left\Vert \mathcal{L}_{N}\mathcal{W}_{B-1}^{(k)}\xi^{(k+2-B)}\right\Vert_{L_{2}}+\cdots+\left\Vert \mathcal{L}_{N}\mathcal{W}_{1}^{(k)}\xi^{(k)}\right\Vert_{L_{2}}+\left\Vert \mathcal{L}_{N}\mathcal{W}_{0}^{(k)}\xi^{(k+1)}\right\Vert_{L_{2}},
\end{align*}
where $\mathcal{L}_{N}:=I_{Nd}-\frac{1}{N}\left(\left(1_{N}1_{N}^{\top}\right)\otimes I_{d}\right)$. 
By applying Lemma~\ref{lem:iterate} in Appendix~\ref{app:additional.lemmas} and Assumption~\ref{assumption:noise}, we get
\begin{equation}
\left\Vert\tilde{y}^{(k+1)}\right\Vert_{L_{2}}
\leq
\delta\left\Vert\tilde{y}^{(k+1-B)}\right\Vert_{L_{2}}+\sum_{t=1}^{B}\left\Vert z^{(k+2-t)}\right\Vert_{L_{2}}
+2B\sigma\sqrt{N}.
\end{equation}
Therefore, for any $k=B-1,B,\ldots,$ we have
\begin{equation}
\lambda^{-(k+1)}\left\Vert\tilde{y}^{(k+1)}\right\Vert_{L_{2}}
\leq
\frac{\delta}{\lambda^{B}}\lambda^{-(k+1-B)}\left\Vert\tilde{y}^{(k+1-B)}\right\Vert_{L_{2}}
+\sum_{t=1}^{B}\frac{1}{\lambda^{t-1}}\lambda^{-(k+2-t)}\left\Vert z^{(k+2-t)}\right\Vert_{L_{2}}
+2B\sigma\sqrt{N}\mgb{\lambda^{-(k+1)}}.\label{ineq-intermediate}
\end{equation} 

By following a similar argument as in the proof of Lemma 3.10 in \cite{nedic2017}, we obtain
that for every $K$:
\begin{equation}
\Vert\tilde{y}\Vert_{L_{2}}^{\lambda,K}
\leq
\frac{\delta}{\lambda^{B}}\Vert\tilde{y}\Vert_{L_{2}}^{\lambda,K}
+\sum_{t=1}^{B}\frac{1}{\lambda^{t-1}}\Vert z\Vert_{L_{2}}^{\lambda,K}
+\frac{2B\sigma\sqrt{N}}{\lambda^{K}}
+\sum_{t=1}^{B}\lambda^{1-t}\left\Vert\tilde{y}^{(t-1)}\right\Vert_{L_{2}}.
\end{equation}
This implies that
\begin{equation}
\Vert\tilde{y}\Vert_{L_{2}}^{\lambda,K}
\leq
\frac{\lambda(1-\lambda^{B})}{(\lambda^{B}-\delta)(1-\lambda)}\Vert z\Vert_{L_{2}}^{\lambda,K}
+\frac{\lambda^{B}}{\lambda^{B}-\delta}\frac{2B\sigma\sqrt{N}}{\lambda^{K}}
+\frac{\lambda^{B}}{\lambda^{B}-\delta}\sum_{t=1}^{B}\lambda^{1-t}\left\Vert\tilde{y}^{(t-1)}\right\Vert_{L_{2}}.
\end{equation}
This completes the proof.
\end{proof}

\subsection{Proof of Lemma~\ref{lem:third:arrow}}
\label{appendix-proof-lem:third:arrow}
\begin{proof}
Recall from \eqref{x:vector:form} that
$x^{(k+1)}=\mathcal{W}^{(k)}x^{(k)}-\eta y^{(k)}+\sqrt{2\eta}w^{(k+1)}.$
This recursion is structurally similar to the one in \eqref{ineq-y-recursion}, which appears in the proof of Lemma~\ref{lem:second:arrow}. The proof therefore proceeds by following the same steps as in the proof of Lemma~\ref{lem:second:arrow}
and the fact that 
$\mathbb{E}\left\Vert\sqrt{2\eta}w^{(k+1)}\right\Vert^{2}=2\eta Nd$.
\end{proof}


\subsection{Proof of Lemma~\ref{lem:last:arrow}}
\label{app-section-lem:last:arrow}
\begin{proof}
First, let us recall from \eqref{q:z:k:defn} that
$q^{(k)}=x^{(k)}-\mathbf{x}_{\ast}=x^{(k)}-\mathbf{\bar{x}}^{(k)}+\mathbf{\bar{x}}^{(k)}-\mathbf{x}_{\ast}$.
Next, by taking the average of $N$ nodes in \eqref{DIGing:x}--\eqref{DIGing:y}, and using the fact that $W^{(k)}$ is doubly stochastic, we obtain:
\begin{equation}
\bar{x}^{(k+1)}=\bar{x}^{(k)}-\eta\bar{y}^{(k)}+\sqrt{2\eta}\bar{w}^{(k+1)},
\end{equation}
where for any $k=0,1,2,\ldots$,
\begin{equation}
\bar{y}^{(k+1)}=\bar{y}^{(k)}+\frac{1}{N}\sum_{i=1}^{N}\nabla f_{i}\left(x_{i}^{(k+1)}\right)-\frac{1}{N}\sum_{i=1}^{N}\nabla f_{i}\left(x_{i}^{(k)}\right)
+\bar{\xi}^{(k+2)}-\bar{\xi}^{(k+1)},
\end{equation}
which implies that for any $k=0,1,2,\ldots$,
\begin{equation}
\bar{y}^{(k)}=\frac{1}{N}\sum_{i=1}^{N}\nabla f_{i}\left(x_{i}^{(k)}\right)+\bar{\xi}^{(k+1)}.
\end{equation}
Therefore, we have
\begin{equation}\label{to:be:rewritten}
\bar{x}^{(k+1)}=\bar{x}^{(k)}-\eta\frac{1}{N}\sum_{i=1}^{N}\nabla f_{i}\left(x_{i}^{(k)}\right)-\eta\bar{\xi}^{(k+1)}+\sqrt{2\eta}\bar{w}^{(k+1)}.
\end{equation}
By Lemma~\ref{lemma-surjective-gradients} in Appendix~\ref{app:additional.lemmas}, we can re-write the equation \eqref{to:be:rewritten} as:
\begin{equation}
\bar{x}^{(k+1)}=\bar{x}^{(k)}-\eta\frac{1}{N}\sum_{i=1}^{N}\nabla f_{i}\left(s_{i}^{(k)}\right),
\end{equation}
where $s_{i}^{(k)}$ is defined implicitly via:
\begin{equation}
\nabla f_{i}\left(s_{i}^{(k)}\right)=\nabla f_{i}\left(x_{i}^{(k)}\right)+\bar{\xi}^{(k+1)}-\sqrt{\frac{2}{\eta}}\bar{w}^{(k+1)}.
\end{equation}
Since $f_{i}$ is $\mu$-strongly convex, we have
$\left\Vert s_{i}^{(k)}-x_{i}^{(k)}\right\Vert\leq\frac{1}{\mu}\left\Vert\bar{\xi}^{(k+1)}-\sqrt{\frac{2}{\eta}}\bar{w}^{(k+1)}\right\Vert$,
which implies that
\begin{equation}
\left\Vert s_{i}^{(k)}-x_{i}^{(k)}\right\Vert_{L_{2}}
\leq\frac{1}{\mu}\left(\frac{\sigma}{\sqrt{N}}+\sqrt{\frac{2d}{\eta N}}\right).
\end{equation}
By applying Lemma~\ref{lem:inexact} in Appendix~\ref{app:additional.lemmas}, under the assumption that
$\sqrt{1-\frac{\eta\mu\beta}{\beta+1}}\leq\lambda<1$
and $\eta\leq\frac{1}{(1+\alpha)L}$,
where $\alpha,\beta>0$, we have
\begin{equation}
\Vert\bar{x}-x_{\ast}\Vert_{L_{2}}^{\lambda,K}
\leq
\mgb{\left(\mathbb{E}\left\Vert\bar{x}^{(0)}-x_{\ast}\right\Vert^2\right)^{1/2}}
+\left(\lambda\sqrt{N}\right)^{-1}\left(\sqrt{\frac{L(1+\alpha)}{\mu\alpha}+\beta}\right)
\sum_{i=1}^{N}\Vert\bar{x}-s_{i}\Vert_{L_{2}}^{\lambda,K},
\end{equation}
for any $K=0,1,2,\ldots$ where $x_{\ast}$ is the minimizer of $f$ and 
$\Vert \bar{x}-s_{i} \Vert_{L_{2}}^{\lambda,K}
=\max_{0,1,\ldots,K}\frac{1}{\lambda^{k}}\left(\mathbb{E}\left\Vert \bar{x}-s_{i}^{(k)}\right\Vert^{2}\right)^{1/2}.$
Therefore, we have 
\begin{eqnarray*}
\Vert\bar{x}-x_{\ast}\Vert_{L_{2}}^{\lambda,K}
&\leq&
\mgb{\left(\mathbb{E}\left\Vert\bar{x}^{(0)}-x_{\ast}\right\Vert^2\right)^{1/2}}
+\left(\lambda\sqrt{N}\right)^{-1}\left(\sqrt{\frac{L(1+\alpha)}{\mu\alpha}+\beta}\right)
\sum_{i=1}^{N}\Vert\bar{x}-x_{i}\Vert_{L_{2}}^{\lambda,K}
\\
&&\qquad\qquad\qquad\qquad
+\left(\lambda\sqrt{N}\right)^{-1}\left(\sqrt{\frac{L(1+\alpha)}{\mu\alpha}+\beta}\right)\frac{1}{\mu}\left(\frac{\sigma}{\sqrt{N}}+\sqrt{\frac{2d}{\eta N}}\right)\frac{N}{\lambda^{K}}
\\
&\leq&
\mgb{\left(\mathbb{E}\left\Vert\bar{x}^{(0)}-x_{\ast}\right\Vert^2\right)^{1/2}}+(\lambda)^{-1}\left(\sqrt{\frac{L(1+\alpha)}{\mu\alpha}+\beta}\right){\sqrt{N}}
\Vert\tilde{x}\Vert_{L_{2}}^{\lambda,K}\\
&&\qquad\qquad\qquad\qquad+(\lambda)^{-1}\left(\sqrt{\frac{L(1+\alpha)}{\mu\alpha}+\beta}\right)\frac{1}{\mu}\left(\sigma+\sqrt{\frac{2d}{\eta}}\right)\frac{1}{\lambda^{K}},
\end{eqnarray*}
where we used \begin{equation} \Vert \bar{x}-x_{i} \Vert_{L_{2}}^{\lambda,K} \leq \|\tilde{x}\|_{L_2}^{\lambda,K}.\label{ineq-to-prove-last}\end{equation}
{To see why the last inequality holds, note that
\begin{equation}\label{eq-sum-identity}
  \sum_{i=1}^{N} \mathbb{E}\big\Vert \bar{x}^{(k)}-x_i^{(k)}\big\Vert^{2} =\mathbb{E}\big\Vert \tilde{x}^{(k)}\big\Vert^{2},
  \qquad k=0,1,\ldots,K .
\end{equation}
}
{Fix $i\in\{1,\ldots,N\}$ and $k\in\{0,1,\ldots,K\}$. 
Then, 
\begin{align}
  \big\Vert \bar{x}-x_{i}\big\Vert_{L_{2}}^{\lambda,K}
  &=\max_{0\leq k\leq K}\frac{1}{\lambda^{k}}
    \Big(\mathbb{E}\big\Vert \bar{x}^{(k)}-x_{i}^{(k)}\big\Vert^{2}\Big)^{1/2}
  \leq\;\,\max_{0\leq k\leq K}\frac{1}{\lambda^{k}}
    \Big(\mathbb{E}\big\Vert \tilde{x}^{(k)}\big\Vert^{2}\Big)^{1/2}
  \;=\;\,\big\Vert \tilde{x}\big\Vert_{L_{2}}^{\lambda,K},
\end{align}
 which holds for every $i=1,\ldots,N$. This implies \eqref{ineq-to-prove-last}.}

Finally,
$q^{(k)}=\tilde{x}^{(k)}+\mathbf{\bar{x}}^{(k)}-\mathbf{x}_{\ast}$,
and it follows that
\begin{equation}
\left\Vert q\right\Vert_{L_{2}}^{\lambda,K}
\leq\left\Vert\tilde{x}\right\Vert_{L_{2}}^{\lambda,K}+\sqrt{N}\Vert\bar{x}-x_{\ast}\Vert_{L_{2}}^{\lambda,K}.
\end{equation}
Hence, we conclude that
\begin{align*}
\left\Vert q\right\Vert_{L_{2}}^{\lambda,K}
&\leq
\sqrt{N} \mgb{\left(\mathbb{E}\left\Vert\bar{x}^{(0)}-x_{\ast}\right\Vert^2\right)^{1/2}}
+\left(1+\frac{{N}}{\lambda}\left(\sqrt{\frac{L(1+\alpha)}{\mu\alpha}+\beta}\right)\right)
\Vert\tilde{x}\Vert_{L_{2}}^{\lambda,K}
\\
&\qquad\qquad
+\frac{\sqrt{N}}{\lambda}\left(\sqrt{\frac{L(1+\alpha)}{\mu\alpha}+\beta}\right)\frac{1}{\mu}\left(\sigma+\sqrt{\frac{2d}{\eta}}\right)\frac{1}{\lambda^{K}}.
\end{align*}
This completes the proof.
\end{proof}

\subsection{Proof of Lemma~\ref{lem:y:q}}\label{app-sec-proof-lem:y:q}

\begin{proof}
Under the assumption $\delta<\lambda^{B}<1$ with $\delta(k)$ defined in Assumption~\ref{assump:W}
and \eqref{assump:alpha:beta}, Lemma~\ref{lem:first:arrow}, Lemma~\ref{lem:second:arrow}, Lemma~\ref{lem:third:arrow}
and Lemma~\ref{lem:last:arrow} hold,
and it follows from \eqref{recursive:ineq:1}--\eqref{recursive:ineq:4} that
\begin{equation}
\Vert\tilde{y}\Vert_{L_{2}}^{\lambda,K}\leq\gamma_{1}\gamma_{2}\gamma_{3}\gamma_{4}\Vert\tilde{y}\Vert_{L_{2}}^{\lambda,K}
+\gamma_{1}\gamma_{2}\gamma_{3}\omega_{4}(K)+\gamma_{1}\gamma_{2}\omega_{3}(K)
+\gamma_{1}\omega_{2}(K)+\omega_{1}(K),
\end{equation}
and if $0<\gamma_{1}\gamma_{2}\gamma_{3}\gamma_{4}<1$, we obtain:
\begin{equation}
\Vert\tilde{y}\Vert_{L_{2}}^{\lambda,K}\leq\frac{\gamma_{1}\gamma_{2}\gamma_{3}\omega_{4}(K)+\gamma_{1}\gamma_{2}\omega_{3}(K)
+\gamma_{1}\omega_{2}(K)+\omega_{1}(K)}{1-\gamma_{1}\gamma_{2}\gamma_{3}\gamma_{4}}.
\end{equation}
Similarly, one can show that
\begin{equation}
\Vert q\Vert_{L_{2}}^{\lambda,K}\leq\gamma_{1}\gamma_{2}\gamma_{3}\gamma_{4}\Vert q\Vert_{L_{2}}^{\lambda,K}
+\gamma_{3}\gamma_{4}\gamma_{1}\omega_{2}(K)+\gamma_{3}\gamma_{4}\omega_{1}(K)
+\gamma_{3}\omega_{4}(K)+\omega_{3}(K),
\end{equation}
and if $0<\gamma_{1}\gamma_{2}\gamma_{3}\gamma_{4}<1$, we obtain:
\begin{equation}
\Vert q\Vert_{L_{2}}^{\lambda,K}\leq\frac{\gamma_{3}\gamma_{4}\gamma_{1}\omega_{2}(K)+\gamma_{3}\gamma_{4}\omega_{1}(K)
+\gamma_{3}\omega_{4}(K)+\omega_{3}(K)}{1-\gamma_{1}\gamma_{2}\gamma_{3}\gamma_{4}}.
\end{equation}
This completes the proof.
\end{proof}

\subsection{Proof of Lemma~\ref{lem:1}}\label{app-sec-lem:1}

\begin{proof}
By the iterates of $x^{(k)}$ given in \eqref{x:vector:form}, we get
\begin{equation*}
x^{(k+1)}=\left(W^{(k)}\otimes I_{d}\right)x^{(k)}-\eta y^{(k)}+\sqrt{2\eta}w^{(k+1)}.
\end{equation*}
It follows that for any $k\geq 1$,
\begin{align}
x^{(k)}=\left(W_{k}^{(k-1)}\otimes I_{d}\right)x^{(0)}-\eta\sum_{s=0}^{k-1}\left(W_{k-1-s}^{(k-1)}\otimes I_{d}\right)y^{(s)}
+\sqrt{2\eta}\sum_{s=0}^{k-1}\left(W_{k-1-s}^{(k-1)}\otimes I_{d}\right)w^{(s+1)}.\label{x:k:expression}
\end{align}
Let us define $\mathbf{\bar{x}}^{(k)}:=\left[\left(\bar{x}^{(k)}\right)^{\top},\cdots,\left(\bar{x}^{(k)}\right)^{\top}\right]^{\top}\in\mathbb{R}^{Nd}$ where $\bar{x}^{(k)}=\frac{1}{N}\sum_{i=1}^{N}x_{i}^{(k)}$.
Notice that
$\mathbf{\bar{x}}^{(k)}=\frac{1}{N}\left(\left(1_{N}1_{N}^{\top}\right)\otimes I_{d}\right)x^{(k)}$.
Therefore, we get
\begin{equation*}
\sum_{i=1}^{N}\left\Vert x_{i}^{(k)}-\bar{x}^{(k)}\right\Vert^{2}
=\left\Vert x^{(k)}-\mathbf{\bar{x}}^{(k)}\right\Vert^{2}
=\left\Vert x^{(k)}-\frac{1}{N}\left(\left(1_{N}1_{N}^{\top}\right)\otimes I_{d}\right)x^{(k)}\right\Vert^{2}.
\end{equation*}
Note that it follows from \eqref{x:k:expression} that
\begin{align*}
&x^{(k)}-\frac{1}{N}\left(\left(1_{N}1_{N}^{\top}\right)\otimes I_{d}\right)x^{(k)}
\\
&=\left(W_{k}^{(k-1)}\otimes I_{d}\right)x^{(0)}-\frac{1}{N}\left(\left(1_{N}1_{N}^{\top}W_{k}^{(k-1)}\right)\otimes I_{d}\right)x^{(0)}
\\
&\quad-\eta\sum_{s=0}^{k-1}\left(W_{k-1-s}^{(k-1)}\otimes I_{d}\right)y^{(s)}
+\eta\sum_{s=0}^{k-1}\frac{1}{N}\left(\left(1_{N}1_{N}^{\top}W_{k-1-s}^{(k-1)}\right)\otimes I_{d}\right)y^{(s)}
\\
&\quad\quad
+\sqrt{2\eta}\sum_{s=0}^{k-1}\left(W_{k-1-s}^{(k-1)}\otimes I_{d}\right)w^{(s+1)}
-\sqrt{2\eta}\sum_{s=0}^{k-1}\frac{1}{N}\left(\left(1_{N}1_{N}^{\top}W_{k-1-s}^{(k-1)}\right)\otimes I_{d}\right)w^{(s+1)}.
\end{align*}
By Cauchy-Schwarz inequality, we have
\begin{align*}
&\left\Vert x^{(k)}-\frac{1}{N}\left(\left(1_{N}1_{N}^{\top}\right)\otimes I_{d}\right)x^{(k)}\right\Vert^{2}
\\
&\leq 
3\left\Vert\left(W_{k}^{(k-1)}\otimes I_{d}\right)x^{(0)}-\frac{1}{N}\left(\left(1_{N}1_{N}^{\top}W_{k}^{(k-1)}\right)\otimes I_{d}\right)x^{(0)}\right\Vert^{2}
\\
&\quad
+3\left\Vert-\eta\sum_{s=0}^{k-1}\left(W_{k-1-s}^{(k-1)}\otimes I_{d}\right)y^{(s)}
+\eta\sum_{s=0}^{k-1}\frac{1}{N}\left(\left(1_{N}1_{N}^{\top}W_{k-1-s}^{(k-1)}\right)\otimes I_{d}\right)y^{(s)}\right\Vert^{2}
\\
&\quad
+3\left\Vert\sqrt{2\eta}\sum_{s=0}^{k-1}\left(W_{k-1-s}^{(k-1)}\otimes I_{d}\right)w^{(s+1)}
-\sqrt{2\eta}\sum_{s=0}^{k-1}\frac{1}{N}\left(\left(1_{N}1_{N}^{\top}W_{k-1-s}^{(k-1)}\right)\otimes I_{d}\right)w^{(s+1)}\right\Vert^{2}
\\
&=3\left\Vert\left(W_{k}^{(k-1)}\otimes I_{d}\right)x^{(0)}-\frac{1}{N}\left(\left(1_{N}1_{N}^{\top}\right)\otimes I_{d}\right)x^{(0)}\right\Vert^{2}
\\
&\quad
+3\left\Vert-\eta\sum_{s=0}^{k-1}\left(W_{k-1-s}^{(k-1)}\otimes I_{d}\right)y^{(s)}
+\eta\sum_{s=0}^{k-1}\frac{1}{N}\left(\left(1_{N}1_{N}^{\top}\right)\otimes I_{d}\right)y^{(s)}\right\Vert^{2}
\\
&\quad
+3\left\Vert\sqrt{2\eta}\sum_{s=0}^{k-1}\left(W_{k-1-s}^{(k-1)}\otimes I_{d}\right)w^{(s+1)}
-\sqrt{2\eta}\sum_{s=0}^{k-1}\frac{1}{N}\left(\left(1_{N}1_{N}^{\top}\right)\otimes I_{d}\right)w^{(s+1)}\right\Vert^{2},
\end{align*}
where we used the property that $W^{(k)}$ is doubly stochastic for every $k$.
Therefore, we get
\begin{align} \label{eq:inter-estimate}
\left\Vert x^{(k)}-\frac{1}{N}\left(\left(1_{N}1_{N}^{\top}\right)\otimes I_{d}\right)x^{(k)}\right\Vert^{2}
&\leq 3\left\Vert\left(\left(W_{k}^{(k-1)}-\frac{1}{N}1_{N}1_{N}^{\top}\right)\otimes I_{d}\right)x^{(0)}\right\Vert^{2}
\nonumber \\
&+3\eta^{2}\left\Vert\sum_{s=0}^{k-1}\left(\left(W_{k-1-s}^{(k-1)}-\frac{1}{N}1_{N}1_{N}^{\top}\right)\otimes I_{d}\right)y^{(s)}\right\Vert^{2}
\nonumber \\
&\qquad
+6\eta\left\Vert\sum_{s=0}^{k-1}\left(\left(W_{k-1-s}^{(k-1)}-\frac{1}{N}1_{N}1_{N}^{\top}\right)\otimes I_{d}\right)w^{(s+1)}\right\Vert^{2}. \nonumber\\
\end{align}
Note that
\begin{align*}
3\eta^{2}\left\Vert\sum_{s=0}^{k-1}\left(\left(W_{k-1-s}^{(k-1)}-\frac{1}{N}1_{N}1_{N}^{\top}\right)\otimes I_{d}\right)y^{(s)}\right\Vert^{2}
&\leq
3\eta^{2}\left(\sum_{s=0}^{k-1}\left\Vert\left(W_{k-1-s}^{(k-1)}-\frac{1}{N}1_{N}1_{N}^{\top}\right)\otimes I_{d}\right\Vert
\cdot\left\Vert y^{(s)}\right\Vert\right)^{2}
\\
&\leq
3\eta^{2}\left(\sum_{s=0}^{k-1}\left\Vert W_{k-1-s}^{(k-1)}-\frac{1}{N}1_{N}1_{N}^{\top}\right\Vert
\cdot\left\Vert y^{(s)}\right\Vert\right)^{2}
\\
&=
3\eta^{2}\left(\sum_{s=0}^{k-1}\bar{\gamma}_{k-1-s}^{(k-1)}
\cdot\left\Vert y^{(s)}\right\Vert\right)^{2}
\\
&=3\eta^{2}\left(\sum_{s=0}^{k-1}\bar{\gamma}_{k-1-s}^{(k-1)}\right)^{2}\left(\frac{\sum_{s=0}^{k-1}\bar{\gamma}_{k-1-s}^{(k-1)}
\cdot\left\Vert y^{(s)}\right\Vert}{\sum_{s=0}^{k-1}\bar{\gamma}_{k-1-s}^{(k-1)}}\right)^{2}
\\
&\leq
3\eta^{2}\left(\sum_{s=0}^{k-1}\bar{\gamma}_{k-1-s}^{(k-1)}\right)^{2}
\sum_{s=0}^{k-1}\frac{\bar{\gamma}_{k-1-s}^{(k-1)}}{\sum_{s=0}^{k-1}\bar{\gamma}_{k-1-s}^{(k-1)}}\left\Vert y^{(s)}\right\Vert^{2},
\end{align*}
where we used Jensen's inequality in the last step above.
Recall from Lemma~\ref{lem:0} that for every $k=0,1,2,\ldots$, 
$\mathbb{E}\left[\left\Vert y^{(k)}\right\Vert^{2}\right]
\leq D^{2}$,
where $D$ is defined in \eqref{eqn:D1:main}.
Therefore, we have
\begin{align*}
&3\eta^{2}\mathbb{E}\left[\left\Vert\sum_{s=0}^{k-1}\left(\left(W_{k-1-s}^{(k-1)}-\frac{1}{N}1_{N}1_{N}^{\top}\right)\otimes I_{d}\right)y^{(s)}\right\Vert^{2}\right]
\\
&\qquad\qquad \leq
3D^{2}\eta^{2}\left(\sum_{s=0}^{k-1}\bar{\gamma}_{k-1-s}^{(k-1)}\right)^{2}
\sum_{s=0}^{k-1}\frac{\bar{\gamma}_{k-1-s}^{(k-1)}}{\sum_{s=0}^{k-1}\bar{\gamma}_{k-1-s}^{(k-1)}}
\leq 
3D^{2}\eta^{2}\left(\sum_{s=0}^{k-1}\bar{\gamma}_{k-1-s}^{(k-1)}\right)^{2}.
\end{align*}
Similarly, we can show that
\begin{equation*}
3\left\Vert\left(\left(W_{k}^{(k-1)}-\frac{1}{N}1_{N}1_{N}^{\top}\right)\otimes I_{d}\right)x^{(0)}\right\Vert^{2}
\leq 3\left(\bar{\gamma}_{k}^{(k-1)}\right)^{2}\left\Vert x^{(0)}\right\Vert^{2}.
\end{equation*}
It follows from \eqref{eq:inter-estimate} that
\begin{align*}
\sum_{i=1}^{N}&\mathbb{E}\left\Vert x_{i}^{(k)}-\bar{x}^{(k)}\right\Vert^{2}
\\
&\leq
3\left(\bar{\gamma}_{k}^{(k-1)}\right)^{2}\mathbb{E}\left\Vert x^{(0)}\right\Vert^{2}
+3D^{2}\eta^{2}\left(\sum_{s=0}^{k-1}\bar{\gamma}_{k-1-s}^{(k-1)}\right)^{2}
\\
&\qquad\qquad\qquad
+6\eta\sum_{s=0}^{k-1}\mathbb{E}\left\Vert\left(\left(W_{k-1-s}^{(k-1)}-\frac{1}{N}1_{N}1_{N}^{\top}\right)\otimes I_{d}\right)w^{(s+1)}\right\Vert^{2}
\\
&\leq
3\left(\bar{\gamma}_{k}^{(k-1)}\right)^{2}\mathbb{E}\left\Vert x^{(0)}\right\Vert^{2}
+3D^{2}\eta^{2}\left(\sum_{s=0}^{k-1}\bar{\gamma}_{k-1-s}^{(k-1)}\right)^{2}
+6\eta\sum_{s=0}^{k-1}\left\Vert W_{k-1-s}^{(k-1)}-\frac{1}{N}1_{N}1_{N}^{\top}\right\Vert^{2}
\mathbb{E}\left\Vert w^{(s+1)}\right\Vert^{2}
\\
&\leq
3\left(\bar{\gamma}_{k}^{(k-1)}\right)^{2}\mathbb{E}\left\Vert x^{(0)}\right\Vert^{2}
+3D^{2}\eta^{2}\left(\sum_{s=0}^{k-1}\bar{\gamma}_{k-1-s}^{(k-1)}\right)^{2}
+6dN\eta\sum_{s=0}^{k-1}\left(\bar{\gamma}_{k-1-s}^{(k-1)}\right)^{2}.
\end{align*} 
The proof is complete.
\end{proof}

\subsection{Proof of Lemma~\ref{cor:1}}\label{app-sec-proof-lem-cor:1}

\begin{proof}

For $k=0$, the bound holds trivially. Assume $k\geq 1$. It follows from Lemma~\ref{lem:1} that for any $k\geq 1$, 
\begin{equation*}
\sum_{i=1}^{N}\mathbb{E}\left\Vert x_{i}^{(k)}-\bar{x}^{(k)}\right\Vert^{2}
\leq
3\left(\bar{\gamma}_{k}^{(k-1)}\right)^{2}\mathbb{E}\left\Vert x^{(0)}\right\Vert^{2}
+3D^{2}\eta^{2}\left(\sum_{s=0}^{k-1}\bar{\gamma}_{k-1-s}^{(k-1)}\right)^{2}
+6dN\eta\sum_{s=0}^{k-1}\left(\bar{\gamma}_{k-1-s}^{(k-1)}\right)^{2},
\end{equation*}
where $D$ is defined in \eqref{eqn:D1:main} and
\begin{equation}
\bar{\gamma}_{k-1-s}^{(k-1)}=\left\Vert W_{k-1-s}^{(k-1)}-\frac{1}{N}1_{N}1_{N}^{\top}\right\Vert.
\end{equation}
Under Assumption~\ref{assump:W}, 
there exists some positive integer $B$ such that
$\delta:=\sup_{k\geq B-1}\delta(k)<1$, where $\delta(k):=\sigma_{\max}\left\{W_{B}^{(k)}-\frac{1}{N}1_{N}1_{N}^{\top}\right\}$ for every $k=0,1,2,\ldots$.
For every $k$, $W^{(k)}$ is doubly stochastic. That means that for every $k$, we have $W^{(k)}J = J W^{(k)}$ with $J := \frac{1}{N}1_{N}1_{N}^{\top}$ where $J^2 = J$. Note that if $A$ and $M$ are $N\times N$ doubly stochastic matrices that are not necessarily symmetric, they satisfy $AJ = JA = J$ and $MJ=JM=J$ and we always have
$(A-J)(M-J)=AM-AJ-JM+J^2=AM-J$.
Furthermore, the product $AM$ is always double stochastic, even if it is not necessarily symmetric. Therefore, writing $j=mB+r$ with $m=\lfloor j/B\rfloor$ and $0\le r<B$, we obtain
\[
W^{(k-1)}_{j}-J
=\left(\prod_{\ell=0}^{m-1}\left(W^{(k-1-\ell B)}_{B}-J\right)\right)\,\left(W^{(k-1-mB)}_{r}-J\right),
\]
where matrices $W^{(k-1-mB)}_{r}$ and $W^{(k-1-\ell B)}_{B}$ are all double stochastic as products of double stochastic matrices. By part (iii) of Assumption~\ref{assump:W}, $\left\|W^{(k)}_{B}-J\right\|\le \delta$ for all $k\ge B-1$, and because $W^{(k)}$ is non-expansive
on the orthogonal complement of $1_N$ (i.e. the singular values of the symmetric matrix $W^{(k)}-J$ are at most 1),\footnote{Indeed, $W^{(k)}$ is a symmetric double stochastic matrix with all the eigenvalues lying in the interval $[-1,1]$ and admits $1$ as an eigenvalue with the eigenvector $1_N$. 
Since $W^{(k)}-J$ is symmetric, its singular values coincide with the absolute values of its eigenvalues, and hence $\|W^{(k)}-J\|_2 = \max_i |\lambda_i(W^{(k)}-J)| \le 1$.} we also have
$\left\|W^{(k)}_{r}-J\right\| \leq \left\|W^{(k)}-J\right\| \cdot \left\|W^{(k-1)}-J\right\| \cdots \left\|W^{(k-r+1)}-J\right\|
\le 1,$ for any $r\ge 0$ and $k$. 
Hence
\begin{equation}
\bar\gamma^{(k-1)}_{k-1-s}\ \le\ \delta^{\left\lfloor\frac{k-1-s}{B}\right\rfloor} \leq \delta^{\frac{k-1-s}{B}-1}, 
\qquad s=0,1,\ldots,k-1,
\end{equation}
and 
$\bar{\gamma}_{k}^{(k-1)}\leq 
\delta^{\frac{k}{B}-1}$. Therefore, we get
\begin{align*}
\sum_{i=1}^{N}\mathbb{E}\left\Vert x_{i}^{(k)}-\bar{x}^{(k)}\right\Vert^{2}
&\leq
3\left(\delta^{\frac{k}{B}-1}\right)^{2}\mathbb{E}\left\Vert x^{(0)}\right\Vert^{2}
+3D^{2}\eta^{2}\left(\sum_{s=0}^{k-1}\delta^{\frac{k-1-s}{B}-1}\right)^{2}
+6dN\eta\sum_{s=0}^{k-1}\left(\delta^{\frac{k-1-s}{B}-1}\right)^{2}
\\
&\leq
3\cdot \delta^{-2}\left(\delta^{\frac{2}{B}}\right)^{k}\mathbb{E}\left\Vert x^{(0)}\right\Vert^{2}
+\frac{3D^{2}\eta^{2}\delta^{-2}}{\left(1-\delta^{\frac{1}{B}}\right)^{2}}
+\frac{6dN\eta\cdot \delta^{-2}}{1-\delta^{\frac{2}{B}}}.
\end{align*}
The proof is complete.
\end{proof}


\subsection{Proof of Lemma~\ref{lem:2}}\label{app-sec-lem:2}

\begin{proof}
First, we can compute that
\begin{align*}
\mathbb{E}\left\Vert\mathcal{E}_{k}\right\Vert^{2}
=\mathbb{E}\left\Vert
\frac{1}{N}\sum_{i=1}^{N}\left(\nabla f_{i}\left(x_{i}^{(k)}\right)
-\nabla f_{i}\left(\bar{x}^{(k)}\right)\right)\right\Vert^{2}
\leq
\mathbb{E}\left\Vert
\frac{1}{N}\sum_{i=1}^{N}\left(\nabla f_{i}\left(x_{i}^{(k)}\right)
-\nabla f_{i}\left(\bar{x}^{(k)}\right)\right)\right\Vert^{2}.
\end{align*}
By Lemma~\ref{cor:1}, we can compute that
\begin{align*}
\mathbb{E}\left\Vert
\frac{1}{N}\sum_{i=1}^{N}\left(\nabla f_{i}\left(x_{i}^{(k)}\right)
-\nabla f_{i}\left(\bar{x}^{(k)}\right)\right)\right\Vert^{2}
&\leq
\frac{1}{N^{2}}\sum_{i=1}^{N}
N\mathbb{E}\left\Vert
\left(\nabla f_{i}\left(x_{i}^{(k)}\right)
-\nabla f_{i}\left(\bar{x}^{(k)}\right)\right)\right\Vert^{2}
\\
&\leq\frac{1}{N}L^{2}\sum_{i=1}^{N}
\mathbb{E}\left\Vert
x_{i}^{(k)}
-\bar{x}^{(k)}\right\Vert^{2}
\\
&\leq
\frac{3L^{2} \delta^{-2}}{N}\left(\delta^{\frac{2}{B}}\right)^{k}\mathbb{E}\left\Vert x^{(0)}\right\Vert^{2}
+\frac{3L^{2}D^{2}\eta^{2} \delta^{-2}}{N(1-\delta^{\frac{1}{B}})^{2}}
+\frac{6dL^{2}\eta\cdot \delta^{-2}}{1-\delta^{\frac{2}{B}}}.
\end{align*} 
The proof is complete.
\end{proof}


\section{Additional Technical Lemmas}\label{app:additional.lemmas}

\begin{lemma}[Lemma 3.4 in \cite{nedic2017}]\label{lem:iterate}
Under Assumption~\ref{assump:W}, for any $k=B-1,B,\ldots$ and any $Nd$-dimensional vector $b$,
we have 
$\left\Vert \mathcal{L}_{N}\mathcal{W}_{B}^{(k)}b\right\Vert\leq\delta(k)\left\Vert \mathcal{L}_{N}b\right\Vert$,
where $\delta(k)$ is defined in Assumption~\ref{assump:W}
and $\mathcal{L}_{N}:=I_{Nd}-\frac{1}{N}\left(\left(1_{N}1_{N}^{\top}\right)\otimes I_{d}\right)$.
\end{lemma}




Next, we consider
$\min_{x\in\mathbb{R}^{d}}g(x):=\frac{1}{N}\sum_{i=1}^{N}g_{i}(x)$,
where $g_{i}$ are $\mu$-strongly convex and $L$-smooth.
Consider the iterates:
\begin{equation}
p^{(k+1)}=p^{(k)}-\eta\frac{1}{N}\sum_{i=1}^{N}\nabla g_{i}\left(s_{i}^{(k)}\right),\label{eq-p-recursion}
\end{equation}
where the decision vector \(s_i^{(k)} \in \mathbb{R}^d\) of agent $i$ at time $k$ can be arbitrary. When iteration superscripts are omitted, \(p:=\{p^{(k)}\}_{k\geq 0}\) and \(s_i:=\{s_i^{(k)}\}_{k\geq 0}\) denote the corresponding sequences. Then, we have the following technical lemma. 
\begin{lemma}\label{lem:inexact}
Assume 
$\sqrt{1-\frac{\eta\mu\beta}{\beta+1}}\leq\lambda<1$
and $\eta\leq\frac{1}{(1+\alpha)L}$,
where $\alpha,\beta>0$.
Let the square-integrable random sequences $\{p^{(k)}\}_{k\geq0}$ and
$\{s_i^{(k)}\}_{k\geq0}$ satisfy \eqref{eq-p-recursion} almost surely where \(\nabla g_i\) is evaluated at \(s_i^{(k)}\) at iteration \(k\) for each \(i=1,\ldots,N\). \mgb{Consider the weighted $L_2$ norm
\[
\|p-p_*\|_{L_2}^{\lambda,K}:=
\max_{0\leq k\leq K}\lambda^{-k}
\bigl(\mathbb E\|p^{(k)}-p_*\|^2\bigr)^{1/2},
\]
where $\|p-s_i\|_{L_2}^{\lambda,K}$ is defined analogously. Then, for every
$K\geq0$,
\begin{equation}\label{eq:stochastic-inexact-L2}
\|p-p_*\|_{L_2}^{\lambda,K}
\leq
\bigl(\mathbb E\|p^{(0)}-p_*\|^2\bigr)^{1/2}
+\frac{1}{\lambda\sqrt N}
\sqrt{\frac{L(1+\alpha)}{\mu\alpha}+\beta}
\sum_{i=1}^N\|p-s_i\|_{L_2}^{\lambda,K},
\end{equation}}
where $p_{\ast}$ is the minimizer of  $g:=\frac{1}{N}\sum_{i=1}^{N}g_{i}(x)$.
\end{lemma}
\begin{proof} \mgb{The recursion \eqref{eq-p-recursion} was studied in Lemma 3.12 of \cite{nedic2017} in the deterministic DIGing setting, where a weighted-norm estimate for the corresponding inexact-gradient iteration was established. In the present work, however, the sequences $\{p^{(k)}\}_k$ and $\{s_i^{(k)}\}_k$ are stochastic, so the deterministic weighted-maximum argument used in \cite{nedic2017} does not directly provide the $L_2$ estimate required for our analysis. We therefore retain the underlying one-step strong-convexity and smoothness argument of Lemma 3.12 in \cite{nedic2017}, but reformulate the subsequent analysis at the level of expected squared errors. This yields a deterministic scalar recursion for the second moments, which can then be unrolled and converted into the desired weighted $L_2$ bound without interchanging expectation and a pathwise temporal maximum. The complete proof is deferred to the supplementary material as it parallels that of Lemma 3.12 in \cite{nedic2017}.}
\end{proof}

\begin{lemma}\label{lemma-surjective-gradients} 
For any function $f \in \mathcal{S}_{\mu,L}(\mathbb{R}^d)$, the gradient operator $\nabla f:\mathbb{R}^d\to \mathbb{R}^d$ is surjective, i.e. for every $v\in \mathbb{R}^d$, there exists some $x\in \mathbb{R}^d$ such that $\nabla f(x) = v$.
\end{lemma}

\begin{proof} 
This is a direct consequence of \cite[Theorem 1]{chiappinelli2020remarks}. Indeed, for $f \in \mathcal{S}_{\mu,L}(\mathbb{R}^d)$, the gradient operator $\nabla f$ is strongly coercive around $x_*$, i.e. it satisfies 
$\langle \nabla f(x)- \nabla f(x_*), x - x_* \rangle \geq \mu \|x-x_*\|^p$,
for $p=2$. As a consequence, it is coercive, i.e., 
$\frac{\langle \nabla f(x), x\rangle}{\|x\|} \to \infty$, as $\|x\|\to \infty$. 
It is also a bounded operator, i.e., $\|\nabla f(x) - \nabla f(x_*)\| \leq L \|x-x_*\|$. Finally, it is a proper operator, i.e., the preimage $\nabla f^{-1}(K)$ is a compact subset of $\mathbb{R}^d$ whenever $K \subset \mathbb{R}^d$ is compact. To see this, note that by strong convexity, $\|\nabla f(x)-\nabla f(x_*)\|\geq \mu \|x-x_*\|$; hence, the preimage $\nabla f^{-1}(K)$ of a compact set $K$ should be bounded. In addition, because $K$ is closed, such a preimage should also be closed by the continuity of $\nabla f$, which implies that $\nabla f^{-1}(K)$ is indeed compact. Therefore, \cite[Theorem 1]{chiappinelli2020remarks} is applicable and this completes the proof.
\end{proof}

\bibliographystyle{plain}
\bibliography{DIGing}

\newpage
\section*{Supplementary Materials}
In this supplementary material, we provide proofs of Lemma~\ref{lem-bound-D}, Corollary~\ref{coro-param}, Lemma~\ref{lem:3} and Lemma~\ref{lem:inexact}, which have been omitted in the main text for brevity.
{
\section{Proof of Lemma~\ref{lem-bound-D}}\label{app-sec-proof-lem:lem-bound-D}
We first provide a helper lemma that shows that the proposed parameters satisfy the conditions \eqref{assump:alpha:beta:main} and then prove a bound on \(D\sqrt\eta\) that is uniform in stepsize.

\begin{lemma}\label{lemma-lambda-lb-bis}
Let \(\alpha=1\), \(\beta=2L/\mu\), and let \(J_1,\eta_{\rm sg},\bar\eta\) and \(\lambda(\eta)\) be defined by
\eqref{def:J1-corrected-main}--\eqref{def-opt-lambda-eta}. Then, for every \(\eta\in(0,\bar\eta]\),
\[
\lambda(\eta)\in(\delta^{1/B},1),\qquad
\sqrt{1-\frac{\eta\mu\beta}{\beta+1}}\leq\lambda(\eta),\qquad
\eta\leq\frac1{2L},
\]
and
\[
0<\gamma_1\gamma_2\gamma_3\gamma_4\leq q_0<1,
\]
where \(q_0\) is defined in \eqref{def:lambda-rho-q-main}.
\end{lemma}

\begin{proof}
Since every \(\mu\)-strongly convex and \(L\)-smooth function satisfies \(L\geq\mu\), we have \(\kappa=L/\mu\geq1\). By definition,
\[
\bar\eta\leq\frac1{2L}
\quad\text{and}\quad
\bar\eta\leq\frac{3}{4\mu}(1-\delta^{2/B}).
\]
Thus, for \(0<\eta\leq\bar\eta\),
\begin{equation}\label{eq:network-gap-bis}
1-\frac{\eta\mu}{1.5}
\geq1-\frac{\bar\eta\mu}{1.5}
\geq\frac{1+\delta^{2/B}}{2}>\delta^{2/B}.
\end{equation}
In particular, \(0<\lambda(\eta)<1\), and
\[
\lambda(\eta)^B=\sqrt{1-\frac{\eta\mu}{1.5}}
>\delta^{1/B}\geq\delta,
\]
which implies \(\lambda(\eta)>\delta^{1/B}\). Moreover,
\[
\frac{\beta}{\beta+1}=\frac{2\kappa}{2\kappa+1}\geq\frac23,
\]
so
\[
\sqrt{1-\frac{\eta\mu\beta}{\beta+1}}
\leq\sqrt{1-\frac{\eta\mu}{1.5}}
=\lambda(\eta)^B\leq\lambda(\eta).
\]
This proves the first two conditions in \eqref{assump:alpha:beta:main}; the stepsize condition follows from \(\bar\eta\leq1/(2L)\).

It remains to verify the small-gain condition $\gamma_1 \gamma_2 \gamma_3 \gamma_4 \leq q_0<1$. With \(\alpha=1\), \(\beta=2L/\mu\), the definition \eqref{defn:gamma:3:4:main} gives
\[
\gamma_3=1+\frac{2N\sqrt\kappa}{\lambda}.
\]
Writing \(r:=\lambda^B=\sqrt{1-\eta\mu/1.5}\) and using
\[
\frac{1-\lambda^B}{1-\lambda}=\sum_{j=0}^{B-1}\lambda^j\leq B,
\]
we obtain
\begin{align}
\gamma_1\gamma_2\gamma_3\gamma_4
&=\eta L(1+\lambda)\left(1+\frac{2N\sqrt\kappa}{\lambda}\right)
\frac{\left(\frac{1-\lambda^B}{1-\lambda}\right)^2}{(\lambda^B-\delta)^2}\nonumber\\
&\leq 2\eta L B^2\left(1+\frac{2N\sqrt\kappa}{\underline\lambda}\right)
\frac1{(r-\delta)^2}.
\label{eq:small-gain-step-bis}
\end{align}
We next consider the equation
\begin{equation} r-\delta = \sqrt{1-\eta \mu/1.5}-\delta =\sqrt{\frac{\eta\mu J_1}{1.5}}, 
\label{eq-stepsize}
\end{equation}
we will argue that $\eta = \eta_{\rm sg}$ satisfies this equation. Indeed, if we define
\begin{equation}
y:=\sqrt{\frac{\eta\mu J_1}{1.5}}.
\end{equation}
Then 
the equation \eqref{eq-stepsize} is equivalent to
\begin{equation}
\sqrt{1-\frac{y^2}{J_1}}
=
\delta+y.
\end{equation}
Squaring gives
\begin{equation}
(J_1+1)y^2
+2\delta J_1y
-J_1(1-\delta^2)
=0.
\end{equation}
The positive root is
\begin{equation}
\label{eq:y-root}
y
=
\frac{
\sqrt{J_1^2+(1-\delta^2)J_1}
-\delta J_1
}{J_1+1}.
\end{equation}
Finally, since
\begin{equation}
y^2
=
\frac{\eta\mu J_1}{1.5},
\end{equation}
substituting \eqref{eq:y-root} and solving for $\eta$ gives exactly the expression \eqref{bar:eta:defn}. We conclude that $\eta_{\rm sg}$ satisfies  \eqref{eq-stepsize}, in other words, 
\begin{equation}\label{eq:etasg-identity-bis}
r_{\rm sg}-\delta=\sqrt{1-\eta_{\rm sg}\mu/1.5}-\delta
=\sqrt{\frac{\eta_{\rm sg}\mu J_1}{1.5}},
\end{equation}
where \(r_{\rm sg}:=\sqrt{1-\eta_{\rm sg}\mu/1.5}\). Because \(\eta\leq\bar\eta\leq\eta_{\rm sg}\), we have \(r\geq r_{\rm sg}\), and hence
\[
\frac{\eta}{(r-\delta)^2}
\leq\frac{\eta_{\rm sg}}{(r_{\rm sg}-\delta)^2}
=\frac{1.5}{\mu J_1}.
\]
Substituting this estimate into \eqref{eq:small-gain-step-bis} yields
\[
\gamma_1\gamma_2\gamma_3\gamma_4
\leq\frac{3\kappa B^2}{J_1}
\left(1+\frac{2N\sqrt\kappa}{\underline\lambda}\right)
=q_0.
\]
Finally, \(\bar\eta\leq1/(2L)\) implies
\[
\underline\lambda
=\left(1-\frac{\bar\eta\mu}{1.5}\right)^{1/(2B)}
\geq\left(1-\frac{1}{3\kappa}\right)^{1/(2B)}
\geq\sqrt{\frac23}>\frac12.
\]
Therefore
\[
1+\frac{2N\sqrt\kappa}{\underline\lambda}
<1+4N\sqrt\kappa,
\]
and thus \(q_0<1\). This proves all the conditions in \eqref{assump:alpha:beta:main}.
\end{proof}

\begin{lemma}[Restatement of Lemma~\ref{lem-bound-D}]\label{lem-bound-D-bis}
Under the setting and parameter choices of Lemma~\ref{lem-bound-D}, the estimate
\[
D\sqrt\eta\leq\overline D
\]
holds for every \(\eta\in(0,\bar\eta]\), where \(\overline D\) is given by \eqref{eqn:D1:main:upperbound}.
\end{lemma}

\begin{proof}
By Lemma~\ref{lemma-lambda-lb-bis}, the hypotheses of Theorems~\ref{thm:average} and \ref{thm:undirected} hold and
\begin{equation}\label{eq:gain-margin-bis}
1-\gamma_1\gamma_2\gamma_3\gamma_4\geq1-q_0>0.
\end{equation}
We first establish bounds on the finite collection of initial transient terms appearing in \(\tilde\omega_1\) and \(\tilde\omega_4\). Since every \(W^{(k)}\) is symmetric and doubly stochastic with nonnegative entries, \(\|\mathcal W^{(k)}\|_{\rm op}\leq1\). Assumption~\ref{assumption:noise}, the Lipschitz continuity of the gradients, and \(\|x^{(0)}\|_{L_2}<\infty\) imply \(Y_0=\|y^{(0)}\|_{L_2}<\infty\). Also,
\[
\|w^{(k)}\|_{L_2}=\sqrt{Nd},\qquad
\|\xi^{(k)}\|_{L_2}\leq\sqrt N\,\sigma.
\]
Using \eqref{x:vector:form}--\eqref{y:vector:form}, the \(L\)-Lipschitz continuity of \(\nabla F\), and the triangle inequality in \(L_2\), we obtain, for every \(0<\eta\leq\bar\eta\),
\begin{align*}
\|x^{(t+1)}\|_{L_2}
&\leq\|x^{(t)}\|_{L_2}+\bar\eta\|y^{(t)}\|_{L_2}+\sqrt{2\bar\eta Nd},\\
\|y^{(t+1)}\|_{L_2}  
&\leq \|y^{(t)}\|_{L_2} +  L \|x^{(t+1)}-x^{(t)}\|_{L_2} + 2\sqrt N\,\sigma \\
&\leq(1+L\bar\eta)\|y^{(t)}\|_{L_2}
+2L\|x^{(t)}\|_{L_2}
+L\sqrt{2\bar\eta Nd}+2\sqrt N\,\sigma,
\end{align*}
where we used the bound
\begin{eqnarray*}\|x^{(t+1)}-x^{(t)}\|_{L_2} &\leq& \|(W^{t} - I) x^{(t)}\|_{L_2} + \eta \| y^{(t)}\|_{L_2} + \sqrt{2\eta}\| w^{t+1}\|_{L_2} \\ 
&\leq& 2\|x^{(t)}\|_{L_2} + \eta \| y^{(t)}\|_{L_2} +  \sqrt{2\eta Nd}~.
\end{eqnarray*}
Consequently, the recursively defined quantities \(X_t,Y_t\) in \eqref{def:finite-horizon-bounds-main} satisfy
\begin{equation}\label{eq:finite-horizon-domination-bis}
\|\tilde x^{(t)}\|_{L_2}\leq\|x^{(t)}\|_{L_2}\leq X_t,
\qquad
\|\tilde y^{(t)}\|_{L_2}\leq\|y^{(t)}\|_{L_2}\leq Y_t,
\quad 0\leq t\leq B-1.
\end{equation}
In particular, these bounds are finite and independent of \(\eta\) and \(\epsilon\).

Next, define the scaled quantities
\[
\tilde\omega_i^{(\eta)}:=\sqrt\eta\,\tilde\omega_i,
\qquad
\hat\omega_i^{(\eta)}:=\sqrt\eta\,\hat\omega_i,
\qquad i\in\{1,3,4\}.
\]
Because \(\lambda\geq\underline\lambda\), the function \(s\mapsto s/(s-\delta)\) is decreasing for \(s>\delta\), and \(\lambda^{1-t}\leq\underline\lambda^{1-t}\) for \(t\geq1\). Combining these observations with \eqref{eq:finite-horizon-domination-bis} gives
\begin{align}
\tilde\omega_1^{(\eta)}+\hat\omega_1^{(\eta)}&\leq\Omega_1,\nonumber\\
\tilde\omega_3^{(\eta)}+\hat\omega_3^{(\eta)}&\leq\Omega_3,\nonumber\\
\tilde\omega_4^{(\eta)}+\hat\omega_4^{(\eta)}&\leq\Omega_4,
\label{eq:scaled-omega-corrected-bis}
\end{align}
where \(\Omega_1,\Omega_3,\Omega_4\) are defined in \eqref{def:Omega-main}. For example,
\[
\sqrt\eta\,\hat\omega_3
=\frac{2\sqrt{N\kappa}}{\mu\lambda}
\left(\sigma\sqrt\eta+\sqrt{2d}\right)
\leq\frac{2\sqrt{N\kappa}}{\mu\underline\lambda}
\left(\sigma\sqrt{\bar\eta}+\sqrt{2d}\right),
\]
which is the second term in \(\Omega_3\).

The same elementary bounds give
\begin{equation}\label{eq:gamma-corrected-bis}
\gamma_1\leq\Gamma_1,
\qquad \gamma_2\leq\Gamma_2,
\qquad \gamma_3\leq\Gamma_3,
\qquad \gamma_4\leq\Gamma_4,
\end{equation}
with \(\Gamma_i\) defined in \eqref{def:Gamma-main}. Indeed,
\[
\gamma_1
=\frac{\lambda}{\lambda^B-\delta}\sum_{j=0}^{B-1}\lambda^j
\leq\frac{B}{\underline\lambda^B-\delta},
\qquad
\gamma_4
=\frac{\eta}{\lambda^B-\delta}\sum_{j=0}^{B-1}\lambda^j
\leq\frac{B\bar\eta}{\underline\lambda^B-\delta}.
\]

Finally, multiply the definition \eqref{eqn:D1:main} by \(\eta\), use \eqref{eq:gain-margin-bis}, \eqref{eq:scaled-omega-corrected-bis}, \eqref{eq:gamma-corrected-bis}, and \(\eta\leq\bar\eta\). This gives exactly
\[
D\sqrt\eta\leq\overline D
\]
with \(\overline D\) defined in \eqref{eqn:D1:main:upperbound}. Since all quantities entering \(\overline D\) are finite and independent of \(\epsilon\), the same is true of \(\overline D\). This completes the proof.
\end{proof}
}
{
\section{Proof of Corollary~\ref{coro-param}}\label{coro-param-proof}
We restate the corollary and then verify the error bound.

\begin{corollary}[Restatement of Corollary~\ref{coro-param}]\label{coro-param-bis}
Use the parameter choices and constants in Lemma~\ref{lem-bound-D} and Corollary~\ref{coro-param}. Then, with \(\eta_*\) defined in \eqref{def:eta-star-corrected-main},
\[
\frac1N\sum_{i=1}^N
\mathcal W_2\!\left(\operatorname{Law}(x_i^{(k)}),\pi\right)\leq\epsilon
\]
for every \(k\) satisfying \eqref{def:k-star-corrected-main}.
\end{corollary}

\begin{proof}
Fix \(\eta=\eta_*\), and let
\[
a:=\delta^{2/B},\qquad
b_1:=1-\eta\mu\left(1-\frac{\eta L}{2}\right),\qquad
b_2:=1-\frac{\eta\mu}{1.5}.
\]
By Lemma~\ref{lem-bound-D}, \(\eta\leq\bar\eta\leq1/(2L)\), \(D\sqrt\eta\leq\overline D\), and
\begin{equation}\label{eq:rho-gap-cor-bis}
b_2-a\geq\rho_0>0.
\end{equation}
Moreover,
\[
b_2-b_1
=\eta\mu\left(\frac13-\frac{\eta L}{2}\right)
\geq\frac{\eta\mu}{12}>0.
\]
Recall from \eqref{eq:def:Hk} that
\[
\mathsf H_k(r,q)=\sum_{j=0}^{k-1}r^{k-1-j}q^j.
\]
Since this is a sum of nonnegative monomials, $\mathsf H_k$ is nondecreasing in each argument on \([0,1)\). Hence
\[
\mathsf H_k(b_1,a)\leq \mathsf H_k(b_2,a)
=\frac{b_2^k-a^k}{b_2-a}
\leq\frac{b_2^k}{\rho_0}.
\]
This also covers the case \(k=0\), as $H_0(b_1,a) = 0$. Moreover, because \(\eta\leq\bar\eta\),
\begin{align*}
P(\eta)&:=\left[\eta\left(\eta+
\frac{(1+\eta L)^2}{\mu(1-\eta L/2)}\right)\right]^{1/2}\\
&\leq\left[\bar\eta\left(\bar\eta+
\frac{(1+\bar\eta L)^2}{\mu(1-\bar\eta L/2)}\right)\right]^{1/2}
=P_0.
\end{align*}
Thus the last term in \(E_2\) is bounded by
\[
\frac{P_0}{\sqrt{\rho_0}}
\frac{\sqrt3L\delta^{-1}}{\sqrt N}\|x^{(0)}\|_{L_2}
\,b_2^{k/2}.
\]

Starting from \(E_1,E_2,E_3\) in Theorem~\ref{thm:undirected}, use \(D\sqrt\eta\leq\overline D\),
\((u+v)^{1/2}\leq u^{1/2}+v^{1/2}\), \(1-\eta L/2\geq3/4\), and the preceding estimates. A term-by-term calculation yields
\begin{equation}\label{ineq-E2-bound-final}
E_2\leq 
\frac{P_0}{\sqrt{\rho_0}}
\frac{\sqrt3L\delta^{-1}}{\sqrt N}\|x^{(0)}\|_{L_2}\,b_2^{k/2}
+C_{3,E_2}\sqrt\eta+C_{4,E_2}\eta,
\end{equation}
and
\begin{equation}\label{ineq-E3-bound-final}
E_3\leq
\frac{\sqrt3\delta^{-1}}{\sqrt N}\|x^{(0)}\|_{L_2}\,\delta^{k/B}
+\frac{\sqrt3 \cdot\overline D\delta^{-1}}{\sqrt N(1-\delta^{1/B})}\sqrt\eta
+\frac{\sqrt{6d}\,\delta^{-1}}{\sqrt{1-\delta^{2/B}}}\sqrt\eta,
\end{equation}
where, with
\[
A_D:=\left(\frac{3L^2\overline D^2\delta^{-2}}{N(1-\delta^{1/B})^2}\right)^{1/2},
\qquad
A_d:=\left(\frac{6dL^2\delta^{-2}}{1-\delta^{2/B}}\right)^{1/2},
\]
we may take
\[
C_{3,E_2}:=\frac{2}{\mu}(A_D+A_d)+\frac{2\sigma}{\sqrt{3\mu N}},
\qquad
C_{4,E_2}:=\frac{2}{\sqrt{3\mu}}(A_D+A_d).
\]
Indeed,
\[
\frac{1+\eta L}{\mu(1-\eta L/2)}\leq\frac2\mu,
\qquad
\frac1{\sqrt{\mu(1-\eta L/2)}}\leq\frac2{\sqrt{3\mu}},
\]
so these bounds follow directly by splitting each square root in the first term of \(E_2\). 

Next, \eqref{eq:network-gap-bis} implies
\[
\delta^{1/B}<\sqrt{b_2}.
\]
Also, since \(\eta\mu\leq\mu/(2L)\leq1/2\),
\[
0<1-\eta\mu\leq\sqrt{1-\frac{\eta\mu}{1.5}}=\sqrt{b_2}.
\]
Consequently, combining \(E_1\), \eqref{ineq-E2-bound-final}, and \eqref{ineq-E3-bound-final} gives
\begin{equation}\label{ineq-final-error-bound}
E_1+E_2+E_3
\leq(\overline C_1+\overline C_2)(\sqrt{b_2})^k
+\overline C_3\sqrt\eta+\overline C_4\eta,
\end{equation}
where $\overline C_3$ is obtained by adding the \(\sqrt\eta\)-bias in \(E_1\) and the two \(\sqrt\eta\)-terms in \eqref{ineq-E3-bound-final} to $C_{3,E_2}$ and $\overline C_4 = C_{4,E_2}$. By the definition of \(\eta_{\rm noise}(\epsilon)\),
\[
\overline C_3\sqrt\eta\leq\frac\epsilon3,
\qquad
\overline C_4\eta\leq\frac\epsilon3.
\]
Furthermore,
\[
(\sqrt{b_2})^k
=\left(\sqrt{1-\frac{\eta\mu}{1.5}}\right)^k
\leq\exp\left(-\frac{\eta\mu k}{3}\right).
\]
Therefore, if
\[
k\geq\frac{3}{\mu\eta}
\log\left(\frac{3(\overline C_1+\overline C_2)}{\epsilon}\right),
\]
the first term in \eqref{ineq-final-error-bound} is at most \(\epsilon/3\). With \(\eta=\eta_*\), this is precisely \eqref{def:k-star-corrected-main}, and hence \(E_1+E_2+E_3\leq\epsilon\). Theorem~\ref{thm:undirected} now gives the claimed Wasserstein bound.
\end{proof}
}
\section{Proof of Lemma~\ref{lem:3}}\label{app-sec-lem:3-bis} We next provide a restatement of Lemma~\ref{lem:3} followed by its proof. 
\begin{lemma}[Restatement of Lemma~\ref{lem:3}]\label{lem:3-bis}
Assume that $\delta<\lambda^{B}<1$ with $\delta(k)$ defined in Assumption~\ref{assump:W}
and \eqref{assump:alpha:beta} holds.
We also assume $\mathbb{E}\left\Vert x^{(0)}\right\Vert^{2}<\infty$, where the expectation is taken with respect to the distribution of the initial iterate $x^{(0)}$ at time \(k=0\).
For any stepsize $\eta \in (0,2/L)$, we have for every $k$,
{\small
\begin{align*}
\mathbb{E}\left\Vert\bar{x}^{(k)}-x_{k}\right\Vert^{2}
&\leq
\eta\left(\frac{\eta}{\mu(1-\frac{\eta L}{2})}+\frac{(1+\eta L)^{2}}{\mu^{2}(1-\frac{\eta L}{2})^{2}}\right)
\left(\frac{3L^{2}D^{2}\eta  \delta^{-2}}{N(1-\delta^{\frac{1}{B}})^{2}}
+\frac{6dL^{2}\cdot  \delta^{-2}}{1-\delta^{\frac{2}{B}}}\right)
\\
&\quad+\frac{\eta\sigma^{2}}{\mu(1-\frac{\eta L}{2})N}
\nonumber
\\
&\quad
+\mgb{\eta\left(\eta+\frac{(1+\eta L)^{2}}{\mu(1-\frac{\eta L}{2})}\right)
\mathsf H_k\left(1-\eta\mu\left(1-\frac{\eta L}{2}\right),\delta^{\frac{2}{B}}\right)}
\frac{3L^{2} \delta^{-2}}{N}\mathbb{E}\left\Vert x^{(0)}\right\Vert^{2}.
\end{align*}}
\end{lemma}
\begin{proof}
The proof is similar to the proof of Lemma~7 in \cite{gurbuzbalaban2021decentralized} 
and for the sake of completeness we include all the details here.
From \eqref{bar:x:iterates} and \eqref{overdamped:iterates}, 
we can compute that
\begin{equation*}
\bar{x}^{(k+1)}-x_{k+1}=\bar{x}^{(k)}-x_{k}
-\frac{\eta}{N}\left[\nabla f\left(\bar{x}^{(k)}\right)-\nabla f(x_{k})\right]
+\eta\mathcal{E}_{k}
-\eta\bar{\xi}^{(k+1)},
\end{equation*}
where we recall from \eqref{E:k:defn} that
$\mathcal{E}_{k}
=\frac{1}{N} \nabla f\left(\bar{x}^{(k)}\right)-\frac{1}{N}\sum_{i=1}^{N}\nabla f_{i}\left(x_{i}^{(k)}\right)$,
and this implies that
{\small
\begin{align}
\left\Vert\bar{x}^{(k+1)}-x_{k+1}\right\Vert^{2}
&=\left\Vert\bar{x}^{(k)}-x_{k}
-\frac{\eta}{N}\left[\nabla f\left(\bar{x}^{(k)}\right)-\nabla f(x_{k})\right]\right\Vert^{2}
+\eta^{2}\left\Vert\mathcal{E}_{k}-\bar{\xi}^{(k+1)}\right\Vert^{2}
\nonumber
\\
&\qquad\qquad
+2\left\langle\bar{x}^{(k)}-x_{k}
-\frac{\eta}{N}\left[\nabla f\left(\bar{x}^{(k)}\right)-\nabla f(x_{k})\right]
,\eta\mathcal{E}_{k}-\eta\bar{\xi}^{(k+1)}\right\rangle
\nonumber
\\
&=\left\Vert\bar{x}^{(k)}-x_{k}\right\Vert^{2}
+\eta^{2}\left\Vert\frac{1}{N}\left[\nabla f\left(\bar{x}^{(k)}\right)-\nabla f(x_{k})\right]\right\Vert^{2}
\nonumber
\\
&\qquad
-2\left\langle\bar{x}^{(k)}-x_{k},\eta\frac{1}{N}\left[\nabla f\left(\bar{x}^{(k)}\right)-\nabla f(x_{k})\right]\right\rangle
+\eta^{2}\left\Vert\mathcal{E}_{k}-\bar{\xi}^{(k+1)}\right\Vert^{2}
\nonumber
\\
&\qquad\qquad
+2\left\langle\bar{x}^{(k)}-x_{k}
-\eta\frac{1}{N}\left[\nabla f\left(\bar{x}^{(k)}\right)-\nabla f(x_{k})\right]
,\eta\mathcal{E}_{k}-\eta\bar{\xi}^{(k+1)}\right\rangle
\nonumber
\\
&\leq
\left\Vert\bar{x}^{(k)}-x_{k}\right\Vert^{2}
+\eta^{2}L\left\langle\bar{x}^{(k)}-x_{k},\frac{1}{N}\left[\nabla f\left(\bar{x}^{(k)}\right)-\nabla f(x_{k})\right]\right\rangle
\nonumber
\\
&\qquad
-2\left\langle\bar{x}^{(k)}-x_{k},\eta\frac{1}{N}\left[\nabla f\left(\bar{x}^{(k)}\right)-\nabla f(x_{k})\right]\right\rangle
+\eta^{2}\left\Vert\mathcal{E}_{k}-\bar{\xi}^{(k+1)}\right\Vert^{2}
\nonumber
\\
&\qquad\qquad
+2\left\langle\bar{x}^{(k)}-x_{k}
-\eta\frac{1}{N}\left[\nabla f\left(\bar{x}^{(k)}\right)-\nabla f(x_{k})\right]
,\eta\mathcal{E}_{k}-\eta\bar{\xi}^{(k+1)}\right\rangle
\nonumber
\\
&\leq
\left(1-2\eta\mu\left(1-\frac{\eta L}{2}\right)\right)\left\Vert\bar{x}^{(k)}-x_{k}\right\Vert^{2}
+\eta^{2}\left\Vert\mathcal{E}_{k}-\bar{\xi}^{(k+1)}\right\Vert^{2}
\nonumber
\\
&\qquad\qquad
+2\left\langle\bar{x}^{(k)}-x_{k}
-\eta\frac{1}{N}\left[\nabla f\left(\bar{x}^{(k)}\right)-\nabla f(x_{k})\right]
,\eta\mathcal{E}_{k}-\eta\bar{\xi}^{(k+1)}\right\rangle,\label{take:expect}
\end{align}} 
where we used $L$-smoothness of $\frac{1}{N}f$ to obtain the second term after
the first inequality above 
and $\mu$-strongly convexity of $\frac{1}{N}f$ and the assumption that $\eta< 2/L$ to obtain the first term after the second inequality above.
\mgb{Set
\[
A_k:=\bar x^{(k)}-x_k-\frac{\eta}{N}\bigl(\nabla f(\bar x^{(k)})-\nabla f(x_k)\bigr).
\]
Under Assumption~\ref{assumption:noise}, the quantities $x^{(k)}$, $\bar x^{(k)}$,
$\mathcal E_k$, and the comparison iterate $x_k$ are $\mathcal H_k$-measurable; the
last assertion holds because $x_k$ uses only $\bar w^{(1)},\ldots,\bar w^{(k)}$.
Moreover, $\mathbb E[\bar\xi^{(k+1)}\mid\mathcal H_k]=0$. Hence
\begin{align*}
\mathbb E\langle A_k,\bar\xi^{(k+1)}\rangle
&=\mathbb E\left\langle A_k,
\mathbb E[\bar\xi^{(k+1)}\mid\mathcal H_k]\right\rangle=0,\\
\mathbb E\langle\mathcal E_k,\bar\xi^{(k+1)}\rangle&=0.
\end{align*}
Conditional independence across agents also gives $
\mathbb E\|\bar\xi^{(k+1)}\|^2\leq\frac{\sigma^2}{N}.
$}
\mgb{Note that $\bar{\xi}^{(k+1)}$ has mean zero and by the conditional-orthogonality calculation above, every cross term involving $\bar\xi^{(k+1)}$ and an $\mathcal H_k$-measurable multiplier has zero expectation.}
In addition, by Lemma~\ref{lem:2},
\begin{equation}\label{apply:tilde}
\mathbb{E}\left\Vert\mathcal{E}_{k}\right\Vert^{2}
\leq
\frac{3L^{2} \delta^{-2}}{N}\left(\delta^{\frac{2}{B}}\right)^{k}\mathbb{E}\left\Vert x^{(0)}\right\Vert^{2}
+\frac{3L^{2}D^{2}\eta^{2} \delta^{-2}}{N(1-\delta^{\frac{1}{B}})^{2}}
+\frac{6dL^{2}\eta\cdot  \delta^{-2}}{1-\delta^{\frac{2}{B}}}.
\end{equation}
Taking expectations in \eqref{take:expect}, therefore gives
{\small
\begin{align*}
\mathbb{E}\left\Vert\bar{x}^{(k+1)}-x_{k+1}\right\Vert^{2}
&\leq
\left(1-2\eta\mu\left(1-\frac{\eta L}{2}\right)\right)\mathbb{E}\left\Vert\bar{x}^{(k)}-x_{k}\right\Vert^{2}
+\eta^{2}\mathbb{E}\left\Vert\mathcal{E}_{k}-\bar{\xi}^{(k+1)}\right\Vert^{2}
\nonumber
\\
&\qquad
+\mathbb{E}\left[2\left\langle\bar{x}^{(k)}-x_{k}
-\eta\frac{1}{N}\left[\nabla f\left(\bar{x}^{(k)}\right)-\nabla f(x_{k})\right]
,\eta\mathcal{E}_{k}-\eta\bar{\xi}^{(k+1)}\right\rangle\right]
\\
&=\left(1-2\eta\mu\left(1-\frac{\eta L}{2}\right)\right)\mathbb{E}\left\Vert\bar{x}^{(k)}-x_{k}\right\Vert^{2}
+\eta^{2}\mathbb{E}\left\Vert\mathcal{E}_{k}\right\Vert^{2}
+\eta^{2}\mathbb{E}\left\Vert\bar{\xi}^{(k+1)}\right\Vert^{2}
\nonumber
\\
&\qquad
+\mathbb{E}\left[2\left\langle\bar{x}^{(k)}-x_{k}
-\eta\frac{1}{N}\left[\nabla f\left(\bar{x}^{(k)}\right)-\nabla f(x_{k})\right]
,\eta\mathcal{E}_{k}\right\rangle\right]
\\
&\leq
\left(1-2\eta\mu\left(1-\frac{\eta L}{2}\right)\right)\mathbb{E}\left\Vert\bar{x}^{(k)}-x_{k}\right\Vert^{2}
+\eta^{2}\mathbb{E}\left\Vert\mathcal{E}_{k}\right\Vert^{2}
+\eta^{2}\frac{\sigma^{2}}{N}\\
&\qquad+2(1+\eta L)
\eta\mathbb{E}\left[\left\Vert\bar{x}^{(k)}-x_{k}\right\Vert\cdot\left\Vert\mathcal{E}_{k}\right\Vert\right],
\end{align*}}
where we used $L$-smoothness of $\frac{1}{N}f$.

For any $x,y\geq 0$ and $c>0$, we have the inequality $2xy\leq cx^{2}+\frac{y^{2}}{c}$,
which implies that
{\small
\begin{align*}
\mathbb{E}\left\Vert\bar{x}^{(k+1)}-x_{k+1}\right\Vert^{2}
&\leq
\left(1-2\eta\mu\left(1-\frac{\eta L}{2}\right)\right)\mathbb{E}\left\Vert\bar{x}^{(k)}-x_{k}\right\Vert^{2}
+\eta^{2}\mathbb{E}\left\Vert\mathcal{E}_{k}\right\Vert^{2}
+\eta^{2}\frac{\sigma^{2}}{N}
\\
&\qquad
+(1+\eta L)\eta
\left(\frac{\mu(1-\frac{\eta L}{2})}{1+\eta L}\mathbb{E}\left\Vert\bar{x}^{(k)}-x_{k}\right\Vert^{2}
+\frac{1+\eta L}{\mu(1-\frac{\eta L}{2})}\mathbb{E}\left\Vert\mathcal{E}_{k}\right\Vert^{2}\right)
\\
&=\left(1-\eta\mu\left(1-\frac{\eta L}{2}\right)\right)\mathbb{E}\left\Vert\bar{x}^{(k)}-x_{k}\right\Vert^{2}
+\eta\left(\eta+\frac{(1+\eta L)^{2}}{\mu(1-\frac{\eta L}{2})}\right)\mathbb{E}\left\Vert\mathcal{E}_{k}\right\Vert^{2}\\
&\qquad
+\eta^{2}\frac{\sigma^{2}}{N}.\end{align*}}By applying \eqref{apply:tilde}, we get
{\small
\begin{align*}
&\mathbb{E}\left\Vert\bar{x}^{(k+1)}-x_{k+1}\right\Vert^{2}
\leq\left(1-\eta\mu\left(1-\frac{\eta L}{2}\right)\right)
\mathbb{E}\left\Vert\bar{x}^{(k)}-x_{k}\right\Vert^{2} +\eta^{2}\frac{\sigma^{2}}{N}
\\
&\qquad\qquad\qquad\qquad
+\eta\left(\eta+\frac{(1+\eta L)^{2}}{\mu(1-\frac{\eta L}{2})}\right)
\Bigg(\frac{3L^{2} \delta^{-2}}{N}\left(\delta^{\frac{2}{B}}\right)^{k}\mathbb{E}\left\Vert x^{(0)}\right\Vert^{2} \\
&\qquad\qquad\qquad\qquad\qquad\qquad
+\frac{3L^{2}D^{2}\eta^{2} \delta^{-2}}{N(1-\delta^{\frac{1}{B}})^{2}}
+\frac{6dL^{2}\eta\cdot  \delta^{-2}}{1-\delta^{\frac{2}{B}}}\Bigg),
\end{align*}}for every $k$. Note that $\mathbb{E}\left\Vert\bar{x}^{(0)}-x_{0}\right\Vert^{2}=0$.
By iterating the above equation, we get for $k\geq 1$,
{\fontsize{8pt}{10pt}\selectfont
\begin{align*}
\mathbb{E}\left\Vert\bar{x}^{(k)}-x_{k}\right\Vert^{2}
&\leq
\sum_{i=0}^{k-1}
\left(1-\eta\mu\left(1-\frac{\eta L}{2}\right)\right)^{i}
\\
&\qquad\qquad
\cdot\left(\eta\left(\eta+\frac{(1+\eta L)^{2}}{\mu(1-\frac{\eta L}{2})}\right)
\left(\frac{3L^{2}D^{2}\eta^{2} \delta^{-2}}{N(1-\delta^{\frac{1}{B}})^{2}}
+\frac{6dL^{2}\eta\cdot  \delta^{-2}}{1-\delta^{\frac{2}{B}}}
\right)+\eta^{2}\frac{\sigma^{2}}{N}\right)
\\
&\qquad
+\sum_{i=0}^{k-1}
\left(1-\eta\mu\left(1-\frac{\eta L}{2}\right)\right)^{i}
\eta\left(\eta+\frac{(1+\eta L)^{2}}{\mu(1-\frac{\eta L}{2})}\right)
\frac{3L^{2} \delta^{-2}}{N}\left(\delta^{\frac{2}{B}}\right)^{\mgb{k-1-i}}\mathbb{E}\left\Vert x^{(0)}\right\Vert^{2}
\\
&=\frac{1-
\left(1-\eta\mu\left(1-\frac{\eta L}{2}\right)\right)^{k}}
{1-\left(1-\eta\mu\left(1-\frac{\eta L}{2}\right)\right)}
\cdot\Bigg(\eta\left(\eta+\frac{(1+\eta L)^{2}}{\mu(1-\frac{\eta L}{2})}\right)
\left(\frac{3L^{2}D^{2}\eta^{2} \delta^{-2}}{N(1-\delta^{\frac{1}{B}})^{2}}
+\frac{6dL^{2}\eta\cdot  \delta^{-2}}{1-\delta^{\frac{2}{B}}}\right)\\
&\qquad\qquad\qquad\qquad\qquad\qquad\qquad\qquad+\eta^{2}\frac{\sigma^{2}}{N}\Bigg)
\\
&\qquad
+
\mgb{\eta\left(\eta+\frac{(1+\eta L)^{2}}{\mu(1-\frac{\eta L}{2})}\right)
\mathsf H_k\left(1-\eta\mu\left(1-\frac{\eta L}{2}\right),\delta^{\frac{2}{B}}\right)}
\frac{3L^{2} \delta^{-2}}{N}\mathbb{E}\left\Vert x^{(0)}\right\Vert^{2}.
\end{align*}}
Under our assumption, the stepsize $\eta<2/L$, such that
$1-\eta\mu\left(1-\frac{\eta L}{2}\right)\in[0,1)$.
Hence, we conclude that for every $k$,
{\small
\begin{align*}
\mathbb{E}\left\Vert\bar{x}^{(k)}-x_{k}\right\Vert^{2}
&\leq\frac{\eta\left(\eta+\frac{(1+\eta L)^{2}}{\mu(1-\frac{\eta L}{2})}\right)
\left(\frac{3L^{2}D^{2}\eta^{2} \delta^{-2}}{N(1-\delta^{\frac{1}{B}})^{2}}
+\frac{6dL^{2}\eta\cdot  \delta^{-2}}{1-\delta^{\frac{2}{B}}}
\right)+\eta^{2}\frac{\sigma^{2}}{N}}{
1-\left(1-\eta\mu\left(1-\frac{\eta L}{2}\right)\right)}
\\
&\qquad
+
\mgb{\eta\left(\eta+\frac{(1+\eta L)^{2}}{\mu(1-\frac{\eta L}{2})}\right)
\mathsf H_k\left(1-\eta\mu\left(1-\frac{\eta L}{2}\right),\delta^{\frac{2}{B}}\right)}
\frac{3L^{2} \delta^{-2}}{N}
\mathbb{E}\left\Vert x^{(0)}\right\Vert^{2}
\\
&=\frac{\eta\left(\eta+\frac{(1+\eta L)^{2}}{\mu(1-\frac{\eta L}{2})}\right)
\left(\frac{3L^{2}D^{2}\eta  \delta^{-2}}{N(1-\delta^{\frac{1}{B}})^{2}}
+\frac{6dL^{2}\cdot  \delta^{-2}}{1-\delta^{\frac{2}{B}}}
\right)+\eta\frac{\sigma^{2}}{N}}{
\mu\left(1-\frac{\eta L}{2}\right)}\\
&\qquad+
\mgb{\eta\left(\eta+\frac{(1+\eta L)^{2}}{\mu(1-\frac{\eta L}{2})}\right)
\mathsf H_k\left(1-\eta\mu\left(1-\frac{\eta L}{2}\right),\delta^{\frac{2}{B}}\right)}\frac{3L^{2} \delta^{-2}}{N}
\mathbb{E}\left\Vert x^{(0)}\right\Vert^{2}.
\end{align*}}
The proof is complete.
\end{proof}
\section{Proof of Lemma \ref{lem:inexact}}
\begin{lemma}[Restatement of Lemma \ref{lem:inexact}]
Assume 
$\sqrt{1-\frac{\eta\mu\beta}{\beta+1}}\leq\lambda<1$
and $\eta\leq\frac{1}{(1+\alpha)L}$,
where $\alpha,\beta>0$.
Let the square-integrable random sequences $\{p^{(k)}\}_{k\geq0}$ and
$\{s_i^{(k)}\}_{k\geq0}$ satisfy \eqref{eq-p-recursion} almost surely where \(\nabla g_i\) is evaluated at \(s_i^{(k)}\) at iteration \(k\) for each \(i=1,\ldots,N\).
\mgb{Consider the weighted $L_2$ norm
\[
\|p-p_*\|_{L_2}^{\lambda,K}:=
\max_{0\leq k\leq K}\lambda^{-k}
\bigl(\mathbb E\|p^{(k)}-p_*\|^2\bigr)^{1/2},
\]
where $\|p-s_i\|_{L_2}^{\lambda,K}$ is defined analogously. Then, for every
$K\geq0$,
\begin{equation}\label{eq:stochastic-inexact-L2}
\|p-p_*\|_{L_2}^{\lambda,K}
\leq
\bigl(\mathbb E\|p^{(0)}-p_*\|^2\bigr)^{1/2}
+\frac{1}{\lambda\sqrt N}
\sqrt{\frac{L(1+\alpha)}{\mu\alpha}+\beta}
\sum_{i=1}^N\|p-s_i\|_{L_2}^{\lambda,K},
\end{equation}}
where $p_{\ast}$ is the minimizer of  $g:=\frac{1}{N}\sum_{i=1}^{N}g_{i}(x)$.
\end{lemma}
\mgb{\begin{proof} We start by noting that by $L$-smoothness,  $\|\nabla g_i(s_i^{(k)})\| \leq  \|\nabla g_i(0)\| + L \| s_i^{(k)}\| $, therefore square-integrability of $s_i^{(k)}$ implies the square-integrability of the gradients $\nabla g_i(s_i^{(k)})$.
Now, for each $k\geq 0$, define
\[
h^{(k)}
:=
\frac{1}{N}\sum_{i=1}^{N}\nabla g_i\!\left(s_i^{(k)}\right),
\qquad
\varepsilon_k
:=
\frac{1}{N}\sum_{i=1}^{N}
\left\|p^{(k)}-s_i^{(k)}\right\|^2.
\]
Then the recursion \eqref{eq-p-recursion} can be written as
\[
p^{(k+1)}=p^{(k)}-\eta h^{(k)}.
\]
We also set
\[
a^2
:=
1-\frac{\eta\mu\beta}{\beta+1},
\qquad
H^2
:=
\frac{L(1+\alpha)}{\mu\alpha}+\beta.
\]

We first derive a deterministic one-step inequality. This part of the
argument is pointwise, so no probabilistic assumptions or independence
conditions are needed.

Fix $k\geq0$ and, to simplify notation, write
\[
p:=p^{(k)},\qquad
p^+:=p^{(k+1)},\qquad
s_i:=s_i^{(k)},\qquad
h:=h^{(k)}.
\]
Since $p^+=p-\eta h$, we have
\[
p-p_*=p^+-p_*+\eta h,
\]
and hence
\begin{equation}\label{eq:L2-proof-distance-identity}
\|p^+-p_*\|^2
=
\|p-p_*\|^2
-2\eta\langle h,p^+-p_*\rangle
-\eta^2\|h\|^2.
\end{equation}

We next lower-bound the inner product in
\eqref{eq:L2-proof-distance-identity}. Since each $g_i$ is
$\mu$-strongly convex,
\[
g_i(p_*)
\geq
g_i(s_i)
+
\left\langle
\nabla g_i(s_i),p_*-s_i
\right\rangle
+
\frac{\mu}{2}\|p_*-s_i\|^2.
\]
Equivalently,
\begin{equation}\label{eq:L2-proof-sc}
\left\langle
\nabla g_i(s_i),s_i-p_*
\right\rangle
\geq
g_i(s_i)-g_i(p_*)
+
\frac{\mu}{2}\|s_i-p_*\|^2.
\end{equation}
On the other hand, $L$-smoothness of $g_i$ gives
\[
g_i(p^+)
\leq
g_i(s_i)
+
\left\langle
\nabla g_i(s_i),p^+-s_i
\right\rangle
+
\frac{L}{2}\|p^+-s_i\|^2,
\]
so that
\begin{equation}\label{eq:L2-proof-smooth}
\left\langle
\nabla g_i(s_i),p^+-s_i
\right\rangle
\geq
g_i(p^+)-g_i(s_i)
-
\frac{L}{2}\|p^+-s_i\|^2.
\end{equation}
Adding \eqref{eq:L2-proof-sc} and
\eqref{eq:L2-proof-smooth} yields
\[
\left\langle
\nabla g_i(s_i),p^+-p_*
\right\rangle
\geq
g_i(p^+)-g_i(p_*)
+
\frac{\mu}{2}\|s_i-p_*\|^2
-
\frac{L}{2}\|p^+-s_i\|^2.
\]
Averaging this inequality over $i=1,\ldots,N$ and recalling that
\[
g=\frac1N\sum_{i=1}^N g_i,
\qquad
h=\frac1N\sum_{i=1}^N\nabla g_i(s_i),
\]
we obtain
\begin{align}
\langle h,p^+-p_*\rangle
\geq{}&
g(p^+)-g(p_*)
+
\frac{\mu}{2N}\sum_{i=1}^N\|s_i-p_*\|^2
\nonumber\\
&\quad
-
\frac{L}{2N}\sum_{i=1}^N\|p^+-s_i\|^2.
\label{eq:L2-proof-inner}
\end{align}
Substituting \eqref{eq:L2-proof-inner} into
\eqref{eq:L2-proof-distance-identity} gives
\begin{align}
\|p^+-p_*\|^2
\leq{}&
\|p-p_*\|^2
-2\eta\bigl(g(p^+)-g(p_*)\bigr)
-\frac{\eta\mu}{N}
\sum_{i=1}^N\|s_i-p_*\|^2
\nonumber\\
&\quad
+\frac{\eta L}{N}
\sum_{i=1}^N\|p^+-s_i\|^2
-\eta^2\|h\|^2.
\label{eq:L2-proof-basic}
\end{align}

We now bound the last two terms in
\eqref{eq:L2-proof-basic}. Since
\[
p^+-s_i
=
-\eta h+(p-s_i),
\]
the elementary inequality
\[
\|u+v\|^2
\leq
(1+\alpha)\|u\|^2
+
\left(1+\frac1\alpha\right)\|v\|^2,
\qquad \alpha>0,
\]
gives
\[
\|p^+-s_i\|^2
\leq
(1+\alpha)\eta^2\|h\|^2
+
\left(1+\frac1\alpha\right)\|p-s_i\|^2.
\]
For completeness, this elementary inequality follows from
\[
2\langle u,v\rangle
\leq
\alpha\|u\|^2+\frac1\alpha\|v\|^2.
\]
Averaging over $i$ therefore gives
\begin{equation}\label{eq:L2-proof-forward}
\frac1N\sum_{i=1}^N\|p^+-s_i\|^2
\leq
(1+\alpha)\eta^2\|h\|^2
+
\frac{1+\alpha}{\alpha}\varepsilon_k.
\end{equation}
Consequently,
\begin{align}
&\frac{\eta L}{N}
\sum_{i=1}^N\|p^+-s_i\|^2
-\eta^2\|h\|^2
\nonumber\\
&\qquad\leq
\eta^2\bigl(\eta L(1+\alpha)-1\bigr)\|h\|^2
+
\eta L\frac{1+\alpha}{\alpha}\varepsilon_k.
\end{align}
By the stepsize assumption
\[
\eta\leq\frac{1}{(1+\alpha)L},
\]
we have $\eta L(1+\alpha)-1\leq0$. Hence
\begin{equation}\label{eq:L2-proof-forward-final}
\frac{\eta L}{N}
\sum_{i=1}^N\|p^+-s_i\|^2
-\eta^2\|h\|^2
\leq
\eta L\frac{1+\alpha}{\alpha}\varepsilon_k.
\end{equation}

It remains to control the strong-convexity term in
\eqref{eq:L2-proof-basic}. For every $\beta>0$ and every
$u,v\in\mathbb R^d$,
\begin{equation}\label{eq:L2-proof-reverse-young}
\|u-v\|^2
\geq
\frac{\beta}{\beta+1}\|u\|^2
-\beta\|v\|^2.
\end{equation}
Indeed,
\begin{align*}
\|u-v\|^2
-\frac{\beta}{\beta+1}\|u\|^2
+\beta\|v\|^2
&=
\frac{1}{\beta+1}\|u\|^2
-2\langle u,v\rangle
+(\beta+1)\|v\|^2\\
&=
\left\|
\frac{u}{\sqrt{\beta+1}}
-\sqrt{\beta+1}\,v
\right\|^2
\geq0.
\end{align*}
Apply \eqref{eq:L2-proof-reverse-young} with
\[
u=p-p_*,
\qquad
v=p-s_i.
\]
Since
\[
s_i-p_*=(p-p_*)-(p-s_i),
\]
we obtain
\[
\|s_i-p_*\|^2
\geq
\frac{\beta}{\beta+1}\|p-p_*\|^2
-\beta\|p-s_i\|^2.
\]
Averaging over $i$ yields
\begin{equation}\label{eq:L2-proof-backward}
\frac1N\sum_{i=1}^N\|s_i-p_*\|^2
\geq
\frac{\beta}{\beta+1}\|p-p_*\|^2
-\beta\varepsilon_k.
\end{equation}

Substituting \eqref{eq:L2-proof-forward-final} and
\eqref{eq:L2-proof-backward} into
\eqref{eq:L2-proof-basic}, we obtain
\begin{align}
\|p^+-p_*\|^2
\leq{}&
\left(
1-\frac{\eta\mu\beta}{\beta+1}
\right)\|p-p_*\|^2
-2\eta\bigl(g(p^+)-g(p_*)\bigr)
\nonumber\\
&\quad
+
\eta
\left(
\mu\beta+
\frac{L(1+\alpha)}{\alpha}
\right)
\varepsilon_k.
\end{align}
By the definitions of $a$ and $H$, this is exactly
\begin{equation}\label{eq:L2-proof-one-step}
\|p^{(k+1)}-p_*\|^2
\leq
a^2\|p^{(k)}-p_*\|^2
-2\eta
\bigl(g(p^{(k+1)})-g(p_*)\bigr)
+\eta\mu H^2\varepsilon_k.
\end{equation}

We now pass from the pointwise estimate to an $L_2$ estimate. Since each
$g_i$ is $\mu$-strongly convex, their average $g$ is also
$\mu$-strongly convex. Moreover, since $p_*$ minimizes the differentiable
function $g$, $\nabla g(p_*)=0$. Therefore
\[
g(z)-g(p_*)
\geq
\frac{\mu}{2}\|z-p_*\|^2,
\qquad z\in\mathbb R^d.
\]
Applying this with $z=p^{(k+1)}$ in
\eqref{eq:L2-proof-one-step} gives
\[
(1+\eta\mu)\|p^{(k+1)}-p_*\|^2
\leq
a^2\|p^{(k)}-p_*\|^2
+\eta\mu H^2\varepsilon_k.
\]
Taking expectations, define
\[
u_k
:=
\mathbb E\|p^{(k)}-p_*\|^2,
\qquad
v_k
:=
\mathbb E\varepsilon_k,
\]
and
\[
\rho^2
:=
\frac{a^2}{1+\eta\mu},
\qquad
c_0
:=
\frac{\eta\mu H^2}{1+\eta\mu}.
\]
Then
\begin{equation}\label{eq:L2-proof-scalar-recursion}
u_{k+1}
\leq
\rho^2u_k+c_0v_k.
\end{equation}
Notice that no independence assumption has been used: the preceding
inequality follows from a pointwise inequality followed by expectation.

Repeated substitution in
\eqref{eq:L2-proof-scalar-recursion} gives, for every $k\geq1$,
\begin{equation}\label{eq:L2-proof-unroll}
u_k
\leq
\rho^{2k}u_0
+
c_0\sum_{t=0}^{k-1}
\rho^{2(k-1-t)}v_t.
\end{equation}
The assumption
\[
a
=
\sqrt{1-\frac{\eta\mu\beta}{\beta+1}}
\leq\lambda
\]
implies
\[
\rho^2
=
\frac{a^2}{1+\eta\mu}
<
\lambda^2.
\]
For a fixed $k\in [1,K]$, let
\[
M_k
:=
\max_{0\leq t<k}
\lambda^{-2t}v_t.
\]
Since $t<k\le K$, every index $t$ entering $M_k$ is smaller than $K$. Then $v_t\leq\lambda^{2t}M_k$ for $0\leq t<k$, and
\eqref{eq:L2-proof-unroll} yields
\begin{align*}
\lambda^{-2k}u_k
&\leq
\left(\frac{\rho^2}{\lambda^2}\right)^k u_0
+
c_0M_k
\sum_{t=0}^{k-1}
\lambda^{-2k}
\rho^{2(k-1-t)}
\lambda^{2t}\\
&=
\left(\frac{\rho^2}{\lambda^2}\right)^k u_0
+
\frac{c_0}{\lambda^2}M_k
\sum_{j=0}^{k-1}
\left(\frac{\rho^2}{\lambda^2}\right)^j\\
&\leq
u_0
+
\frac{c_0}{\lambda^2-\rho^2}M_k.
\end{align*}
Furthermore,
\begin{align}
\frac{c_0}{\lambda^2-\rho^2}
&=
\frac{\eta\mu H^2}
{(1+\eta\mu)\lambda^2-a^2}.
\label{eq:L2-proof-coeff}
\end{align}
Since $a^2\leq\lambda^2$,
\[
(1+\eta\mu)\lambda^2-a^2
\geq
\eta\mu\lambda^2,
\]
and therefore
\begin{equation}\label{eq:L2-proof-coeff-final}
\frac{c_0}{\lambda^2-\rho^2}
\leq
\frac{H^2}{\lambda^2}.
\end{equation}
Combining the preceding inequalities gives
\begin{equation}\label{eq:L2-proof-weighted-uk}
\lambda^{-2k}u_k
\leq
u_0+\frac{H^2}{\lambda^2}M_k.
\end{equation}

It remains only to express $M_k$ in terms of the weighted $L_2$ norms
appearing in the statement. Define
\[
R_i
:=
\|p-s_i\|_{L_2}^{\lambda,K}
=
\max_{0\leq t\leq K}
\lambda^{-t}
\left(
\mathbb E\|p^{(t)}-s_i^{(t)}\|^2
\right)^{1/2}.
\]
For every $0\leq t\leq K$,
\begin{align*}
\lambda^{-2t}v_t
&=
\frac1N\sum_{i=1}^N
\lambda^{-2t}
\mathbb E\|p^{(t)}-s_i^{(t)}\|^2\\
&\leq
\frac1N\sum_{i=1}^N R_i^2\\
&\leq
\frac1N
\left(\sum_{i=1}^N R_i\right)^2.
\end{align*}
Hence for $0 < k \leq K$,
\begin{equation}\label{eq:L2-proof-input}
M_k^{1/2}
\leq
\frac1{\sqrt N}
\sum_{i=1}^N
\|p-s_i\|_{L_2}^{\lambda,K}.
\end{equation}
Taking square roots in
\eqref{eq:L2-proof-weighted-uk} and using
$\sqrt{r+s}\leq\sqrt r+\sqrt s$ for $r,s\geq0$ gives
\begin{align*}
\lambda^{-k}
\left(
\mathbb E\|p^{(k)}-p_*\|^2
\right)^{1/2}
&\leq
\left(
\mathbb E\|p^{(0)}-p_*\|^2
\right)^{1/2}
+
\frac{H}{\lambda}M_k^{1/2}\\
&\leq
\left(
\mathbb E\|p^{(0)}-p_*\|^2
\right)^{1/2}
+
\frac{H}{\lambda\sqrt N}
\sum_{i=1}^N
\|p-s_i\|_{L_2}^{\lambda,K}.
\end{align*}
The same inequality is trivially valid for $k=0$. Maximizing over
$0\leq k\leq K$ therefore yields the stronger estimate
\[
\|p-p_*\|_{L_2}^{\lambda,K}
\leq
\left(
\mathbb E\|p^{(0)}-p_*\|^2
\right)^{1/2}
+
\frac{1}{\lambda\sqrt N}
\sqrt{
\frac{L(1+\alpha)}{\mu\alpha}+\beta
}
\sum_{i=1}^N
\|p-s_i\|_{L_2}^{\lambda,K}.
\]
Then, the stated bound \eqref{eq:stochastic-inexact-L2} follows.
\end{proof}}

\end{document}